\documentclass[a4paper, 11pt]{amsart}
\usepackage{amssymb,graphicx}
\usepackage{graphicx} 
\usepackage[utf8]{inputenc}
\usepackage[T1]{fontenc}
\usepackage{float}
\usepackage{amsmath, amssymb, amscd, amsthm,epsfig,epic,eepic}
\usepackage[english]{babel}
\usepackage{mathrsfs}
\usepackage{tikz-cd}
\usepackage{quiver}
\usepackage{bbm}
\usepackage{verbatim}
\usepackage{color}
\usepackage[all]{xy}
\usepackage{tikz-cd}
\usepackage{enumerate}
\usepackage[justify]{ragged2e}
\usepackage[normalem]{ulem}
\usepackage[dvipsnames]{xcolor}
\usepackage{hyperref,cleveref}
\hypersetup{
    colorlinks=true, linkcolor=RedViolet, citecolor=NavyBlue, filecolor=teal, urlcolor=RedViolet,
    linktocpage=true}
\usepackage{mathbbol}
\usepackage{bbm}
\definecolor{green}{rgb}{0,0.6,0} 

\theoremstyle{plain}

\newtheorem{theorem}{Theorem}[section]
\newtheorem{lemma}[theorem]{Lemma}
\newtheorem{question}[theorem]{Question}
\newtheorem{proposition}[theorem]{Proposition}
\newtheorem{corollary}[theorem]{Corollary}
\newtheorem{defn}[theorem]{Definition}
\newtheorem{theoremA}{Theorem}

\newtheorem{propositionA}[theoremA]{Proposition}

\usepackage[left=3cm,right=3cm,top=3cm,bottom=3cm]{geometry}

\makeatletter
     \let\oldfootnote\footnote
     \def\footnote{\@ifstar\footnote@star\footnote@nostar}
     \def\footnote@star#1{{\let\thefootnote\relax\footnotetext{#1}}}
     \def\footnote@nostar{\oldfootnote}
     \makeatother

\theoremstyle{definition}
\newtheorem{definition}[theorem]{Definition}
\newtheorem{remark}[theorem]{Remark}
\newtheorem{remarks}[theorem]{Remarks}
\newtheorem{example}[theorem]{Example}

\newcommand{\et}{\text{\'et}}
\newcommand{\id}{{\rm id}}
\newcommand{\pr}{{\rm pr}}
\newcommand{\red}{{\rm red}}
\newcommand{\sic}{{\rm sc}}
\newcommand{\sep}{{\rm sep}}
\newcommand{\sm}{{\rm sm}}

\newcommand{\Aut}{\operatorname{Aut}}
\newcommand{\Br}{\operatorname{Br}}
\newcommand{\Cl}{\operatorname{Cl}}
\newcommand{\Eff}{\operatorname{Eff}}
\newcommand{\Gal}{\operatorname{Gal}}

\newcommand{\Grass}{\operatorname{Grass}}
\newcommand{\Int}{\operatorname{Int}}
\newcommand{\Ker}{\operatorname{Ker}}
\newcommand{\Lie}{\operatorname{Lie}}
\newcommand{\Pic}{\operatorname{Pic}}
\newcommand{\Spec}{\operatorname{Spec}}
\newcommand{\Stab}{\operatorname{Stab}}
\newcommand{\Supp}{\operatorname{Supp}}
\newcommand{\Sym}{\operatorname{Sym}}

\newcommand{\Or}{\mathbf{O}}
\newcommand{\GL}{\mathbf{GL}}
\newcommand{\PGL}{\mathbf{PGL}}
\newcommand{\SL}{\mathbf{SL}}
\newcommand{\SU}{\mathbf{SU}}
\newcommand{\Sp}{\mathbf{Sp}}
\newcommand{\PSp}{\mathbf{PSp}}
\newcommand{\bG}{\mathbf{G}}

\newcommand{\bF}{\mathbb{F}}
\newcommand{\bN}{\mathbb{N}}
\newcommand{\bP}{\mathbb{P}}
\newcommand{\bQ}{\mathbb{Q}}
\newcommand{\bZ}{\mathbb{Z}}

\newcommand{\bfAut}{\mathbf{Aut}}
\newcommand{\bfalpha}{\boldsymbol{\alpha}}

\newcommand{\cD}{\mathcal{D}}
\newcommand{\cG}{\mathcal{G}}
\newcommand{\cH}{\mathcal{H}}
\newcommand{\cL}{\mathcal{L}}
\newcommand{\cN}{\mathcal{N}}
\newcommand{\cO}{\mathcal{O}}

\newcommand{\cT}{\mathcal{T}}
\newcommand{\cX}{\mathcal{X}}
\newcommand{\cU}{\mathcal{U}}

\newcommand{\fg}{\mathfrak{g}}
\newcommand{\fm}{\mathfrak{m}}
\newcommand{\fp}{\mathfrak{p}}

\newcommand{\rN}{{\rm N}}
\newcommand{\rR}{{\rm R}}

\newcommand{\uAut}{\underline{\Aut}}

\tikzset{
  root/.style={circle, draw, thick, minimum size=5pt, inner sep=0pt},
  dots/.style={draw=none, fill=none}
}
\title[Complete homogeneous varieties over a field]{Complete homogeneous varieties\\ 
over a field}

\author{Michel Brion}
\address{Institut Fourier, UMR CNRS 5582, Universit\'e de Grenoble, France}
\email{Michel.Brion@univ-grenoble-alpes.fr}

\author{Matilde Maccan}
   \address{Institut Camille Jordan, UMR CNRS 5208, Universit\'e Lyon 1, Villeurbanne, France}
\email{maccan@math.univ-lyon1.fr} 

\author{Srinivasan Srimathy}
 \address{School of Mathematics\\
Tata Institute of Fundamental Research\\
Mumbai\\
India}
\email{srimathy@math.tifr.res.in}
\date{\today}

\begin{document}
\begin{abstract}
The main objects of this paper are complete homogeneous varieties under the action of a semisimple group $G$ over field $k$ of  characteristic $p>0$. 
When the isotropy groups of the geometric points are smooth parabolics, these are classically known as twisted flag varieties and are well studied in the literature.  If the smoothness assumption is removed, we get a larger class of varieties  which we call as \emph{twisted parabolic varieties}.  The main goal of this paper is to study the structure, classification and properties  of these varieties.  We first understand their structure over a separable closure of $k$ by classifying non-smooth parabolics containing a Borel; this builds on the recently completed classification
of such parabolics over an algebraically closed field.
Then, via  Galois descent we get a  combinatorial criterion for such a variety 
to descend over $k$. We compute the Picard group of twisted parabolic varieties
by embedding them  into products of twisted flag varieties, study their local rigidity properties and show that they are rational if they have a $k$-rational point 
for  $p \geq 3$.  
When $k$ is perfect, we show many interesting properties such as existence of a universal purely inseparable cover, transitivity of $k$-rational points 
under $G(k)$ as well as the existence of a smooth parabolic 
in the stabilizer of every such point. 
However, these properties fail if $k$ is imperfect and we demonstrate this 
 with illuminating examples. 
\end{abstract}
\maketitle

\footnote*{\scriptsize{Keywords: homogeneous variety, group scheme, positive characteristic, (non-smooth) parabolic subgroup.}}

\footnote*{\scriptsize{2020 Mathematics Subject Classification: 14M15, 14M17, 14L15, 14L30, 14G17, 17B22}}

\addtocontents{toc}{\setcounter{tocdepth}{1}}
\tableofcontents

\section{Introduction}

The focus of the present work is the structure of complete algebraic 
varieties $X$ that are homogeneous under a smooth connected algebraic 
group $G$, over an arbitrary ground field $k$ of characteristic $p > 0$. 
By \emph{homogeneous}, we mean that $G$ acts on $X$ and the induced 
action of $G(\bar{k})$ on $X(\bar{k})$ is transitive, where $\bar{k}$ 
denotes an algebraic closure of $k$. Examples include abelian varieties
and their torsors, projective spaces $\bP^n$ and their forms (Severi--Brauer 
varieties), quadrics, as well as varieties that are special to positive 
characteristics and whose geometry presents new phenomena.  

\smallskip

For instance, the hypersurface
\begin{equation}\label{eqn:tiv}
X_m := ( x_0^{p^m} y_0 + \cdots + x_n^{p^m} y_n = 0 )
\subset \bP^n \times \bP^n 
\end{equation}
is homogeneous under the special linear group 
$\SL_{n+1}$ acting on $x_0,\ldots,x_n$ via matrix
multiplication, and on $y_0,\ldots,y_n$ via multiplication by
the inverse transpose of the matrix of $p^m$-th powers. 
The two projections $\pi_1,\pi_2 : X_m \to \bP^n$ are equivariant
and the pullbacks $\pi_1^* \cO(1), \pi_2^* \cO(1)$ yield a basis 
of $\Pic(X_m)$. If $m = 0$ then both projections are $\bP^{n-1}$-bundles.
But only $\pi_1$ is a $\bP^{n-1}$-bundle for $m \geq 1$: then 
$\pi_2$ is everywhere non-smooth. The geometry of $X_m$ for 
$m \geq 1$ differs from that of the incidence variety $X_0$
in many other aspects, e.g., the tangent bundle of $X_m$
is globally generated if and only if $m = 0$, and $X_m$ is Fano
if and only if $n \geq p^m$. The main reason is that the stabilizer 
of $X_m$ is a non-smooth parabolic subgroup of $\SL_{n+1}$ for $m > 0$ 
(it will be described explicitly below, see Example \ref{ex infinite non isomorphic}). 
Similar constructions 
of such varieties yield counterexamples to the Kodaira vanishing 
theorem for ample line bundles, see \cite{LauritzenRao, Totaro}. 

\smallskip

The classification of complete homogeneous varieties over 
an \emph{algebraically closed} field $k$ is now well understood,
with the small characteristics $p = 2,3$
being brought to an end very recently.
A general structure result asserts that every such variety
admits a decomposition $X \simeq A \times Y$, where $A$ is an abelian 
variety and $Y$ is a complete homogeneous variety under a 
\emph{semisimple} group $H$ (this originates in the work of 
Borel and Remmert on homogeneous compact K\"ahler manifolds,
see \cite{BR}. The present version is due to Sancho de Salas
\cite[Theorem 5.2]{SS}). 
Choosing a point $y \in Y(k)$ and denoting by $P$ its stabilizer, 
we may identify $Y$ with the homogeneous space $H/P$, where 
$P \subset H$ is a \emph{parabolic subgroup}. If $P$ is smooth 
then $Y$ is called a \emph{flag variety}, the basic examples 
being Grassmannians and varieties of complete flags inside 
of some vector space. Smooth parabolics and flag varieties 
are classical objects whose structure and classification 
are well known, see e.g.~\cite{BT}. 
The study of non-smooth parabolics started much more recently; 
these have been classified by Wenzel if $p \geq 5$ (see \cite{Wenzel}), 
and by the second author (\cite{Maccan1, Maccan2}) for $p = 2,3$. 
Wenzel's result (recalled below) implies that most parabolics are 
non-smooth. Equivalently, flag varieties form a restricted class 
of complete homogeneous varieties under semisimple groups.

\smallskip 

To the authors' knowledge, the only previous work dealing with complete
homogeneous varieties with non-smooth stabilizers over a 
\emph{non-algebraically closed} 
field is that of \cite{sri_motives}, where the focus is on their Chow motives 
with finite coefficients. The present article is instead concerned with 
the structure and classification of such varieties. The results presented 
here should be regarded as a first contribution towards their study, 
especially over imperfect fields. Still, many questions remain open. 
Some of these are collected at the end of the present introduction, 
indicating that much remains to be understood about complete homogeneous 
varieties and parabolic subgroups of semisimple groups in this generality.

\smallskip

We now give an overview of our main results. We work over a ground
field $k$, with algebraic closure $\bar{k}$ and separable
closure $k^{\sep} \subset \bar{k}$. We first consider 
(not necessarily complete) homogeneous varieties under 
a smooth connected algebraic group $G$. Such a variety $X$ 
may well have no $k$-rational point, but it has a 
$k^{\sep}$-point $x$. The base change $X^{\sep} := X_{k^{\sep}}$ 
is then identified with the homogeneous space $G^{\sep}/H$,
where $H$ denotes the stabilizer of $x$.   

\smallskip

The homogeneous varieties most commonly studied in the literature 
are those whose stabilizer is smooth. We refer to them as \emph{hss} 
varieties (shorthand for ``homogeneous with smooth stabilizer'') 
and provide several equivalent characterizations of this 
property. In particular, a smooth $G$-variety is hss if and only if 
its tangent bundle is generated by the global sections induced by 
the Lie algebra of $G$.

\smallskip

When the ground field $k$ is \emph{perfect}, we show the existence of 
a \emph{universal purely inseparable cover} of a homogeneous $G$-variety $X$, 
denoted by $\hat{X}$. That is, we show  that there exists a purely inseparable 
$G$-equivariant morphism $f \colon \hat{X} \longrightarrow X$, which is
universal for these properties (Proposition \ref{prop:hss}). When there
is a rational point $x \in X(k)$, this morphism coincides with 
the canonical map $G/H_{\red} \rightarrow G/H$, where $H$ is the stabilizer 
of $x$ and $H_{\red}$ its reduced subscheme (which is the largest smooth 
subgroup of $H$ in this setting). A particularly useful property is 
that $X$ has a $k$-rational point if and only if so does $\hat{X}$.

\medskip

We next prove the following generalization of Sancho de Salas' result
mentioned above: any complete homogeneous variety $X$ over an
arbitrary field $k$ admits a decomposition 
$X \simeq Y \times Z$, where $Y$ is a torsor under an abelian variety 
and $Z$ is homogeneous under a semisimple group. 
This decomposition is canonical and equivariant. Moreover, it yields 
the following description at the level of Picard groups: 
$\Pic(Z)$ is a free abelian group of finite rank, and the pull-back 
map induces an isomorphism over the rationals,
\[ \Pic(Y)_{\bQ} \times \Pic(Z)_{\bQ} 
\stackrel{\sim}{\longrightarrow} 
\Pic(X)_{\bQ}. \]
In addition, the effective cone of $Z$ is generated by a basis of 
$\Pic(Z)_{\bQ}$ and every effective divisor on $Z$ is base-point-free.
In particular, $\Pic(Z)_{\bQ}$ is a finite-dimensional vector
space equipped with a canonical basis.
If $k$ is separably closed then these results hold over the integers, 
and the Schubert divisors form the canonical basis of $\Pic(Z)$.
Over an arbitrary field $k$, the basis of $\Pic(Z)_{\bQ}$ is 
obtained from that of $\Pic(Z^{\sep})$ by averaging over the Galois
group. It is important here that we work with rational 
coefficients: for example, if $Z$ is a conic without rational
point then $\Pic(Z)$ is a subgroup of index $2$ of $\Pic(Z^\sep)$.

\smallskip

The above direct product decomposition reduces the study of complete 
homogeneous varieties to two fundamentally distinct families. 
In the present work, we focus on one of them. The theory of torsors 
under abelian varieties is already well developed in the literature, 
beginning with \cite{LangTate}. So from now on, we consider 
complete homogeneous varieties under a \emph{semisimple} group $G$.

\smallskip

We now recall and adopt the following definition from \cite{SS}.

\begin{definition}\label{def:parabolic}
Let $G$ be a semisimple group over a field $k$.
A \emph{parabolic subgroup} of $G$ is a subgroup $P \subset G$ 
such that the homogeneous space $G/P$ is complete. 
A \emph{parabolic variety} is a complete $G$-variety of the form $G/P$. 
The latter is called a \emph{flag variety} if $P$ is smooth.
\end{definition}

In other words, the parabolic varieties are exactly the complete 
homogeneous $G$-varieties that have a $k$-rational point. 

\medskip

\noindent\textbf{Over a separably closed field.} 
We first consider the case where $k = k^{\sep}$, so that 
$X \simeq G/P$ as above. Moreover, the semisimple group $G$ 
is split, and hence admits a split Borel subgroup $B$. Then
$B$ fixes a unique $k$-rational point of $X$ by Proposition 
\ref{prop:rational}. Thus, we may assume that $P$ contains $B$ 
(but over an imperfect field, some parabolics contain no Borel 
subgroup as the examples discussed later will illustrate).  
By \cite{BT}, there exist finitely many maximal smooth parabolic
subgroups $P_1,\ldots,P_n$ containing $B$; they are indexed by 
the set of simple roots, and $n$ is the rank of $G$. Moreover, 
any \emph{smooth} parabolic containing $B$ is of the form
\[
P = P_{i_1} \cap \ldots \cap P_{i_r}
\]
where $i_1,\ldots, i_r$ are unique up to reordering.

For non-smooth parabolics, we obtain the following structure result,
first proved by Wenzel \cite[Theorem 10]{Wenzel} over an algebraically 
closed field:

\begin{theoremA}[Theorem \ref{thm:wenzel}]
\label{thm:wenzelintro}
Assume that $k$ is separably closed and let $P$ be a parabolic subgroup 
of $G$ containing $B$. Denote by $Q \supset B$ the largest smooth subgroup 
of $P$. Then the multiplication map of $G$ induces an isomorphism of schemes
\[
\prod_\gamma (U_{-\gamma} \cap P) \times Q \stackrel{\sim}{\longrightarrow} P,
\]
where $U_{-\gamma}$ varies among the root subgroups not contained in $Q$. 
For any such $\gamma$, there exists a unique non-negative integer 
$m = m(\gamma)$ such that 
$U_{- \gamma} \cap P \simeq \bfalpha_{p^m}$ (the $m$-th Frobenius kernel of $\bG_a$). In particular we have $Q = P_\red$.
\end{theoremA}

The key point is that, unlike over the algebraic closure, we cannot 
assume a priori that the largest smooth subgroup of $P$ coincides 
with its reduced subscheme; this is instead part of the statement.

If $p \geq 5$, this takes the following form: every parabolic 
containing $B$ writes uniquely as 
\begin{align}
\label{eqn:wenzelintro}
P = {}_{m_1}G \cdot P_{i_1} \cap \ldots \cap {}_{m_r}G \cdot P_{i_r},
\end{align}
where ${}_m G$ denotes the kernel of the relative Frobenius 
homomorphism $F^m : G \to G^{(m)}$ for any $m \geq 0$, the
$m_i$s are non-negative integers, and the pairs 
$(m_1,i_1),\ldots, (m_r,i_r)$ are unique up to reordering. 
In other words, we can write $P$ as an intersection of smooth parabolics, 
each fattened by some Frobenius kernel. 

\smallskip

Following \cite{Wenzel}, we associate to each $P \supset B$ 
its \emph{Wenzel function}, namely the numerical function 
$\varphi$ that assigns to each positive root $\gamma$ the integer 
$m=m(\gamma)$ defined above. Then $P$, and hence $X$, is uniquely 
determined by $\varphi$; we write $P = P_{\varphi}$ and $X = X_{\varphi}$.
If one restricts to characteristic $p \geq 5$, the Wenzel function of $P$ 
is in turn completely determined by the integers $m_i$ in 
(\ref{eqn:wenzelintro}). Nevertheless, it is often convenient, and 
in some cases necessary, to work with Wenzel functions in order 
to obtain uniform statements as new phenomena arise in characteristics 
$2$ and $3$ (see \cite{Maccan1, Maccan2}).

\smallskip

The Wenzel function, assuming that  the acting group is fixed
and the action is faithful, can actually be recovered from 
the isomorphism class of $X$ as a variety. More precisely, 
we prove the following (Proposition \ref{prop:characterisationvarphi}). 
For a semisimple group $G$ over $k = k^\sep$, let $X$ be a parabolic 
variety on which $G$ acts faithfully. Then, up to automorphisms 
of the Dynkin diagram of $G$, there is a unique Wenzel function 
$\varphi$ such that $X = G/P_\varphi$ as varieties.

\medskip

\noindent\textbf{Parabolic varieties.}
The next step is to classify the parabolic varieties, 
now over an arbitrary ground field $k$. As expected, this requires 
a Galois descent argument. If $p \geq 5$, we obtain the following 
description: for a parabolic $k^\sep$-subgroup $P$ of the form 
(\ref{eqn:wenzelintro}), the homogeneous space $G^\sep/P$ is defined 
over $k$ if and only if both the set $\{P_{i_1},\ldots,P_{i_r}\}$ 
and the integers $m_1,\ldots,m_r$ are invariant under the Galois action. 
To obtain a characteristic-free statement, we use 
the Wenzel function; as one might expect, the corresponding criterion 
is the following (Proposition \ref{prop:definedfunction}). 
Let $P$ be a parabolic $k^\sep$-subgroup containing a Borel subgroup. 
Then the corresponding homogeneous space is defined over $k$ 
if and only if the associated Wenzel function $\varphi$ is invariant 
under the Galois action.

\smallskip

If the field is \emph{perfect}, we also get the following 
analogue of the big cell for flag varieties:

 \begin{propositionA}[Proposition \ref{prop:perfect}]
 \label{prop:rational intro}
Let $X$ be a parabolic $G$-variety over a perfect field $k$. 
Let $P_{\min} \subset G$ be a minimal smooth parabolic,
and $U^-$ the unipotent radical of an opposite parabolic. 
Then 

\begin{enumerate}
\item[{\rm (i)}]
$X(k)$ is the unique orbit of $G(k)$.
\item[{\rm (ii)}]
The stabilizer of every $k$-rational point of $X$ contains
a smooth parabolic subgroup. Moreover, $X$ has a unique 
$k$-rational $P_{\min}$-fixed point $x$. 
\item[{\rm (iii)}]
The orbit $V := U^- \cdot x$ is open 
in $X$ and isomorphic to an affine space. 
\item[{\rm (iv)}]
$X$ is covered by the translates $g \cdot V$,
where $g \in G(k)$.
\end{enumerate}
\end{propositionA}

This is due to Borel and Tits \cite{BT} in the setting of flag varieties
over an \emph{arbitrary} field $k$. Our statement follows easily from this
in view of our earlier result on the universal purely inseparable covers, 
which uses the perfectness assumption. Let us also mention that, 
if $k$ is instead assumed to be \emph{separably closed}, then 
part of the statement also holds: there is a unique $k$-rational point 
fixed by the Borel subgroup $B$, and the corresponding $R_u(B^-)$-orbit 
is open and isomorphic to an affine space (Proposition \ref{prop:rational}). 

\smallskip

The situation is very different over an \emph{arbitrary} field $k$:
in Section \ref{sec:example}, we will construct examples of parabolic varieties 
under a faithful action of a semisimple group, such that the stabilizer 
of every $k$-rational point contains no smooth parabolic. Still, all parabolic 
varieties turn out to be rational in a strong sense:

\begin{theoremA}[Theorem \ref{thm:krational}]
Every parabolic variety over a field of characteristic 
$p \geq 3$ has an open covering by affine spaces.
\end{theoremA}

\medskip

\noindent\textbf{Forms of parabolic varieties.} 
At this stage, what is left is to consider complete homogeneous 
varieties over an arbitrary field $k$, without assuming that they admit 
a $k$-rational point.

\begin{definition}\label{def:twistedparabolic}
Let $G$ be a semisimple group over a field $k$. 
A homogeneous $G$-variety $X$ is called a 
\emph{twisted parabolic variety} (resp.~a \emph{twisted flag variety}) 
if $X^\sep$ is a parabolic variety (resp.~a flag variety). 
\end{definition}

In other words, we add the adjective twisted to indicate that we are looking 
at \emph{forms}: varieties that become isomorphic to a given one 
(in this case, quotients of a semisimple group by a parabolic subgroup) 
after extending scalars to the separable closure. Let us mention that, 
in the article \cite{sri_motives}, twisted parabolic varieties are referred 
by the name \emph{projective pseudo-homogeneous varieties}. We prefer 
to use the term twisted parabolic varieties, as it suggests a natural 
generalization of twisted flag varieties.

 \smallskip
 
 We now turn to another of the main results of the present text, 
 namely the following structure theorem. Let $X$ be a twisted parabolic 
 $G$-variety and let $\pi_i \colon X \rightarrow X_i^\prime$ 
 be the contractions associated to the canonical basis 
 $(D_1,\ldots,D_s)$ of $\Pic(X)_\bQ$. If $X = G/P$ where $P$ is 
 of the form (\ref{eqn:wenzelintro}), these contractions are just 
 the morphisms $G/P \to G/{}_{m_i}G \cdot P_i$ corresponding to 
 the inclusions $P \subset {}_{m_i}G \cdot P_i$. In general, 
 the maps $\pi_i$ are still well-defined, because the contraction 
 associated with a base-point-free divisor is unchanged when 
 replacing it with a positive multiple.

\begin{theoremA}[Theorem \ref{thm:structure}]
Let $k$ be an arbitrary field of characteristic $p \geq 5$. 
With the above notation, the following hold.
\begin{enumerate}
\item[{\rm (i)}] 
Each $X_i^\prime$ is a twisted parabolic $G$-variety
for a unique action of $G$ such that 
$\pi_i$ is equivariant.
\item[{\rm (ii)}]
The morphism 
\[ (\pi_1,\ldots,\pi_s) : X \longrightarrow 
X_1^\prime \times \cdots \times X_s^\prime \]
is a $G$-equivariant immersion with image
being the smallest closed $G$-subvariety.
\item[{\rm (iii)}]
The $k$-group $G$ acts on each $X_i^\prime$ via a unique
quotient $G^{(m_i)} = G/{_{m_i}} G$ such that
the homogeneous $G^{(m_i)}$-variety $X_i^\prime$ has
a smooth stabilizer.
\item[{\rm (iv)}] 
We have $\Pic(X_i^\prime) \simeq \bZ$ and the pull-back
of the ample generator of $\Pic(X_i^\prime)$ is a
positive multiple of $D_i$.
\end{enumerate}
\end{theoremA}

Now as a consequence, we get a description over an arbitrary 
field $k$ of \emph{all} parabolic subgroups, meaning that 
the assumption that $P$ contains a smooth parabolic is removed. 
Specifically, if $p \geq 5$ any parabolic subgroup $P \subset G$ 
writes uniquely as 
\begin{align}
\label{equation:structure}
P = \bigcap_{i = 1}^s (F^{m_i})^{-1}(Q_i)
\qquad \text{and} \qquad
\Pic(G^{(m_i)}/Q_i) \simeq \bZ
\end{align}
where the $m_i$ are non-negative integers,
and $Q_i \subset G^{(m_i)}$ are smooth parabolics. 
If $k$ is perfect, we recover that 
$(F^{m_i})^{-1}(Q_i) = {_{m_i}}G \cdot P_i  $ 
for unique maximal smooth parabolic subgroups 
$P_i \subset G$; then $Q_i = P_i^{(m_i)}$.

\medskip

Let us emphasize that the introduction, as well as the main structure results in it, are presented under the assumption that $p \geq 5$. This simplifies both the statements and the exposition, as it allows us to work with pullbacks of smooth parabolics and Frobenius kernels. Whenever possible, however, we include a separate subsection treating the cases of small characteristics.

\smallskip

\noindent {\bf Local rigidity.} Flag varieties are locally rigid: 
if a smooth proper family contains a geometric fiber that is 
a flag variety, then there exists an open neighbourhood of 
the corresponding point in the base over which every geometric fiber 
is a flag variety of the same type \cite{Demazure}. This property 
does not extend to parabolic varieties. Indeed, one may check that 
it fails for the varieties $X_m$ defined in (\ref{eqn:tiv}), as soon as 
$m \geq 1$.

\smallskip

On the other hand, parabolic varieties are locally rigid \emph{as $G$-varieties} over an algebraically closed field \cite[Theorems 2.2 and 3.3]{SS}. We obtain the following generalization to parabolic varieties 
over an arbitrary field $k$ (Theorem \ref{thm:rigid}). Let 
$f \colon \mathcal{X} \rightarrow \mathcal S$ be
a proper map of schemes, where $G$ is split semisimple over $k$ and acts on 
$\mathcal X$. Assume moreover that the fiber above every $s \in \mathcal S$ 
is homogeneous under $G_{\kappa(s)}$. If the base $\mathcal S$ is reduced, 
connected and locally noetherian and $f$ has a section, then there exists 
a unique parabolic subgroup $P \subset G$ containing $B$ and a unique 
$G$-equivariant isomorphism of $\mathcal S$-schemes 
\[ F : G/P \times \mathcal S \stackrel{\sim}{\longrightarrow} \cX. \]

 A second rigidity result (Theorem \ref{thm:closed}) establishes an equivariant product decomposition for normal varieties whose orbits are complete and
 have the same dimension. More precisely, under these assumptions, such a variety is equivariantly isomorphic to the product of a parabolic variety and a closed subvariety. This extends \cite[Theorem 3.3]{SS}  by removing one of its technical assumptions.

\smallskip

\noindent {\bf Examples and questions over imperfect fields.} 
In the final section of this paper, we examine pathological phenomena 
that arise over imperfect ground fields. The properties stated in Proposition \ref{prop:rational intro} do not extend to this setting. Indeed, 
considering the example of $X_m$ defined in (\ref{eqn:tiv}), with $n = 2$,
one can show that, as soon as $m \geq 2$, the set $X_m(k)$ contains 
infinitely many $\SL_3(k)$-orbits. We then generalize this 
to any split adjoint semisimple group $G$ of rank at least two, 
obtaining parabolic varieties $X=G/P$ for which $G(k)$ has infinitely many orbits in $X(k)$ (Proposition \ref{prop:imperfect no torus}). Moreover, there exists a $k$-rational point of $X$ whose stabilizer $Q$ contains no maximal $k$-torus of $G$; in particular, the corresponding stabilizers are conjugate over the algebraic closure but not over $k$, even if the latter is 
separably closed. Finally, these varieties admit no universal purely inseparable cover.

A second construction (Proposition \ref{prop:isotropic}) produces, for certain 
simple adjoint groups $G$, parabolic varieties whose stabilizer of any $k^\sep$-point contains no Borel subgroup. We then obtain two families of examples: the first arises from a central division algebra $D$ of degree $p$, while the second occurs in characteristic $2$, starting from an anisotropic quadratic space of odd dimension.

\medskip

We conclude by outlining several open questions, though this list is by no means exhaustive. To complete the picture of twisted parabolic $G$-varieties over
an arbitrary field $k$, it would be of interest to characterize those having 
a $k$-rational point and to describe the $G(k)$-orbits of such points.
The latter question amounts to classifying the $G(k)$-conjugacy classes 
of parabolic subgroups. In addition, one would like to determine the 
largest smooth subgroups of parabolic subgroups in this setting. 
In particular, the goal is to establish structural results that yield 
a explicit geometric description of these subgroups, rather than merely 
providing counterexamples to properties familiar over algebraically closed 
or perfect fields.

Another direction concerns the second class of examples mentioned above. 
We expect that such phenomena should occur only in small characteristics, 
possibly in relation with numerical invariants of the group such as its 
Coxeter number. Understanding this relationship and determining the precise 
conditions under which these examples can occur remain open questions.

\color{black}
\medskip \noindent
{\bf Outline.}
The present text is organized as follows. Section \ref{sec:hv} introduces the terminology and preliminaries concerning homogeneous varieties. Section \ref{sec:spic} deals with their structure, separating the semisimple case from the case of torsors under abelian varieties, and with the description of their Picard groups. Sections \ref{sec:pvscf} and \ref{sec:pvaf} classify parabolic varieties, first over a separably closed field and then over an arbitrary field. Their twisted forms are described in Section \ref{sec:tpvaf}. The equivariant rigidity results are established in Section \ref{sec:rigidity}. Finally, Section \ref{sec:example} constructs the classes of examples over imperfect fields discussed above.

\medskip \noindent
{\bf Notation and conventions.}
We consider schemes over a ground field $k$ of 
characteristic $p \geq 0$. Products and morphisms 
of schemes are understood to be over $k$, unless 
mentioned otherwise. 
The underlying reduced subscheme of a scheme 
$X$ is denoted by $X_{\red}$.

We choose an algebraic closure $\bar{k}$  
of $k$ and denote by 
$k^{\sep} \subset \bar{k}$ 
the corresponding separable closure, 
with Galois group
$\Gamma_k := \Gal(k^{\sep}/k) 
= \Aut_k(\bar{k}).$ 
Given a scheme $X$ and a field extension 
$k \subset K$, we denote by $X_K$ the
$K$-scheme $X \times_{\Spec(k)} \Spec(K)$.
When $K = \bar{k}$, we also use the notation
$\bar{X}$ for $X_{\bar{k}}$. Likewise,
we write $X^{\sep}$ for $X_{k^{\sep}}$. 

A \emph{variety} is a separated,
geometrically integral scheme $X$
of finite type over $k$. We say that $X$ is 
\emph{complete} if the structure morphism 
$X \to \mathrm{Spec}(k)$ is proper.

An \emph{algebraic group} is a
group scheme $G$ of finite type over $k$.
We also say that $G$ is a $k$-\emph{group},
especially if we want to keep track of the
ground field in the presence of field
extensions. A $k$-\emph{subgroup} $H$ of $G$  
is a subgroup scheme, so $H$ is closed but not 
necessarily smooth. 
If $k$ is perfect then the reduced subscheme 
$G_{\red}$ is a smooth subgroup scheme
(see e.g. \cite[II.5.2.3]{DG}).

Assume $p>0$. For any positive integer $m$, the
$m$-th relative Frobenius morphism
of a scheme $X$ is denoted by
\[ F^m = F^m_{X/k} : 
X \longrightarrow X^{(m)}. \]
If $G$ is a group scheme over $k$, then
so is $G^{(m)}$ and $F^m_{G/k}$ is 
a homomorphism of group schemes.
We denote its kernel by $_mG$
(the $m$th Frobenius kernel).

\medskip

\noindent \textbf{Acknowledgements.} 
	The authors express their sincere gratitude to Charles De Clercq, Jean Fasel and Philippe Gille for insightful feedback and valuable discussions on this topic.
The second author was supported by the DFG through the research grants Le 3093/5-1 and Le
3093/7-1 (project numbers 530132094; 550535392). The third author acknowledges the support of the DAE, Government of India, under Project Identification No.~RTI4014.

\section{Homogeneous varieties}
\label{sec:hv}

\subsection{Definitions and preliminary results}
\label{subsec:dpr}

Throughout this section, we consider a smooth 
algebraic group $G$ and a variety $X$ 
equipped with a $G$-action 
\[ \alpha : G \times X \longrightarrow X,
\quad (g,x) \longmapsto g \cdot x. \]
We then say that $X$ is a $G$-\emph{variety}.
Given $x \in X(k)$, we denote by 
\begin{equation}\label{eqn:orbitmap} 
\alpha_x : G \longrightarrow X, 
\quad g \longmapsto g \cdot x 
\end{equation}
the corresponding \emph{orbit map}, with image 
the \emph{orbit} $G \cdot x$ (a locally closed
subset of $X$, equipped with the structure
of a smooth subscheme). Also, we denote by 
\[
\gamma : G \times X \longrightarrow
X \times X, \quad 
(g,x) \longmapsto (g \cdot x, x) 
\]
the associated \emph{graph morphism}. 
The \emph{stabilizer} $\Stab_G$ is defined
via the cartesian square
\begin{equation}\label{eqn:stabi}
\xymatrix{
\Stab_G \ar[r] \ar[d] & X \ar[d]_-{\Delta} \\
G \times X \ar[r]^-{\gamma} & X \times X, \\
} \end{equation}
where $\Delta$ denotes the diagonal immersion.
Note that $\Stab_G$ is a $X$-subgroup scheme
of the group scheme $\pr_X : G \times X \to X$.
Also, $\gamma$ is equivariant for the action of 
$G \times G$ on $G \times X$ via
\[ 
(g_1,g_2) \cdot (g,x) := (g_1 g g_2^{-1}, g_2 \cdot x)
\]
and its action on $X \times X$ via
\[
(g_1,g_2) \cdot (x_1,x_2) = 
(g_1 \cdot x_1, g_2 \cdot x_2).
\]
Moreover, $\Stab_G \subset G \times X$
is stable under the action of $G$ via 
\[
g_1 \cdot (g,x) := (g_1 g g_1^{-1}, g_1 \cdot x),
\]
and $\pr_X$ is equivariant for this action.

We say that $X$ is \emph{homogeneous under $G$} 
if $\gamma$ is surjective. We now collect 
some preliminary results on this notion, 
starting with criteria for homogeneity due to 
Bruno Laurent:

\begin{lemma}\label{lem:homogeneous}
\cite[Lemma 2.2]{Laurent}

Let $k\subset K$ be a field extension. 
The following assertions are equivalent:

\begin{enumerate}

\item[(i)] 
The variety $X$ is homogeneous under $G$.

\item[(ii)]
The variety $X_K$ is homogeneous under $G_K$.

\item[(iii)]
The action of the abstract group $G(\bar{k})$ 
on the set $X(\bar{k})$ is transitive.

\item[(iv)]
If $X$ has a $K$-rational point $x$, 
then the orbit $G_K \cdot x$ equals $X_K$.

\end{enumerate}

\end{lemma}

In particular, if $X$ is a homogeneous variety
equipped with a $k$-rational point $x$, then
$X = G \cdot x$. Also, recall from
\cite[III.3.2.1]{DG} that the orbit map 
$\alpha_x$ factors as the quotient map 
$q : G \to G/H$ (a fpqc $H$-torsor)
followed by a $G$-equivariant isomorphism 
\[ G/H \stackrel{\sim}{\longrightarrow} X,
\quad o \longmapsto x, \]
where $H := \Stab_G(x)$ is the $k$-subgroup 
of $G$ consisting of those $g$ 
such that $g \cdot x = x$, and $o := q(e)$
is the \emph{base point} of the 
\emph{homogeneous space} $G/H$. 
This identifies the pairs $(X,x)$,
where $X$ is a $G$-homogeneous variety and
$x \in X(k)$, with the homogeneous spaces
$G/H$, where $H \subset G$ is a $k$-subgroup.

\begin{lemma}\label{lem:ff}
Let $X$ be a homogeneous $G$-variety. Then $X$ 
is smooth and quasi-projective. Moreover, the graph 
morphism $\gamma$ is faithfully flat.
\end{lemma}

\begin{proof}
Every homogeneous space is smooth and 
quasi-projective; moreover, both properties
descend under base change by field extensions.
This implies the first assertion.

To show the second assertion, we may assume 
that $k$ is algebraically closed as the
formation of $\gamma$ commutes with
base change by field extensions. Then the above
action of $G \times G$ on $X \times X$ yields 
a transitive action of $G(k) \times G(k)$  
on $X(k) \times X(k)$ in view of Lemma 
\ref{lem:homogeneous}. Moreover, by generic flatness
(see \cite[Th\'eor\`eme 6.9.1]{EGAIV2}), 
there exists a dense open subset
$U \subset X \times X$ such that the pull-back
$\gamma^{-1}(U) \to U$ is flat. As $X \times X$
is covered by the translates $(g_1,g_2) \cdot U$
where $g_1, g_2 \in G(k)$ and as $\gamma$
is $G \times G$-equivariant, this morphism 
is flat everywhere. 
\end{proof}

In particular, every complete homogeneous
variety is smooth and projective.

\begin{lemma}\label{lem:quotient}
Let $N \triangleleft G$ a normal subgroup, 
and $X$ a homogeneous $G$-variety. 
Then the categorical quotient 
$f : X \to X/N$ 
exists, and $X/N$ is a homogeneous 
$G/N$-variety. If $X =G/H$ then
$X/N = G/(N \cdot H)$. 
\end{lemma}

\begin{proof}
Since $X$ is a smooth variety, it has a 
$K$-rational point $x$ for some finite 
Galois extension $k \subset K$. 
Thus, we may identify the pair $(X_K,x)$ 
with the homogeneous space $G_K/H$
equipped with its base point, where 
$H := \Stab_{G_K}(x)$. The quotient morphism
$G_K \to G_K/(N_K \cdot H_K)$ is $G_K$-equivariant
and $H_K$-invariant,
and hence factors uniquely through a
$G_K$-equivariant morphism 
\[ \varphi : X_K = G_K/H_K \to G_K/(N_K \cdot H_K) 
= X_K/N_K. \]
Moreover, $\varphi$ is the categorical quotient
by $N_K$ (see e.g. \cite[Lemma 3.2]{Brion2}).
So the assertion follows by Galois descent, 
which applies as $X_K$ and $X_K/N_K$ are 
quasi-projective.
\end{proof}

\begin{lemma}\label{lem:connected}
Let $X$ be a homogeneous $G$-variety.
Then $X$ is homogeneous under the neutral 
component $G^0$. 
\end{lemma}

\begin{proof}
We may assume that $X = G/H$ in view of
Lemma \ref{lem:homogeneous}.
Then $X/G^0 = G/(G^0 \cdot H)$ by Lemma
\ref{lem:quotient}. Thus, $X/G^0$ is finite, since 
so is $G/G^0$. But $X/G^0$ is also geometrically 
integral, since so is $X$. Therefore, 
$X/G^0 = \Spec(k)$, that is, $G = G^0 \cdot H$,
and $X = G^0/G^0 \cap H$.
\end{proof}

\subsection{Equivariant automorphisms}
\label{subsec:equiaut}

\begin{lemma}\label{lem:autnorm}
Let $H \subset G$ be a $k$-subgroup,
and $N_G(H)$ its (schematic) normalizer. Then 
the $N_G(H)$-action on $G$ by right
multiplication induces a faithful action 
of $N_G(H)/H$ on $G/H$ that commutes with
the $G$-action. Moreover, $N_G(H)/H$
represents the group functor 
\[ \uAut_{G/H}^G : 
S \longmapsto \Aut_S^G(G/H \times S) \]
of $G$-equivariant automorphisms of $G/H$. 
\end{lemma}

\begin{proof}
Composing the multiplication map
$m : G \times N_G(H) \to G$ 
with the quotient map 
$q : G \to G/H$ yields a morphism
$\varphi : G \times N_G(H) \to G/H$ 
which is invariant under the action 
of $H \times H$ via 
$(h_1,h_2) \cdot (g,n) 
= (g h_1^{-1}, n h_2^{-1})$. Since
the quotient map 
$G \times N_G(H) \to G/H \times N_G(H)/H$
is a $H \times H$-torsor, $\varphi$ factors
uniquely through a morphism 
\[ \psi : G/H \times N_G(H)/H \to G/H. \]
As $m$ is a right action of $N_G(H)$
that commutes with the $G$-action by
left multiplication, it follows by
fpqc descent that $\psi$ is a right action 
that commutes with the $G$-action.
Viewing the group scheme $N_G(H)/H$ as 
a group functor, we obtain a homomorphism
of such functors 
\[ \eta : N_G(H)/H \longrightarrow \uAut^G_{G/H}. \]
We check that $\eta$ is an isomorphism:
Let $S$ be a scheme, and 
$u \in \Aut_S^G(G/H \times S)$.
Let $v : G/H \times S \to G/H$ 
be the morphism such that 
$u(x,s) = \big( v(x,s),s \big)$ identically on 
$G/H \times S$. Denoting by 
$o$ the base point of $G/H$, 
we have $v(o,s) \in (G/H)(S)$
and hence $v(o,s) = g \cdot o$
for some $g \in G(S')$, where
$S' \to S$ is a faithfully flat,
quasi-compact morphism, and
$o, v(o,s)$ are viewed as $S'$-points 
of $G/H$. Since $o$ is fixed by $H$ 
and $u$ is $G$-equivariant, $g \cdot o$
is fixed by $H_{S'}$. It follows
that $g^{-1} H_{S'} g \subset H_{S'}$,
i.e., $g \in N_G(H)(S')$. Thus,
$v(o,s) \in (N_G(H)/H)(S)$ and 
hence $u \circ v(o,s)^{-1}$ is
a $G_S$-equivariant automorphism
of $(G/H)_S$ that fixes the base
section $o \times S$.
By homogeneity, it follows that
$u = v(o,s)$. A similar argument 
shows that the kernel of $\eta$
is trivial, completing the proof.
\end{proof}

\begin{remark}\label{rem:norm}
The above lemma holds for any algebraic group 
$G$, not necessarily smooth. Also, note that
the homogeneous space $G/N_G(H)$ may be seen 
as the scheme of conjugates of $H$ in $G$. 
In particular, we may see $G/H$ as the scheme 
of conjugates of $H$, under the assumption
that $N_G(H) = H$. 

The latter assumption is satisfied 
for any parabolic $k$-subgroup $P$ of a 
semisimple group $G$. Indeed, this holds if
$k$ is algebraically closed by 
\cite[Corollary 3.3]{MaccanTerpereau},
and hence over an arbitrary field as
the formation of $N_G(H)$ commutes with 
field extensions.
\end{remark}

\begin{proposition}\label{prop:authom}
Let $X$ be a homogeneous $G$-variety. 

\begin{enumerate}
\item[{\rm (i)}] 
The group functor of $G$-equivariant automorphisms,
\[ \uAut^G_X : 
S \longmapsto \Aut^G_S(X \times S), \] 
is represented by an algebraic group $\bfAut^G_X$.

\item[{\rm (ii)}] 
The formation of $\bfAut^G_X$ commutes with base change
under field extensions.

\item[{\rm (iii)}] If $G$ is semisimple and
$X$ is complete, then $\bfAut^G_X$ is trivial.

\end{enumerate}

\end{proposition}

\begin{proof}
(i) The case where $X$ has a rational point
follows from Lemma \ref{lem:autnorm}.
In the general case, $X$ has a $K$-rational
point for some finite Galois extension
$k \subset K$, and hence the base change functor
$\uAut^{G_K}_{X_K}$
is represented by an algebraic $K$-group
equipped with an action of the Galois group
$\Gamma_{K/k}$. Our assertion now follows
by Galois descent, since every algebraic
group is quasi-projective.

(ii) This follows from the corresponding fact
for the group functor $\uAut^G_X$.

(iii) By (ii), we may assume that $X$
has a $k$-rational point. Then the assertion 
follows by combining Lemma \ref{lem:autnorm}
and Remark \ref{rem:norm}.
\end{proof}

\subsection{Homogeneous varieties with smooth stabilizer}
\label{subsec:hss}

We still consider a smooth algebraic group $G$.

\begin{lemma}\label{lem:smooth}
The following are equivalent
for a homogeneous $G$-variety $X$:

\begin{enumerate}

\item[{\rm (i)}] 
The graph morphism $\gamma$ is smooth.

\item[{\rm (ii)}] 
The stabilizer $\Stab_G$ is smooth over $X$.

\item[{\rm (iii)}] 
The fiber $\Stab_G(x)$ is smooth for some 
$x \in X(\bar{k})$.

\end{enumerate}

\end{lemma}

\begin{proof}
The implication (i) $\Rightarrow$ (ii) follows
from the cartesian square (\ref{eqn:stabi}),
and (ii) $\Rightarrow$ (iii) is obvious.
To show that (iii) $\Rightarrow$ (i), 
recall that $\gamma$ is flat (Lemma 
\ref{lem:ff}) and its fiber at 
$(x,x)$ is isomorphic to $\Stab_G(x)$. Thus,
$\gamma$ is smooth at $(x,x)$, and hence 
everywhere by a homogeneity argument as in
the proof of that lemma.
\end{proof}

\begin{definition}\label{def:hss}
We say that a $G$-variety $X$ is
\emph{homogeneous with smooth stabilizer}
(hss for short) if it satisfies the 
equivalent conditions of Lemma \ref{lem:smooth}.
\end{definition}

The property of being hss can also be characterized 
in terms of the Lie algebra $\Lie(G)$, as follows.

\begin{lemma}\label{lem:lie}
A smooth $G$-variety $X$ is hss if and only if
the tangent sheaf $\cT_X$ is generated by its 
global sections arising from $\Lie(G)$.
\end{lemma}

\begin{proof}
Since the formations of the map
$\Lie(G) \to \Gamma(X,\cT_X)$
and of the associated map 
$\cO_X \otimes_k \Lie(G) \to \cT_X$
commute with base change by field extensions,
we may assume that $k$ is algebraically
closed. For any $x \in X(k)$, consider the
differential at $x$ of the orbit map 
(\ref{eqn:orbitmap}) 
\[ d\alpha_x : \Lie(G) \longrightarrow \cT_x(X), \]
where the right-hand side is the Zariski tangent
space at $x$.
If $X$ is hss, then $\alpha_x$ is smooth and
hence $d\alpha_x$ is surjective. As a consequence,
$\Lie(G)$ generates $\cT_X$. Conversely,
if $\Lie(G)$ generates $\cT_X$ then each
$\alpha_x$ is smooth at $x$, and hence
everywhere by equivariance. Thus, its
fiber $\Stab_G(x)$ is smooth, so that
$X$ is hss by Lemma \ref{lem:smooth}.
\end{proof}

\begin{remark}\label{rem:hss}
As a direct consequence of Lemmas \ref{lem:homogeneous}
and \ref{lem:smooth}, a $G$-variety $X$ is hss
if and only if so is the $G_K$-variety $X_K$
for some field extension $k \subset K$.

Also, the action of $G(k^{\sep})$ on $X(k^{\sep})$
is transitive for any hss $G$-variety $X$. 
Indeed, we may choose $x \in X(k^{\sep})$ 
and the resulting orbit map 
$\alpha_x : G_{k^{\sep}} \to X_{k^{\sep}}$
is smooth, hence surjective on $k^{\sep}$-points.

The above transitivity property does not
characterize hss varieties over an imperfect
field $k$. Indeed, there exist non-trivial $k$-forms
of the additive group $\bG_{a,k}$; moreover, every such
form $G$ admits an equivariant completion by 
a regular curve $X$, and the complement $X \setminus G$ 
consists of a closed point $x$ such that the extension 
$k \subset \kappa(x)$ is non-trivial and purely 
inseparable (see \cite{Russell}). In particular,
$X(k^{\sep}) = G(k^{\sep})$ but $X$ is not $G$-hss.
\end{remark}

We now come to the main result of this subsection:

\begin{proposition}\label{prop:hss}
Let $X$ be a homogeneous $G$-variety
over a perfect field $k$.

\begin{enumerate}

\item[{\rm (i)}] 
There exists a $G$-equivariant morphism 
$f : \widehat{X} \to X$, where $\widehat{X}$
is a hss $G$-variety, which satisfies the 
following universal property: For any 
equivariant morphism $u : Y \to X$, where 
$Y$ is a hss $G$-variety, there exists
a unique $G$-equivariant morphism 
$\widehat{u} : Y \to \widehat{X}$ such that
$u = f \circ \widehat{u}$. 

\item[{\rm (ii)}] 
The morphism $f$ is finite locally free and 
purely inseparable.

\item[{\rm (iii)}]
The map $f(K) : \widehat{X}(K) \to X(K)$
is bijective for any perfect field $K$
containing~$k$. In particular, $X$
has a $k$-rational point if and only if
so does $\widehat{X}$.

\item[{\rm (iv)}]
Assume that $X$ has a $k$-rational point $x$.
Then $f$ is the canonical morphism
$G/H_{\red} \to G/H$, where $H := \Stab_G(x)$. 

\end{enumerate}

\end{proposition}

\begin{proof}
(i) We first show the existence of the 
morphism $f$ in the case where $k$ is
algebraically closed. We may then choose
$x \in X(k)$ and identify $X$ with $G/H$,
where $H := \Stab_G(x)$. Given an equivariant
morphism $u : Y \to X$, where $Y$ is a hss
$G$-variety, there exists $y \in Y(k)$
such that $u(y) = x$. Identifying $Y$ with
$G/I$ where $I := \Stab_G(y)$, we then have 
that $I \subset H$ (as $u$ is equivariant),
and $I$ is smooth (as $Y$ is hss). So
$I \subset H_{\red}$ and hence $u$ factors
uniquely through the canonical morphism
$f : G/H_{\red} \to G/H$. Thus, $f$ is
the desired morphism.

The existence of $f$ in the case of a perfect 
field $k$ follows from this by a standard 
argument of Galois descent, as in 
\cite[Chapitre V, \S 4, no.~22]{Serre}. More 
specifically, consider the universal morphism 
$f : \widehat{\bar{X}} \to \bar{X}$
and set for simplicity $V := \widehat{\bar{X}}$.
For any $\sigma \in \Gamma_k$, let $V^{\sigma}$
be the $\bar{k}$-variety obtained from $V$ 
by base change under 
$\sigma : \Spec(\bar{k}) \to \Spec(\bar{k})$.
Then the composition
\[ V^{\sigma} \stackrel{\sigma}{\longrightarrow}
V \stackrel{f}{\longrightarrow} 
\bar{X} \stackrel{\sigma^{-1}}{\longrightarrow}
\bar{X} \]
is an equivariant $\bar{k}$-morphism from a hss 
$G$-variety to $\bar{X}$, and hence factors as 
$f \circ h_{\sigma}$ for a unique
$\bar{k}$-morphism 
$h_{\sigma} : V^{\sigma} \to V$. 
This uniqueness implies the cocycle condition
$h_{\sigma \tau} = 
(h_{\sigma})^{\tau} \circ h_{\tau}$
for all $\sigma, \tau \in \Gamma$.
By Galois descent (which applies since
$V$ is quasi-projective), it follows
that $V = \bar{W}$ for a unique variety
$W$. Likewise, the $\bar{G}$-action on $V$ 
descends to a $G$-action on $W$ which makes 
it a hss variety, and $f$ descends to the 
desired universal morphism $W \to X$.

(ii) We may assume that $k$ is algebraically 
closed by using fpqc descent 
(see \cite[Proposition 2.5.2]{EGAIV2} for 
``finite locally free'', and loc.~cit.,
Proposition 2.6.1 together with 
\cite[Proposition 3.5.7]{EGAI}
for ``purely inseparable'').  
Then $f : G/H_{\red} \to G/H$ lies in
a diagram with cartesian squares
\[ \xymatrix{ H/H_{\red} \ar[r] & 
\Spec(k)  \\
G \times (H/H_{\red}) \ar[r]^-{\pr_1} 
\ar[d]_m \ar[u]^-{\pr_2} 
& G \ar[d]^q \ar[u] \\
G/H_{\red} \ar[r]^f & G/H, \\
}\]
where $\pr_1$ denotes the projection,
$q$ the quotient map, and $m$
is induced by the multiplication
map $G \times H \to G$
(see e.g.~\cite[(2.8.3)]{Brion1}).
Moreover, the structure map 
$H/H_{\red} \to \Spec(k)$ is finite
locally free and purely inseparable; 
thus, so is $f$ by fpqc descent again.

(iii) Given $x \in X(K)$, the fiber $\widehat{X}_x$ 
is the spectrum of a local $K$-algebra $A$, 
which is finite-dimensional as a $K$-vector 
space and whose residue field $L$ is
a purely inseparable extension of $K$. 
Since $K$ is perfect, $L = K$ and hence 
$\widehat{X}_x$ has a unique $K$-point, 
corresponding to the unique algebra homomorphism 
$A \to K$. 

(iv) Since the formation of $H_{\red}$
commutes with perfect field extensions,
the base change to $\bar{k}$ of
the canonical morphism $G/H_{\red} \to G/H$
is the morphism 
$\bar{G}/\bar{H}_{\red} \to \bar{G}/\bar{H}$.
As seen in the proof of (i), the latter morphism
is the universal morphism for $\bar{X}$,
and descends to the universal morphism for $X$.
This implies our assertion by uniqueness in
Galois descent.
\end{proof}

\begin{definition}
\label{def purely inseparable cover}
With the notation and assumptions of Proposition 
\ref{prop:hss}, we say that $\widehat{X}$ is the 
\emph{universal purely inseparable cover} of $X$. 
\end{definition}

One may show by similar arguments
that $f : \widehat{X} \to X$ factors uniquely 
through any purely inseparable equivariant morphism 
$Z \to X$, where $Z$ is a homogeneous $G$-variety. 
On the other hand, if $k$ is imperfect 
then there exist complete homogeneous varieties
having no universal purely inseparable cover, 
as we will see in Section \ref{sec:example}.

Finally, we record a preliminary result which will 
be used in Section \ref{sec:rigidity}. It is a (small) 
generalization of \cite[Theorem 4.3]{Brosnan}, see also 
\cite[Theorem A]{Richardson}.

\begin{lemma}\label{lem:fixed}
Let $X$ be a homogeneous $G$-variety. Let $T \subset G$
be a $k$-torus, with centralizer $C_G(T) \subset G$ 
and fixed point subscheme $X^T \subset X$. Then $X^T$ 
is a finite union of open and closed homogeneous 
$C_G(T)$-subvarieties. 
\end{lemma}

\begin{proof}
Since $G$ and $X$ are smooth, so are $C_G(T)$ 
(the fixed point subscheme of $G$ under the 
conjugation action of $T$) and $X^T$.
Moreover, $X^T$ is clearly stable by $C_G(T)$.

We first treat the case where $k$ is algebraically
closed. Then every homogeneous variety is an
orbit of a $k$-rational point, and hence it
suffices to show that every $C_G(T)$-orbit in $X^T$ 
is open (taking complements, it follows
that every such orbit is closed, and hence
there are finitely many of them). For this, 
we use the universal purely inseparable cover 
$\widehat{X} \to X$. 
It restricts to a finite $C_G(T)$-equivariant morphism 
$\widehat{X}^T \to X^T$, which is surjective  
by Borel's fixed point theorem. Thus, we may further
assume that the stabilizer of $X$ is smooth. 
For any $x \in X^T(k)$, the differential of 
the orbit map $\Lie(G) \to \cT_x(X)$ 
is surjective (Lemma \ref{lem:lie}), and hence 
restricts to a surjective map
\[ \Lie \big( C_G(T) \big) = \Lie(G)^T \longrightarrow
\cT_x(X)^T = \cT_x(X^T). \]
Using Lemma \ref{lem:lie} again, it follows
that the orbit $C_G(T) \cdot x$ is open in $X^T$
as desired.

In the case of an arbitrary field $k$,
there exist finitely many points
$x_1,\ldots,x_n \in X^T(\bar{k})$
such that $X^T_{\bar{k}}$ is the sum
of the (open and closed) orbits
$C_G(T)_{\bar{k}} \cdot x_i$. Let $k \subset K$
be a finite extension such that 
$x_1,\ldots,x_n \in X(K)$. Since the projection
$X_{\bar{k}} \to X_K$ is open and surjective,
$X_K$ is the sum of the open orbits
$C_G(T)_K \cdot x_i$, which are thus closed.
This yields our assertion by using the
projection $X_K \to X$ (which is finite
locally free, and hence open and closed).
\end{proof}

\section{Structure and Picard group of complete 
homogeneous varieties}
\label{sec:spic}

Throughout this section, we consider a complete
variety $X$, homogeneous under a smooth algebraic
group $G$. We may assume that $G$ acts faithfully 
on $X$ and is connected, by Lemmas \ref{lem:quotient}
and \ref{lem:connected}.

\begin{theorem}\label{thm:nonaffine}
Under the above assumptions, there exists 
a canonical isomorphism of algebraic groups
\begin{equation}\label{eqn:groupis} 
G \stackrel{\sim}{\longrightarrow} 
A \times H, 
\end{equation}
where $A$ is an abelian variety and $H$ is
an adjoint semisimple group,
and a canonical isomorphism of varieties
\begin{equation}\label{eqn:varis}
X \stackrel{\sim}{\longrightarrow} 
Y \times Z, 
\end{equation}
where $Y$ is an $A$-torsor and $Z$ is
a complete variety, homogeneous under
$H$. Moreover, these isomorphisms are
equivariant for the action of $A \times H$
on $Y \times Z$ via 
$(a,h) \cdot (y,z) := (a\cdot h, y \cdot z)$.
\end{theorem}

\begin{proof}
In the case where $k$ is algebraically closed, 
this is a refinement of a result of Sancho de Salas 
(see \cite[Theorem 5.2]{SS} for that result, 
and \cite[Theorem 4.1.1]{BSU} for the refinement).
We may then take $Y = A$ on which $A$ acts
by translations.

In the general case, we consider the center
$Z(G)$, the derived subgroup $\cD(G)$, 
the two quotient homomorphisms by these normal
subgroups, and the resulting homomorphism
\[ u : G \longrightarrow 
G/\cD(G) \times G/Z(G). \]
Then the existence of the isomorphism
(\ref{eqn:groupis}) is equivalent to the 
following claim: $G/\cD(G)$ is an abelian 
variety, $G/Z(G)$ is an adjoint semisimple
group, and $u$ is an isomorphism.
Indeed, (\ref{eqn:groupis}) identifies
$Z(G)$ to $A$, and $\cD(G)$ to $H$.
Thus, the claim holds after base change
to $\bar{k}$ (since the formations of $Z(G)$ 
and $\cD(G)$ commute with base change 
by field extensions). By fpqc descent, it
follows that the claim holds over $k$,
proving the assertion on $G$.

The proof of the assertion on $X$ is similar:
by Lemma \ref{lem:quotient}, we have 
categorical quotients
$X \to X/\cD(G)$, $X \to X/Z(G)$, 
where $X/\cD(G)$ (resp.~$X/Z(G)$) is 
homogeneous under $G/\cD(G)$ (resp.~$G/Z(G)$).
These quotients yield an equivariant morphism
\[ v : X \longrightarrow 
X/\cD(G) \times X/Z(G).\]
Then the existence of a canonical 
isomorphism (\ref{eqn:varis})
is equivalent to the following claim:
$v$ is an isomorphism and $X/\cD(G)$ 
is a torsor under $G/\cD(G)$.
Since the formation of the categorical 
quotient commutes with base change by
field extensions, so does that of $v$.
This reduces the proof of the claim to 
the known case where $k$ is algebraically
closed.
\end{proof}

Next, we will describe the Picard group $\Pic(X)$ in terms
of the product decomposition (\ref{eqn:varis}). 
For this, we will freely use the notation  
of the above theorem. We will also freely identify 
$\Pic(X)$ with the divisor class group $\Cl(X)$.

We begin with the case where $k$ is 
\emph{separably closed}. Then $Y$ has a $k$-rational
point, and hence is a trivial $A$-torsor. 
We now obtain some geometric information
on the $H$-variety $Z$; see Section 
\ref{sec:pvscf} for further developments.
Note that the semisimple $k$-group $H$ is split,
and hence contains two opposite Borel subgroups 
$B,B^-$, with unipotent radicals $U,U^-$.

\begin{proposition}\label{prop:rational}
Keep the above notation and assumptions.

\begin{enumerate}

\item[{\rm (i)}]
$Z$ has a unique $k$-rational $B$-fixed
point $z$. Moreover, the orbit $G(k) \cdot z$ 
consists of the points of $Z(k)$ whose 
stabilizer contains a smooth $k$-parabolic.

\item[{\rm (ii)}]
The orbit $V := U^- \cdot z$ is open in $Z$ and
isomorphic to an affine space. 

\item[{\rm (iii)}]
$Z$ is covered by the translates $h \cdot V$,
where $h \in H(k)$. 

\end{enumerate}

\end{proposition}

\begin{proof}
(i) Since $k$ is separably closed, $Z$ has 
a $k$-rational point. As $B$ is a split solvable
$k$-group, it fixes some point $z \in Z(k)$ 
by a version of the Borel fixed point theorem 
(see \cite[Lemma, p.~35]{Ro57}). If $z' \in Z(k)$
is another $B$-fixed point, then there exists
$h \in H(\bar{k})$ such that $z' = h \cdot z$.
Let $P := \Stab_G(z)$, then both $\bar{B}$ and 
$h^{-1} \bar{B} h$ are contained in $\bar{P}$, 
and hence are Borel subgroups of the
smooth $\bar{k}$-subgroup 
$\bar{P}_{\red} \subset \bar{G}$. So there exists 
$p \in \bar{P}_{\red}(\bar{k}) = P(\bar{k})$
such that $h^{-1} \bar{B} h = p \bar{B} p^{-1}$.
Thus, 
$hp \in N_{\bar{G}}(\bar{B})(\bar{k}) = B(\bar{k})$.
Equivalently, $h \in P(\bar{k})$, and $z' = z$.
Thus, $z$ is the unique $k$-rational $B$-fixed point.

Let $z'' \in Z(k)$. If $\Stab_G(z'')$ contains
a smooth $k$-parabolic, then it also contains
a Borel $k$-subgroup. As any two such subgroups
are conjugate in $G(k)$, there exists $g \in G(k)$
such that $g \cdot z''$ is $B$-fixed, and hence
$z'' \in G(k) \cdot z$ by the above step. 
Conversely, if $z'' \in G(k) \cdot z$ then
its stabilizer contains a Borel $k$-subgroup.

(ii) The multiplication map 
$U^- \times B \to H$ is an open immersion.
Thus, $U^- \cdot z$ is open in $Z = H \cdot z$.
Moreover, the unipotent $k$-group $U^-$ is split 
in view of \cite[Corollaire 3.18]{BT}, and hence 
the $U^-$-orbit $V$ is an affine space 
(as follows from \cite[Theorem 5]{Ro63},
see \cite{Brion2} for a modern presentation
of this result).

(iii) The subset $Z \setminus H(k) \cdot V$ 
is closed in $Z$ and stable by $H(k)$. The latter
is dense in the semisimple group $H$.
Therefore, $Z \setminus H(k) \cdot V$ (equipped
with its closed reduced subscheme structure)
is stable by $H$, and hence empty as $Z$
is homogeneous under $H$.
\end{proof}

Since the open subset $V \subset Z$ is affine, 
the irreducible components of its complement 
$Z \setminus V$ are prime divisors, 
which we denote as $D_1,\ldots,D_r$.

\begin{proposition}\label{prop:picrat}
We still assume that $k$ is separably closed.
Then:

\begin{enumerate}

\item[{\rm (i)}] 
The abelian group $\Pic(Z)$ is freely generated
by the classes of $D_1,\ldots,D_r$, and these generate 
the monoid of effective divisor classes. 

\item[{\rm (ii)}] 
Every effective divisor on $Z$ is base-point-free.

\item[{\rm (iii)}] 
The map 
\begin{equation}\label{eqn:pullback} 
\pr_Y^* \times \pr_Z^* :
\Pic(Y) \times \Pic(Z) \longrightarrow \Pic(X)
\end{equation}
is bijective.

\end{enumerate}

\end{proposition}

\begin{proof}
(i) We adapt a standard approach to 
the Picard group of flag varieties.
Denote by $H_{\sic}$ the simply-connected
cover of $H$. Then the character group of $H_{\sic}$ 
is trivial, and $\Pic(H_{\sic}) = 0$ in view of 
\cite[Lemme 6.9]{Sansuc}. 
By \cite[Proposition 2.10, Lemma 2.13]{Brion0},
it follows that every line bundle $L$ on $X$ 
has a unique $H_{\sic}$-linearization, where 
$H_{\sic}$ acts on $X$ through its quotient
$H$. Thus, the space of global sections 
$\Gamma(Z,L)$ is a finite-dimensional 
$H_{\sic}$-module. If $L$ is effective, 
then of course $\Gamma(Z,L) \neq 0$.
Consider the pull-back of $B^-$ in $H_{\sic}$;
this is a split solvable $k$-group $B^-_{\sic}$, 
and hence has an eigenvector $s$ in $\Gamma(Z,L)$.
In particular, $s$ is fixed by the unipotent
radical $U^-_{\sic} \subset H_{\sic}$, which lifts 
$U^- \subset H$. As $V$ is an orbit of $U^-$, 
the divisor of $s$ is a linear combination of 
$D_1,\ldots,D_r$ with non-negative integer 
coefficients. This shows that the monoid
of effective classes is generated by those
of $D_1,\ldots,D_r$. These are linearly 
independent, since $V$ is an affine space
and hence $\cO(V)^{\times} = k^{\times}$.

(ii) Let again $L$ be an effective line
bundle on $X$. Since $L$ is 
$H_{\sic}$-linearizable, its base locus is 
stable under $H_{\sic}$, and hence empty 
by homogeneity.

(iii) The pull-back map 
$\Cl(Y \times Z) \to \Cl(Y \times V)$
is surjective and its kernel is generated
by the classes of the divisors 
$Y \times D_i = \pr_Z^*(D_i)$, where 
$i = 1,\ldots,r$. Moreover, the map
$\pr_Y^* : \Cl(Y) \to \Cl(Y \times V)$
is bijective, since $V$ is an affine space.
As a consequence, the map (\ref{eqn:pullback})
is surjective. By pulling back line bundles
on $X$ to $\{y \}\times Z$ and $Y \times \{z\}$, 
one checks that this map is injective.
\end{proof}

\begin{remark}\label{rem:can}
The classes of the prime divisors $D_1,\ldots,D_r$
form a \emph{canonical basis} of $\Pic(Z)$, 
which is stable under any automorphism of $Z$.
If $k$ is algebraically closed, then this basis
coincides with that constructed in 
\cite[Theorem 3.12]{Maccan1}.
\end{remark}

We now consider an \emph{arbitrary} field $k$. 
Note that the above assertion (ii) still holds 
in this setting, with the same proof. 
We will obtain versions of assertions (i) 
and (iii), for Picard groups with rational 
coefficients. For this, we need some preliminary 
results on the behavior of Picard groups 
under base change by field extensions. 

\begin{lemma}\label{lem:picinsep}
Let $X$ be a scheme of finite type, and 
$k \subset K$ an algebraic field extension. 

\begin{enumerate}

\item[{\rm (i)}]
For any line bundle $L$ on $X$ such that 
the line bundle $L_K$ on $X_K$ is effective, 
there exists an integer $n \geq 1$ such that 
$L^{\otimes n}$ is effective.

\item[{\rm (ii)}]    
The kernel of the pull-back map 
$\Pic(X) \to \Pic(X_K)$ is torsion.

\item[{\rm (iii)}]
If $k \subset K$ is purely inseparable,
then the cokernel of the above map is torsion.

\end{enumerate}

\end{lemma}

\begin{proof}
(i) Let $s$ be a non-zero global section
of $L_K$. Then there exists a finite subextension
$k \subset K' \subset K$ and a global
section $s'$ of $L_{K'}$ such that $s = s'_K$ 
(as follows e.g.~from \cite[\S 8.8]{EGAIV3}).
The norm $\rN_{X_{K'}/X}(s')$ is a non-zero
section of $\rN_{X_{K'}/X}(L_{K'})$. As the
latter line bundle is isomorphic to $L^{\otimes [K':k]}$,
this completes the proof.

(ii) The argument is similar:
reduce to a finite extension and use the norm.

(iii) Consider a line bundle $M$ on $X_K$. As 
above, we may choose a subextension 
$k \subset K' \subset K$ of finite degree $p^n$, 
and a line bundle $M'$ on $X_{K'}$ 
such that $M \simeq M'_K$. Since the projection
$X_{K'} \to X$ is a finite surjective
morphism, there exists an open 
affine covering $(U_i)_{i \in I}$ of $X$ such that
$M'$ is trivial on each pull-back $(U_i)_{K'}$. 
Thus, $M'$ is defined by a cocycle 
$(a_{ij})_{i,j \in I}$, where 
$a_{ij} \in 
\cO \big( (U_i \cap U_j)\otimes_k K' \big)^{\times}$.
As $(K')^{p^n} \subset k$, we see that 
$a_{ij}^{p^n} \in \cO(U_i \cap U_j)^{\times}$.
This means that $(M')^{\otimes p^n}$ lies in
the image of $\Pic(X) \to \Pic(X_{K'})$. Thus,
$M^{\otimes p^n}$ lies in the image of 
$\Pic(X) \to \Pic(X_K)$.
\end{proof}

Our next preliminary result is due to Sansuc;
we provide a proof for the reader's convenience.

\begin{lemma}\label{lem:picsep}
\cite[Lemma 6.3]{Sansuc}
Let $k \subset K$ be a Galois extension, 
$\Gamma_{K/k}$ its Galois group, and $X$ 
a complete variety. Then the pull-back map 
$\Pic(X) \to \Pic(X_K)^{\Gamma_{K/k}}$
is injective and its cokernel is torsion. 
If $X$ has a $k$-rational point, then this map
is bijective.
\end{lemma}

\begin{proof}
The projection $X_K \to X$ is a 
$\Gamma_{K/k}$-torsor, and hence yields a 
Hochschild-Serre spectral sequence relative to 
the \'etale cohomology of the multiplicative group
\[ H^i \big( \Gamma_{K/k}, H^j_{\et}(X_K, \bG_m) \big)
\Rightarrow H^{i+j}_{\et}(X,\bG_m). \]
The associated exact sequence of low-degree
terms begins with
\[ 0 \longrightarrow 
H^1 \big( \Gamma_{K/k}, \cO(X_K)^{\times} \big)
\longrightarrow \Pic(X)
\longrightarrow \Pic(X_K)^{\Gamma_{K/k}}
\longrightarrow 
H^2 \big( \Gamma_{K/k}, \cO(X_K)^{\times} \big). \]
Since $X$ is a complete variety,  
we have $\cO(X_K)^{\times} = K^{\times}$. 
Moreover,
$H^1(\Gamma_{K/k}, K^{\times}) = 0$
by Hilbert's Theorem 90, and 
the group $H^2(\Gamma_{K/k}, K^{\times})$
is torsion. This yields the first assertion.

To prove the second assertion, 
we recall the notion of rigidified line bundle
at a point $x \in X(k)$: this is a pair 
$(L,u)$, where $L$ is a line bundle on $X$ and 
$u : k \stackrel{\sim}{\longrightarrow} L(x)$
is a $k$-linear isomorphism. 
A morphism of rigidified line bundles 
$(L,u) \to (M,v)$ 
is a morphism of line bundles $f : L \to M$
such that $f(x) \circ u = v$. Note that 
every automorphism of a rigidified line 
bundle is trivial as
$\cO(X)^{\times} = k^{\times}$. Thus,
every rigidified line bundle is
determined by its isomorphism class
up to unique isomorphism. Moreover, these 
isomorphism classes form an abelian group
$\Pic(X,x)$, isomorphic to $\Pic(X)$
via the forgetful map. The Galois group 
$\Gamma_{K/k}$ acts naturally on the line bundles 
on $X_K$ rigidified at $x$, and hence on 
$\Pic(X_K,x)$; furthermore, if the class 
of a pair $(L,u)$ is invariant, then $(L,u)$ 
is equipped with a descent data in view of 
the above uniqueness property. By Galois descent, 
it follows that the pull-back map 
$\Pic(X,x) \to \Pic(X_K,x)^{\Gamma_{K/k}}$
is an isomorphism.
\end{proof}

Let us recall that we denote as $\Gamma_k$ the absolute Galois group of the base field $k$. Combining Lemmas \ref{lem:picinsep} and \ref{lem:picsep} 
(and still considering an arbitrary field $k$), we obtain:

\begin{lemma}\label{lem:piccomp}

\begin{enumerate}

\item[{\rm (i)}]    
For any complete variety $X$, the pull-back map 
$\Pic(X) \to \Pic(X^{\sep})^{\Gamma_k}$ 
induces isomorphisms
\[ \Pic(X)_{\bQ} \stackrel{\sim}{\longrightarrow}
\Pic(X^{\sep})_{\bQ}^{\Gamma_k}, \quad
\Eff(X) \stackrel{\sim}{\longrightarrow}
\Eff(X^{\sep})^{\Gamma_k}, \]
where $\Pic(X)_{\bQ} := \Pic(X) \otimes_{\bZ} \bQ$
and likewise for $\Pic(X^{\sep})_{\bQ}$,
and $\Eff(X) \subset \Pic(X)_{\bQ}$ denotes 
the effective cone.

\item[{\rm (ii)}]
For any variety $X$ and any purely inseparable 
field extension $k \subset K$, the pull-back
map $\Pic(X) \to \Pic(X_K)$ induces isomorphisms
\[ \Pic(X)_{\bQ} \stackrel{\sim}{\longrightarrow}
\Pic(X_K)_{\bQ}, \quad
\Eff(X)_{\bQ} \stackrel{\sim}{\longrightarrow}
\Eff(X_K)_{\bQ}. \]
\end{enumerate}
\end{lemma}

Together with Proposition \ref{prop:picrat} 
and Lemma \ref{lem:picinsep}, 
this readily yields:

\begin{proposition}\label{prop:picQ}
Let $X = Y \times Z$ be a complete homogeneous
variety as in Theorem \ref{thm:nonaffine}. 
 
\begin{enumerate}

\item[{\rm (i)}] 
The abelian group $\Pic(Z)$ is free of finite rank, 
and the pull-back map (\ref{eqn:pullback}) 
induces an isomorphism 
\[ \Pic(Y)_{\bQ} \times \Pic(Z)_{\bQ} 
\stackrel{\sim}{\longrightarrow} 
\Pic(X)_{\bQ}. \]

\item[{\rm (ii)}]
The cone $\Eff(Z)$ is generated by
a basis of $\Pic(Z)_{\bQ}$, and every effective
divisor on $Z$ is base-point-free.

\end{enumerate}

\end{proposition}

\begin{remark}
\label{rem:picQ}
More specifically, the action of $\Gamma_k$ on 
$\Pic(Z^{\sep})$ stabilizes the canonical basis 
$D_1, \ldots,D_r$. Furthermore, the classes 
\[ D_\Omega := \sum_{D_i \in \Omega} D_i, \]
where $\Omega$ runs over the $\Gamma_k$-orbits
in $\{ D_1,\ldots,D_r \}$, form a basis of
$\Pic(Z^{\sep})^{\Gamma_k}$ and generate the
monoid of base-point-free classes in that
group. In view of Lemma \ref{lem:piccomp}, 
it follows that the classes of the $D_{\Omega}$ 
form a basis of $\Pic(Z)_{\bQ}$ that generates
the cone of base-point-free classes.
\end{remark}

To complete the picture, we record the 
following observation on Picard groups
of abelian torsors, which is analogous
to a result of Sansuc on Picard groups
of torsors under reductive groups (see 
\cite[Lemme 6.7]{Sansuc}):

\begin{lemma}
\label{lem:abtors}
Let $A$ be an abelian variety, and $Y$ an
$A$-torsor. Then there is an isomorphism
$\Pic(Y)_{\bQ} \simeq \Pic(A)_{\bQ}$.
\end{lemma}

\begin{proof}
There exists a finite Galois extension
$k \subset K$ such that $Y$ has a $K$-rational
point $y_0$. This yields an isomorphism
\[ f: A_K \longrightarrow Y_K, \quad 
a \longmapsto a \cdot y_0 \]
and hence an isomorphism 
\[ 
f^* : \Pic(Y_K) \longrightarrow \Pic(A_K).
\]
The finite group $\Gamma_{K/k} := \Gal(K/k)$ acts 
on $A_K$ and $Y_K$; moreover, for the induced
actions on $\Pic(A_K)$ and $\Pic(Y_K)$, 
we have isomorphisms
\begin{equation}\label{eqn:isos}
\Pic(A)_{\bQ} 
\stackrel{\sim}{\longrightarrow}
\Pic(A_K)^{\Gamma_{K/k}}_{\bQ}, \quad \Pic(Y)_{\bQ} 
\stackrel{\sim}{\longrightarrow}
\Pic(Y_K)^{\Gamma_{K/k}}_{\bQ}
\end{equation}
in view of Lemma \ref{lem:picsep}. 
Note that $f$ is not equivariant for these 
actions, but is equivariant for the natural action 
on $Y_K$ and the ``dot action'' on $A_K$ defined by 
\[ 
\gamma \cdot a := 
( \tau_{y_0} \circ \gamma \circ \tau_{-y_0})(a)
= (\gamma \circ \tau_{\gamma(y_0) - y_0})(a),
\]
where $\tau_y$ denotes the translation by $y$.
Therefore, $f^*$ identifies $\Pic(Y_K)^{\Gamma_{K/k}}$
(where $\Gamma_{K/k}$ acts on $\Pic(A_K)$ via its
natural action on $A_K$)
and $\Pic(A_K)^{\cdot \Gamma_{K/k}}$ 
(where $\Gamma_{K/k}$ acts via the above
dot action on $A_K$). Moreover, the ``norm'' map
\[ \prod_{\gamma \in \Gamma_{K/k}} (\gamma \cdot)^* :
\Pic(A_K) \longrightarrow \Pic(A_K) \]
takes values in $\Pic(A_K)^{\cdot \Gamma_{K/k}}$
and its restriction to that subgroup is 
the multiplication by 
$ \vert \Gamma_{K/k} \vert = [K:k]$.
Thus, this map induces a surjection
$\Pic(A_K)_{\bQ} \to 
\Pic(A_K)^{\cdot \Gamma_{K/k}}_{\bQ}$.
Since 
$\prod_{\gamma \in \Gamma_{K/k}} (\gamma \cdot)^*
= (\prod_{\gamma \in \Gamma_{K/k}} \gamma^*)
\circ \tau_{y_1}$
for some $y_1 \in A(K)$, and $\tau_{y_1}$
is an automorphism, we see that 
$\prod_{\gamma \in \Gamma_{K/k}} (\gamma \cdot)^*$
and $\prod_{\gamma \in \Gamma_{K:k}} \gamma^*$
have the same image over the rationals.
This implies the assertion in view of the 
isomorphisms (\ref{eqn:isos}).
\end{proof}

\section{Parabolic varieties over a separably closed
field}
\label{sec:pvscf}

The results of Section \ref{sec:spic} reduce 
somehow the description of complete homogeneous 
$G$-varieties to the case where $G$ is semisimple.  
Over an \emph{algebraically closed} field $k$ 
of characteristic $p \geq 0$, these varieties are 
well-understood if $p = 0$; they have been classified 
by Wenzel if $p \geq 5$ \cite{Wenzel}, and by Maccan 
if $p = 2,3$  \cite{Maccan1,Maccan2}. 
Throughout this section, we assume that $k$ is
\emph{separably closed}, and show how to adapt Wenzel's 
and Maccan's results to this setting. We begin by 
fixing notation. From now on, we assume that $p>0$.

\subsection{Root system and Wenzel function}
\label{subsec:rswf}
Let $G$ be a semisimple $k$-group.
Then there exist a maximal $k$-subtorus 
$T \subset G$ and a Borel $k$-subgroup 
$B \supset T$; moreover, all such pairs
$(B,T)$ are conjugate in $G(k)$. 
We denote by $R$ the root system of $(G,T)$, 
by $R^+ \subset R$ the set of positive roots 
consisting of roots of $(B,T)$, and by 
$\Delta \subset R^+$ the corresponding basis of 
$R$. Every $\gamma \in R$ defines a root subgroup 
$U_{\gamma} \subset G$; recall that
$U_\gamma$ is normalized by $T$ 
and isomorphic to the additive group $\bG_a$, 
on which $T$ acts via $t \cdot x := \gamma(t) x$.

For any subset $I \subset \Delta$,
we denote by $P_I \subset G$
the $k$-subgroup generated by $B$ and the
$U_{-\alpha}$, where $\alpha \in I$.
The assignment $I \mapsto P_I$ yields a
bijection between subsets of $\Delta$
and smooth parabolic $k$-subgroups of 
$G$ containing $B$; moreover, every 
smooth parabolic $k$-subgroup of $G$ is 
 conjugate in $G(k)$ to a unique $P_I$. 
The notation 
$P^I := P_{\Delta \setminus I}$
will also be convenient; for instance,
the maximal smooth parabolic $k$-subgroups 
are exactly the conjugates in $G(k)$ of the parabolic
$P^{\alpha} := 
P_{\Delta \setminus \{ \alpha \} }$,
where $\alpha \in \Delta$.

\medskip

Next, let $X$ be a complete homogeneous variety
under $G$. Since $k$ is separably closed, $X$ has 
a $k$-rational point, i.e., it is a parabolic 
variety. Also, recall from Proposition
\ref{prop:rational} that $X$ has a unique 
$k$-rational $B$-fixed point $x$. Thus, we have 
a $G$-equivariant isomorphism
$G/P \stackrel{\sim}{\longrightarrow} X$,
where $P := \Stab_G(x)$ is a parabolic 
$k$-subgroup of $G$ containing $B$.
Moreover, this isomorphism is unique in view
of Proposition \ref{prop:authom}. So throughout 
the present section, \emph{we denote by $P$
a $k$-subgroup of $G$ containing the Borel 
subgroup $B$}. 

If $k$ is algebraically closed, then 
every parabolic $k$-subgroup of $G$ 
is conjugate in $G(k)$ to a $k$-subgroup $P$ 
as above (in fact, this holds over any perfect field by
Proposition \ref{prop:perfect}). But this fails
over every imperfect field, see the examples
in Section \ref{sec:example}.

\medskip

Denote by $Q$ the largest smooth $k$-subgroup of 
$P$, i.e., the closure of $P(k)$ (see 
\cite[Lemma C.4.1]{CGP}).
Then $Q \supset B$ and hence $Q = P_I$ for a unique
subset $I \subset \Delta$. Let $Q^- \subset G$
be the opposite parabolic $k$-subgroup to $Q$ that
contains $T$. The multiplication map of $G$
yields an open immersion $R_u(Q^-) \times Q \to G$,
and an isomorphism
\[ f : \prod_{\gamma \in R^+ \setminus R_I} U_{- \gamma}
\stackrel{\sim}{\longrightarrow} R_u(Q^-), \]
where the product is relative to an ordering
of $R^+ \setminus R_I$ as in \cite[p.~213]{Wenzel}.
We may now state a version of loc.~cit., 
Theorem 10 (i), (ii):

\begin{theorem}\label{thm:wenzel}
Let $U^-_P := R_u(Q^-) \cap P$.

\begin{enumerate}

\item[{\rm (i)}]
The above isomorphism yields an isomorphism
\[ \prod_{\gamma \in R^+ \setminus R_I} 
(U_{- \gamma} \cap P)
\stackrel{\sim}{\longrightarrow} U^-_P. \]

\item[{\rm (ii)}] 
For any $\gamma \in R^+ \setminus R_I$,
there exists a unique non-negative integer 
$m = m(\gamma)$ such that
$U_{- \gamma} \cap P \simeq \bfalpha_{p^m}$
(the $m$-th Frobenius kernel of $\bG_a$).

\item[{\rm (iii)}]
The multiplication map of $G$ yields
an isomorphism of schemes
$U^-_P \times Q 
\stackrel{\sim}{\longrightarrow} P$.

\item[{\rm (iv)}]
The $k$-group $U^-_P$ is infinitesimal, and
we have $Q = P_{\red}$.
\end{enumerate}

\end{theorem}

\begin{proof}
Since $Q \supset P$, the multiplication map of $G$ 
yields an open immersion 
\[ i : U^-_P \times Q \longrightarrow P. \]  
Moreover, $U^-_P$ is a $k$-subgroup of $R_u(Q^-)$, 
normalized by $T$, and $U^-_P(k)$ is trivial 
as $Q(k) = P(k)$.

We now adapt the arguments in \cite[p.~212--214]{Wenzel}
with the reduced subscheme $P_{\red}$ replaced by $Q$.
Given two linearly independent roots
$\alpha,\beta$, there exists $t \in T(k)$
such that $\alpha(t) = 1 \neq \beta(t)$
(indeed, the connected reduced kernel 
$S = \Ker(\alpha)^0_{\red}$ is a $k$-torus 
by \cite[IV.1.3.9]{DG}, and the quotient 
$\Ker(\alpha)/S$ is finite; thus, $\beta$ 
restricts to a nontrivial character of $S$.
Moreover, $S$ is split by loc.~cit., IV.1.1.7.
Since $k$ is infinite, it follows that 
$\beta(t) \neq 1$ for some $t \in S(k)$). 
As a consequence, Propositions 7 and 8 of 
\cite{Wenzel} hold with the same proofs. Now
(i) is a direct consequence of Proposition 8.

Next, note that each $k$-group 
$U_{- \gamma} \cap P$, where 
$\gamma \in R^+ \setminus R_I$, is normalized
by $T$ and has no non-trivial $k$-point.
Since $U_{- \gamma} \simeq \bG_a$ on which
$T$ acts linearly with weight $-\gamma$, it
follows that the ideal of $U_{- \gamma} \cap P$
in $\cO(\bG_a) = k[x]$ is generated by a positive
power of $x$; this implies (ii).

Together with (i), this implies in turn that
$U^-_P$ is infinitesimal; in particular, finite. 
As a consequence, the natural map $U^-_P \to P/Q$ 
is a closed immersion, and hence $i$ is a closed 
immersion as well. Since $P$ is connected and $i$ 
is open, this proves (iii).

As a further consequence, $(U^-_P)_{\red}$ is trivial. 
Together with the smoothness of $Q$, it follows
that $(U^-_P \times Q)_{\red} = Q$, and hence
$P_{\red} = Q$ in view of (iii). 
\end{proof}

As a direct consequence of the above theorem,
we obtain a very useful decomposition of twisted 
parabolic varieties over an \emph{arbitrary} field:

\begin{corollary}\label{cor:product}
Consider a semisimple adjoint group $G$ over a field $k$ 
and write $G = G_1 \times \cdots \times G_n$, where 
each $G_i$ is a simple adjoint $k$-group. Then every twisted 
parabolic $G$-variety $X$ admits a compatible
decomposition $X = X_1 \times \cdots \times X_n$,
where each $X_i$ is a twisted parabolic $G_i$-variety.
\end{corollary}

\begin{proof}
By induction on $n$, it suffices to show the
following claim: For any decomposition of $G$ 
into a product $G' \times G''$ of $k$-subgroups, 
the morphism 
\[ (f',f'') : X \longrightarrow X/G' \times X/G'' \]
is an isomorphism, where $f',f''$ are the categorical
quotients as in Lemma \ref{lem:quotient}. By fpqc descent,
it suffices in turn to prove the claim in the case where
$k$ is separably closed. Then $X \simeq G/P$ for a unique
parabolic $k$-subgroup $P \supset B$; moreover,
$B = (B \cap G') \times (B \cap G'')$. In view of Theorem
\ref{thm:wenzel}, it follows that 
$P = (P \cap G') \times (P \cap G'')$; this implies
our claim.
(Alternatively, we may reduce to the case where $k$ is 
algebraically closed and use \cite[Lemma 1.2]{Maccan2}).
\end{proof}

We now return to the setting of a \emph{separably closed}
field $k$.

\begin{corollary}\label{cor:definedoverfp}
Any parabolic $G$-variety is defined over $\bF_p$.
\end{corollary}

\begin{proof}
Let $P \supset B$ be the corresponding parabolic $k$-subgroup of $G$. By Theorem \ref{thm:wenzel}, from which we keep the notations, the multiplication map $i \colon   U^-_P \times Q \to P$ is an isomorphism, where $Q$ is a smooth parabolic and $U^-_P$ is a product of iterated Frobenius kernels of root subgroups. All of these objects are simultaneously defined over $\bF_p$, since $G$ is split; this concludes the proof.
\end{proof}

\begin{remark}\label{rem:wenzel}
Theorem \ref{thm:wenzel} extends to the setting 
of an arbitrary field $k$, a \emph{split} 
semisimple group $G$ with Borel subgroup 
$B$, and a parabolic subgroup $P \subset G$
containing $B$. Indeed, the proof works
unchanged in the case where $k$ is infinite. 
For a finite field $k$, we argue by Galois 
descent from the algebraic closure, using 
the fact that the formations of the largest 
smooth subgroup and the reduced subscheme 
commute with perfect field extensions.

Consequently, Corollary \ref{cor:definedoverfp} 
also extends to a split semisimple group $G$ 
over an arbitrary field $k$.
\end{remark}

Thanks to Corollary \ref{cor:product}, we may (and will) assume 
$G$ to be \emph{simple}, throughout the rest of the present section. 
Still adapting \cite{Wenzel}, we now define a map
\[ \varphi : R^+ \longrightarrow \bN \cup \{\infty\} \]
as follows: If $\gamma  \in R^+ \setminus R_I$, then
$\varphi(\gamma) := m(\gamma)$ with the notation of 
Theorem \ref{thm:wenzel} (ii). If $\gamma \in R_I ^+$,
we set $\varphi(\gamma) := \infty$. By this theorem,
we have 
\[ \varphi(\gamma) = 
\begin{cases} m & 
\text{ if }U_{-\gamma} \cap P \simeq \alpha_{p^m}, \\
    \infty & \text{ if }U_{-\gamma} \subset P. \\
\end{cases} \]
Moreover, the $k$-subgroup $P \subset G$ containing $B$ 
is uniquely determined by $\varphi$.

\begin{definition}\label{def:wenzel}
A \emph{Wenzel function} is a map 
$\varphi : R^+ \longrightarrow \bN \cup \{\infty\}$
associated to a parabolic $P$ via the above construction.
We then set $P := P_{\varphi}$.
The \emph{support} of $\varphi$ is defined to be 
\begin{align*}
\Supp(\varphi) = 
\{ \gamma \in \Delta ~\vert~ \varphi(\gamma) = \infty \}.
\end{align*}
\end{definition}

\noindent 
Note that  $(P_{\varphi})_{\red} = P_{\Supp(\varphi)}$. 
Moreover, for any two Wenzel functions $\varphi, \psi$, we have
\[ P_\varphi \subset P_\psi \, \Longleftrightarrow \, \varphi \leq \psi, \]
where the inequality on the right-hand side means that 
$\varphi(\gamma) \leq \psi(\gamma)$ for all positive roots 
$\gamma$. As a consequence,
\[ P_\varphi \cap P_\psi = 
P_{\min(\varphi,\psi)}. \]

\begin{remark}\label{rem:bij}
The assignment $\varphi \mapsto G/P_{\varphi}$
yields a bijective correspondence between
Wenzel functions and equivariant isomorphism classes
of parabolic $G$-varieties. Indeed, this follows from 
the discussion at the beginning of this section,
together with the classification of parabolic
$k$-subgroups of $G$ containing $B$ and the triviality 
of the equivariant automorphism groups of 
parabolic varieties (Proposition \ref{prop:authom}).
 
\end{remark}

\begin{definition}\label{def:bij}
We say that a projective homogeneous $G$-variety is 
\emph{of type $\varphi$} if it is equivariantly isomorphic 
to $G/P_{\varphi}$. By Remark \ref{rem:bij}, this 
$G$-variety is unique up to a unique equivariant isomorphism; 
we denote it by $X_{G,\varphi}$, or just by $X_{\varphi}$ if
this yields no confusion. 
\end{definition}

\begin{remark}\label{kvarieties}
    The bijection in Remark \ref{rem:bij} between equivariant isomorphism classes of projective homogeneous $G$-varieties $\{X_{\varphi}\}$ and Wenzel functions $\{\varphi\}$  is not true if we consider  the isomorphism classes of $\{X_{\varphi}\}$ as varieties over $k$ forgetting the $G$-action. For example, let $P$ be a $k$-subgroup containing
    $B$.  Then $G/P$ and $G/({}_mG \cdot P)$ are not isomorphic as $G$-varieties (since $G$ acts faithfully on the former but not on the latter). However, we have  \[
    G/({}_mG \cdot P)\simeq G^{(m)}/P^{(m)}\simeq (G/P)^{(m)}.\]  Since $G/P$ is defined over $\mathbb{F}_p$ (Corollary \ref{cor:definedoverfp}), it is isomorphic to its Frobenius twists  as a $k$-variety and therefore $G/P$ and $G/({}_mG \cdot P)$  are isomorphic as $k$-varieties.
\end{remark}

\subsection{Wenzel type}
\label{subsec:wenzeltype} 
For a root $\gamma$, let us denote as $\Supp(\gamma)$ its support, consisting of those $\alpha \in \Delta$ such that the $\alpha$-coordinate of $\gamma$ in the basis of simple roots is non-zero. Recall moreover that $G$ is said to be \emph{simply laced} if all its simple factors are of type $A,D,E$; equivalently, if its Dynkin diagram does not have any multiple edge.

By Theorem \ref{thm:wenzel}, the Wenzel function of $P$ 
equals that of the base change $\bar{P} = P_{\bar{k}}$.
Together with \cite[Theorem 10 (iii)]{Wenzel},
it follows that 
\[ \varphi(\gamma) =
\min \{ \varphi(\alpha), 
\alpha \in \Supp(\gamma) \} \]
for any $\gamma \in R^+$
under the additional assumptions that
$p > 3$ or $G$ is simply laced. 

\medskip

As a consequence, under these assumptions $\varphi$ is uniquely 
determined by its restriction to $\Delta$. 
On the other hand, we have
$_mG \cdot P^{\alpha} = P_{\varphi}$,
where $\varphi(\alpha) = m$ and
$\varphi(\beta) = \infty$ for any
$\beta \in \Delta \setminus \{ \alpha \}$.
Thus, every parabolic $k$-subgroup
$P \supset B$ satisfies
\begin{equation}\label{eqn:wenzel} 
P = \bigcap_{i = 1}^r 
{}_{m_i}G \cdot P^{\alpha_i},
\end{equation}
where $\alpha_1, \ldots, \alpha_r$
are distinct simple roots and 
$m_1, \ldots, m_r$ are non-negative 
integers; moreover, the pairs 
$(m_1,\alpha_1), \ldots,(m_r,\alpha_r)$ 
are unique up to reordering. We then 
have $P = P_\varphi$, 
where $\varphi(\alpha_i) = m_i$
for $i = 1,\ldots, r$, and
$\varphi(\beta) = \infty$
for any other simple root $\beta$.
As a consequence, the maximal smooth
parabolic subgroups containing $P_{\red}$
are exactly 
$P^{\alpha_1}, \ldots, P^{\alpha_r}$.

\begin{definition}
    We call a parabolic $k$-subgroup $P$ of the form
(\ref{eqn:wenzel}) 
\emph{of Wenzel type}.
\end{definition}

\begin{remark}
\label{rem:wenzeltype}
Throughout the text, we pay attention and adapt the proofs according to whether all the parabolic subgroups (those containing a Borel, over a separably closed field) are of Wenzel type or not. As already mentioned, this is the case when $p >3$ or when $G$ is simply laced. Moreover, this holds for $p =3$ when the group $G$ is not of type $G_2$. If $k$ is algebraically closed, this follows from \cite[Proposition 3.1]{Maccan1} together with \cite[Theorem 1]{Maccan2}. The analogous statement over a separably closed field then follows from Theorem \ref{thm:wenzel}: indeed, any parabolic subgroup containing a Borel is determined by its Wenzel function, which can be computed on the algebraic closure.
On the other hand, if $p=3$ and $G$ is of type $G_2$, or if $p=2$ and $G$ is not simply laced, then there are additional parabolic subgroups beyond the Wenzel type, as illustrated in Subsections \ref{sec:quasiwenzel} and \ref{subsec:G2p2} below.
\end{remark}

By a slight abuse of notation, given a parabolic of Wenzel type, we will identify the associated Wenzel function $\varphi$ with its restriction $\varphi_{\vert \Delta}$. Throughout the present subsection, we assume that 
\emph{all parabolic subgroups of $G^\sep$ containing $B$
are of Wenzel type}. 
Then the inclusion 
$P \subset {}_{m_i}G \cdot P^{\alpha_i}$
yields a $G$-equivariant morphism
\[ \pi_i : X \longrightarrow 
G/({}_{m_i}G \cdot P^{\alpha_i}) =: X'_i \]
for $i = 1, \ldots,r$. Moreover, 
$X'_i \simeq G^{(m_i)}/(P^{\alpha_i})^{(m_i)}$
where $(P^{\alpha_i})^{(m_i)}$ is a 
maximal smooth parabolic $k$-subgroup of the 
semisimple group $G^{(m_i)}$. 

By the results in \cite[\S 3.2]{Maccan1}
(which are obtained over an algebraically closed
field) together with Lemma \ref{lem:piccomp} 
and Proposition \ref{prop:picQ},
we have $\Pic(X'_i) \simeq \bZ$ and the pull-backs
of the ample generators of the $\Pic(X'_i)$
form a basis of $\Pic(X)_{\bQ}$ which generates
the effective cone. Moreover, the product morphism
\[ \prod_{i = 1}^r \pi_i : X 
\longrightarrow \prod_{i=1}^r X'_i \]
is an immersion with image being the closed $G$-orbit, 
which is unique due to Borel's fixed point theorem. 
Also, each $\pi_i$ is a \emph{contraction}, i.e., 
the natural map $\cO_{X'_i} \to (\pi_i)_*\cO_X$
is an isomorphism. We will need the following
generalization of the latter result:

\begin{lemma}\label{lem:contraction}
Assume that all parabolic $k$-subgroups of $G$ 
containing $B$ are of Wenzel type and let 
$P$ be such a $k$-subgroup, of the form
(\ref{eqn:wenzel}). Then the natural morphism
\[ 
\pi_I : X = G/P \longrightarrow 
G/\bigcap_{\alpha_i \in I} {}_{m_i}G \cdot P^{\alpha_i}
=: X'_I \]
is a contraction for any subset 
$I \subset \{ \alpha_1,\ldots, \alpha_r \}$.
Moreover, every contraction $\pi : X \to Y$ 
is the composition of a unique $\pi_I$ with
an isomorphism 
$X'_I \stackrel{\sim}{\longrightarrow} Y$. 
\end{lemma}

\begin{proof}
Consider the Stein factorization of $\pi_I$ as 
\[ X \stackrel{\varphi_I}{\longrightarrow} 
Y_I \stackrel{\psi_I}{\longrightarrow} 
G/\bigcap_{\alpha_i \in I} {}_{m_i}G \cdot P^{\alpha_i}, \]
where $\varphi_I$ is a contraction and $\psi_I$
is finite. By Blanchard's lemma (see e.g.
\cite[Theorem 7.2.1]{Brion1}), there is a unique action 
of $G$ on $Y_I$ such that $\varphi_I$ is equivariant.
Thus, $Y_I = G/Q_I$ for a $k$-subgroup $Q_I \subset G$ 
such that 
$P \subset Q_I \subset 
\bigcap_{\alpha_i \in I} {}_{m_i}G \cdot P^{\alpha_i}$
and $(Q_I)_{\red} = \bigcap_{\alpha_i \in I} P^{\alpha_i}$.
Therefore, 
$Q_I = \bigcap_{\alpha_i \in I} {}_{l_i}G \cdot P^{\alpha_i}$
for unique non-negative integers $l_i$.

Considering the complement 
$J := \{ \alpha_1,\ldots, \alpha_r \} \setminus I$
and the associated morphism $\pi_J$, we obtain
similarly a $k$-subgroup $Q_J \subset G$ such that
$P \subset Q_J \subset 
\bigcap_{\alpha_j \in J} {}_{m_j}G \cdot P^{\alpha_j}$
and 
$Q_J = \bigcap_{\alpha_j \in J} {}_{l_j}G \cdot P^{\alpha_j}$
for unique non-negative integers $l_j$.
As a consequence,
\[ P \subset Q_I \cap Q_J \subset 
\bigcap_{i =1}^r {}_{m_i}G \cdot P^{\alpha_i}, \]
and equality holds everywhere in view of the decomposition (\ref{eqn:wenzel}). By the uniqueness of this
decomposition, it follows that 
$Q_I = \bigcap_{\alpha_i \in I} {}_{m_i}G \cdot P^{\alpha_i}$.
Thus, $\pi_I = \varphi_I$ is a contraction.
This proves the first assertion. The second one
follows from Blanchard's lemma together with
the structure of parabolic $k$-subgroups containing $P$.
\end{proof}

\begin{example} 
Let us consider the simplest example, combinatorially speaking, of projective homogeneous varieties which are not isomorphic to any flag variety. 
Let $G= \SL_{n+1}$ and consider the action 
\[
G \circlearrowright \bP^n \times \mathrm{Grass}(2,n+1)^{(m)} 
= G/P^{\alpha_1} \times G^{(m)}/(P^{\alpha_2})^{(m)}.
\]
Next, consider the twisted incidence subvariety
\[
X = \{ (x,l),\, x \in \mathbb{P}^n, l \subset (\bP^n)^{(m)} \text{ line,} \,\, F^m(x) = (x_0^{p^m}\colon x_1^{p^m}\colon \ldots \colon x_n^{p^m}) \in l\}.
\]
Notice that this variety is preserved by the $G$-action on the product, since on the second factor of the product there is a Frobenius twist of the appropriate order. If we fix as base point 
\[
o = \big( x_o = (1\colon 0 \cdots \colon 0), l_o = (a_0\colon a_1\colon 0 \cdots \colon 0) \big),
\]
which corresponds to the standard choice of diagonal maximal torus and upper triangular Borel subgroup, by an explicit computation we get that 
\[
X = G/P, \qquad P = P^{\alpha_1} \cap ({}_mG \cdot P^{\alpha_2}).
\]
In this case, the two equivariant contractions are the natural projections from 
$X$ 
respectively to $\bP^n$ and to $\mathrm{Grass}(2,n+1)^{(m)}$.
For $n = 2$ we recover the hypersurface $X_m$ in (\ref{eqn:tiv}).
\end{example}


\subsection{Small characteristics: quasi-Wenzel type}
\label{sec:quasiwenzel}
The base field is still assumed to be separably closed, and the group $G$ is simple and simply connected. 

\begin{defn}
\label{def:edgehypothesis}
A simple and simply connected group satisfies the \emph{edge hypothesis} 
if its Dynkin diagram has an edge of multiplicity $p$.
\end{defn}

The definition above is equivalent to being simple of type $B_n$, $C_n$ or $F_4$ when $p=2$, or simple of type $G_2$ when $p=3$.

\begin{lemma}[\cite{CGP}, Lemma 7.1.2]
\label{lem: very special kernel}
If $G$ satisfies the edge hypothesis, there exist a unique $k$-subgroup
    \[
N := N_G \subset {}_1G
\]
which is non-central, normal, killed by Frobenius and minimal for these properties. We call the corresponding quotient $\rho_G \colon G \longrightarrow G^\flat$ the \emph{very special isogeny} of $G$.
\end{lemma}

\begin{corollary}
\label{lem: factorisation isogenies}
    Let $G$ be a simple, simply connected group over $k = k^\sep$. 
    Then any isogeny with source $G$ has a unique factorisation of the form
    \[
    F^m_{G^\flat} \circ \rho_G = \rho_{G^{(m)}} \circ F^m_G \quad \text{or} \quad F^m_G,
    \]
    followed by an isogeny with central kernel. 
\end{corollary}

\begin{proof}
    Over the algebraic closure, this follows from \cite[Proposition 1.12]{Maccan1}. The analogous statement over a separably closed field is then a consequence of the existence of the very special isogeny over $k$, proven in Lemma \ref{lem: very special kernel} above.
\end{proof}

Next, let us assume that either $p \geq 3$, or that $p=2$ and $G$ is \emph{not} of type $G_2$;
this includes both the cases for which all parabolics containing $B$ are of Wenzel type, and those satisfying the edge hypothesis. By Theorem \ref{thm:wenzel}, the Wenzel function of $P$ equals that of the base change $\bar{P} = P_{\bar{k}}$. This, together with \cite[Theorem 1.1, Corollary 1.3]{Maccan2} and Lemma \ref{lem: very special kernel}, yields the following. Every parabolic $k$-subgroup $P \supset B$ satisfies 
\begin{equation}
    \label{eqn:quasistandard}
    P = \bigcap_{i=1}^r (\ker \xi_i) \cdot P^{\alpha_i},
\end{equation}
where the $\xi_i$ are purely inseparable isogenies.  More precisely, 
the map $\xi_i$ is either an iterated Frobenius, or a very special isogeny (in the sense of Lemma \ref{lem: very special kernel}) followed by an iterated Frobenius, where the latter can only happen when the edge hypothesis is satisfied.

\begin{remark}
\label{rem:minimalkernels}
    Given a $k$-parabolic $P \supset B$, the decomposition (\ref{eqn:quasistandard}) is unique, if we assume the kernels $\ker \xi_i$ to be \emph{minimal} with respect to inclusion (among those satisfying the desired equality); see \cite[Remark 1.6]{Maccan2}.
\end{remark}

\begin{definition}
We call these parabolics of \emph{quasi-Wenzel} type.  When writing $P$ as such an intersection, the minimality condition will always be assumed. 
\end{definition}

In a compact form, $P$ is determined by the data 
\[
(\alpha_i, m_i,\epsilon_i), \quad \alpha_i \in \Delta, \, m_i \in \mathbb{N}, \, \epsilon_i \in \{0,1\}
\]
where $\epsilon_i =0$ means that $\ker \xi_i = {}_{m_i}G$ is a Frobenius kernel, while  $\epsilon_i = 1$ means that $\ker \xi_i = \rho^{-1} ({}_{m_i}G)$, where $\rho$ is the very special isogeny of $G$.

\begin{lemma}
\label{lem:supwenzelfunction}
    Let $\varphi$ and $\psi$ be two Wenzel functions and assume that $P_\varphi$ and $P_\psi$ are both of quasi-Wenzel type and that they do not generate the whole of $G$. Let $P=P_\theta$ be the $k$-subgroup generated by $P_\varphi$ and $P_\psi$. Then $\theta = \sup \{\varphi,\psi\}$.
\end{lemma}

\begin{proof}
    The $k$-subgroup $P$ contains both the intersections $U_{-\gamma} \cap P_\varphi$ and $U_{-\gamma}\cap P_\psi$ for any positive root $\gamma$; this means that $\sup \{\varphi,\psi\} \leq \theta$. So in order to conclude it suffices to show that $\sup \{\varphi,\psi\}$ is still a Wenzel function. This is true, thanks to the assumption that both parabolics we started with are of quasi-Wenzel type. Indeed, $\varphi$ is determined by $(m_\alpha,\epsilon_\alpha)_{\alpha \in \Delta}$ where $m_\alpha \in \bN \cup \{\infty\}$ and $\epsilon_\alpha \in \{0,1\}$, and analogously $\psi$ is determined by the values $(r_\alpha, \varepsilon_\alpha)_\alpha$. Then taking the sup of these two functions corresponds to the data $(\sup\{m_\alpha,r_\alpha\}, \nu_\alpha)_\alpha$, where $\nu_\alpha = \epsilon_\alpha$ if $m_\alpha>r_\alpha$, $\nu_\alpha = \varepsilon_\alpha$ if $m_\alpha <r_\alpha$, and $\nu_\alpha = \sup \{\epsilon_\alpha,\varepsilon_\alpha\}$ if $m_\alpha = r_\alpha$. This defines again a Wenzel function corresponding to a parabolic of quasi-Wenzel type.
\end{proof}

\subsection{Type $G_2$ in characteristic two}
\label{subsec:G2p2}
The base field is still assumed to be separably closed. To treat the only remaining case, we assume that $G$ is simple of type $G_2$ and that $p=2$. Denote by $\alpha_1$ the short simple root and by $\alpha_2$ the long one. Since the arguments in the proof of \cite[Corollary 2.24]{Maccan1}, carried out over $\bar{k}$, rely only on root systems and Lie subalgebras, and since all the relevant objects are defined over $k$, we obtain the following statement. There are exactly two restricted Lie subalgebras of $\Lie G$ that strictly contain $\Lie P^{\alpha_1}$, namely:
\[
\mathfrak{q}_1 := \Lie P^{\alpha_1} \oplus \mathfrak{g}_{-2\alpha_1-\alpha_2}, \quad \mathfrak{q}_2 := \Lie P^{\alpha_1}\oplus \mathfrak{g}_{-\alpha_1}\oplus \mathfrak{g}_{-\alpha_1-\alpha_2}.
\]

For the notion and examples of restricted Lie algebras we refer to \cite[Definition 10.39]{Milne}.
Next, lifting these two Lie subalgebras to $k$-subgroups of height one -- namely, whose Frobenius is trivial, see \cite[II, $\S$7, n.4]{DG} -- and then multiplying them by the smooth parabolic $k$-subgroup $P^{\alpha_1}$, one obtains two \emph{exotic} parabolics containing $B$, which we call $\mathcal{Q}_1$ and $\mathcal{Q}_2$, which satisfy 
\[
\Lie \mathcal{Q}_1 = \mathfrak{q}_1, \quad \Lie \mathcal{Q}_2 = \mathfrak{q}_2 \quad \text{and} \quad  (\mathcal{Q}_1)_{\text{red}} = (\mathcal{Q}_2)_{\text{red}} = P^{\alpha_1}\]
and whose unipotent infinitesimal part has height one. In terms of Wenzel functions, they satisfy respectively $\mathcal Q_1 = P_{\varphi_1}$, where
\[
\varphi_1(2\alpha_1+\alpha_2) = 1, \quad \varphi_1(\alpha_1) = \varphi_1(\alpha_1+\alpha_2) = \varphi_1(3\alpha_1+\alpha_2) = \varphi_1(3\alpha_1+2\alpha_2) = 0
\]
and $\mathcal Q_2 = P_{\varphi_2}$, where
\[
\varphi_2(\alpha_1) = \varphi_2(\alpha_1+\alpha_2) = 1, \quad \varphi_2(2\alpha_1+\alpha_2) = \varphi_2(3\alpha_1+\alpha_2) = \varphi_2(3\alpha_1+2\alpha_2) = 0.
\]
These two parabolic $k$-subgroups cannot be written in the form (\ref{eqn:quasistandard}); moreover, they are enough to complete the classification in characteristic two (over a separably closed field). Indeed, this follows from \cite[Section 2.6]{Maccan1}, together with the fact that the $\bar{k}$-parabolic of this reference are defined over $k$.
Any parabolic $k$-subgroup $P\supset B$ is obtained by intersection of ${}_{m_2}G \cdot P^{\alpha_2}$, for some $m_2 \in \bN \cup \{\infty\}$, together with one among
\[
{}_{m_1}G \cdot P^{\alpha_1}, \quad (F^{m_1})^{-1}(\mathcal{Q}_1), \quad (F^{m_1})^{-1}(\mathcal{Q}_2),
\]
for some $m_1 \in \bN \cup \{\infty\}$. For more details on the construction of $\mathcal{Q}_i$ (over the algebraic closure) and the corresponding homogeneous varieties, see \cite[Section 2.6]{Maccan1}.

\begin{example}
\label{ex:supwenzelfunction}
    For these exotic parabolics, the analogous statement as the one in Lemma \ref{lem:supwenzelfunction} fails. Indeed, the smallest parabolic subgroup containing $\mathcal Q_1$ and $\mathcal Q_2$ is 
    $P_\theta = {}_1G \cdot P^{\alpha_1}$. Thus
    \[
    \theta (3\alpha_1+2\alpha_2) = 1> \sup \{\varphi_1(3\alpha_1+2\alpha_2), \varphi_2(3\alpha_1+2\alpha_2)\}  =0
    \]
    which is a strict inequality.
\end{example}

\subsection{Some characteristic-free statements} Let $k = k^\sep$ of arbitrary characteristic. For a parabolic $P \supset B$, let $Q^{\alpha_i}$ denote the smallest 
parabolic $k$-subgroup containing $P$ and $P^{\alpha_i}$; if $p\geq 5$, then this recovers the Wenzel type and $Q^{\alpha_i} = {}_{m_i}G \cdot P^{\alpha_i}$ as in (\ref{eqn:wenzel}); if the edge hypothesis is satisfied, then $Q^{\alpha_i} = (\ker \xi_i) \cdot P^{\alpha_i}$ as in (\ref{eqn:quasistandard}). The inclusion $P \subset Q^{\alpha_i}$ yields a $G$-equivariant morphism
\begin{align}
\label{maximal contraction}
\pi_i \colon X \longrightarrow G/Q^{\alpha_i} =: X_i^\prime
\end{align}
for $i = 1,\ldots,r$. As in the Wenzel type case, each 
$X_i^\prime$ is projective homogeneous of Picard rank one, and the pull-backs of the ample generators of the $\Pic(X_i^\prime)$ form the canonical basis of $\Pic (X)$. Moreover, each $\pi_i$ is a contraction, and the product morphism $\pi_1\times \ldots \times \pi_r$ is an immersion with image being the closed $G$-orbit. This amounts to saying that $P$ is recovered as intersection of the $Q^{\alpha_i}$; indeed, the analogous statement is true over the algebraic closure by \cite[Theorem 1]{Maccan2}, and it stays true over $k$ because the Wenzel function of any $P \supset B$ can be computed after base change to $\bar{k}$. Then we get the following generalization of Lemma \ref{lem:contraction}.

\begin{lemma}
    \label{lem:contraction any char}
    Let $P$ be a parabolic $k$-subgroup containing $B$ and let $Q^{\alpha_i}$ be the $k$-subgroup generated by $P$ and $P^{\alpha_i}$. Then the natural morphism 
    \[\pi_I : X = G/P \longrightarrow 
G/\bigcap_{\alpha_i \in I} Q^{\alpha_i} =: X_I^\prime
\]
is a contraction for any subset 
$I \subset \{ \alpha_1,\ldots, \alpha_r \}$.
Moreover, every contraction $\pi : X \to Y$ 
is the composition of a unique $\pi_I$ with
an isomorphism 
$X_I^\prime \stackrel{\sim}{\longrightarrow} Y$. 
\end{lemma}

\begin{proof}
    We follow exactly the same arguments as the one of Lemma \ref{lem:contraction}, by replacing ${}_{m_i}G \cdot P^{\alpha_i}$ with $Q^{\alpha_i}$. The proof works also in this slightly more general setting, because it relies on Blanchard's lemma, together with the uniqueness of the decomposition $P = \bigcap_{i=1}^r Q^{\alpha_i}$. 
\end{proof}

At this stage, a natural question is to ask whether, for a fixed $G$, there is a natural characterization in terms of Wenzel functions to determine whether $G/P_\varphi$ and $G/P_\psi$ are isomorphic as $k$-varieties. Clearly, $G/({}_mG\cdot P) \simeq G/P$ as illustrated in Remark \ref{kvarieties}; but it is not clear whether these are the only cases where such phenomena happen. We may add the assumption that $G$ acts faithfully, since we are interested only in the corresponding underlying varieties.

\begin{proposition}
\label{prop:characterisationvarphi}
    Let $G$ be a simple adjoint group over $k = k^\sep$ and let $X$ be a parabolic variety on which $G$ acts faithfully. Then there is a unique (up to automorphisms of the Dynkin diagram of $G$) Wenzel function $\varphi$ such that $X = G/P_\varphi$ as $k$-varieties.
\end{proposition}

\begin{proof}
Let us write $X=G/P$ for some parabolic $k$-subgroup $P$ containing a Borel subgroup; it suffices to show that one can recover the function $\varphi$ satisfying $P=P_\varphi$, starting from geometric information on $X$. The Wenzel function is invariant by base change to the algebraic closure, thus we may assume that $k$ is algebraically closed. 
Let us start under the assumption that $G$ is not of type $G_2$ when $p=2$, so that all parabolic subgroups are of quasi-Wenzel type. 
Let $r$ be the Picard rank of $X$, and consider for $1\leq i \leq  r$ the target $X_i^\prime$ 
of the contraction $\pi_i$ having maximal relative Picard rank, as defined in (\ref{maximal contraction}). Each $X_i^\prime$ is a parabolic variety of Picard rank one having smooth stabilizer, either homogeneous under $G$ or under its quotient by the very special isogeny $G^\flat$. Thus, up to a Dynkin diagram automorphism, there are simple roots $\beta_1,\ldots, \beta_r$ of $G$ such that $X_i^\prime \simeq G/P^{\beta_i}$ or (when the edge hypothesis is satisfied) $X_i^\prime \simeq G^\flat/P^{\beta_i^\flat}$, where $\beta_i^\flat$ denotes the image of $\beta_i$ by the very special isogeny $\rho_G$; see \cite[Theorem 2.1]{Maccan1} for a more detailed statement. Now, we have to compute some intersection numbers to take account of the Frobenius kernels: let $C_i \subset X$ be a rational curve (unique up to rational equivalence) not contracted by $\pi_i$, and let $S_i$ be the class of the Schubert curve on $X_i^\prime$; moreover, let $(D_1,\ldots,D_r)$ be the dual basis of the Picard group of $X$. Then there is a unique integer $m_i \geq 0$ such that 
\[
    (\pi_i)_* D_i \cdot S_i = D_i \cdot \pi_i^* S_i = p^{m_i} D_i \cdot C_i = p^{m_i}.
\]
With this information, we recover $\varphi$ as follows: $\varphi(\beta) = \infty$ for any $\beta \notin \{\beta_i,\ldots,\beta_r\}$, and then $\varphi$ corresponds to the data $ (\beta_i,m_i,\epsilon_i)_{i=1\ldots r}$,
where $m_i$ is given above and 
\[
\epsilon_i = \begin{cases}
        0 & \text{if } X_i^\prime = G/P^{\beta_i};\\
        1 & \text{if } X_i^\prime = G^\flat/P^{\beta_i^\flat}.
\end{cases}
\]
In the form (\ref{eqn:quasistandard}), we can write 
\[
P_\varphi = \bigcap_{i=1}^r (\ker \big( \rho_{G^{(m_i)}}^{\epsilon_i} \circ F_G^{m_i}) \big) 
\cdot P^{\beta_i},
\]
where this expression is unique thanks to the assumption that $G$ acts faithfully. It remains only to consider $p=2$ when $G$ is of type $G_2$: the case of Picard rank one is treated in \cite[Theorem 2.17]{Maccan1}, and the case of Picard rank two in \cite[Corollary 3.11]{Maccan2}.
\end{proof}

Notice that the uniqueness statement in Proposition \ref{prop:characterisationvarphi} fails if the group $G$ is not fixed: for instance, the projective space $\mathbb P^{2n-1}$ writes both as a parabolic variety under $\PGL_{2n}$ and under $\PSp_{2n}$. The above result yields the following characterization: up to outer automorphisms, and up to pulling back by an appropriate iterated Frobenius (when the $G$-action is not necessarily assumed to be faithful) the Wenzel function is uniquely determined by the isomorphism class of $X$ \emph{as a $k$-variety.}

\begin{corollary}
    Let $G$ be a simple adjoint group over $k = k^\sep$ and let $X$ be a parabolic $G$-variety. Suppose $\varphi$ and $\psi$ are Wenzel functions such that $G/P_\varphi$ and $G/P_\psi$ are isomorphic as $k$-varieties. Then, up to relabeling the simple roots by an automorphism of the Dynkin diagram, there is a unique integer $m \in \mathbb Z$ such that $\varphi(\gamma) = \psi(\gamma) + m$ for all positive roots $\gamma$. 
\end{corollary}

\begin{example}
\label{ex infinite non isomorphic}
    Let us illustrate Proposition \ref{prop:characterisationvarphi} in a more concrete way, by exhibiting a first infinite family of parabolic varieties of Picard rank two which are not isomorphic to any flag variety (with reduced stabilizer). We choose such a family to be of the smallest dimension possible for each type of $G$; notice that this is equivalent to asking which is the smallest dimension for which we have parabolic varieties that are \emph{not} Frobenius split \cite{Lauritzen}, and in terms of Wenzel functions it means a Wenzel function that is non-constant on the simple roots on which it is not infinite. If $G$ is of type $G_2$ and $p=2$, there is a unique such variety of Picard rank one, which has dimension $5$; see \cite[Theorem 2.17]{Maccan1}. Otherwise, we must consider Picard rank two and a parabolic of the form
    \[
    P_m := P^{\alpha_1} \cap {}_m G \cdot P^{\alpha_n},
    \] 
    where the notation on simple roots is the one of Bourbaki. The corresponding parabolic varieties are pairwise not isomorphic to each other. Moreover, the choice of roots $\alpha_1,\alpha_n$ is such that
    \[
    \dim (G/P_m) = \# \{ \gamma \in R^+ \colon \alpha_1 \in \mathrm{Supp}(\gamma) \text{ or } \alpha_n \in \mathrm{Supp}(\gamma)\}
    \]
    is minimal (among all the possible dimensions in Picard rank two). The corresponding dimensions are  $2n-1$ in type $A_{n}$, $(n^2+3n+2)/2$ in types $B_n$ and $C_n$, $2(n-1)$ in type $D_n$, and finally 24, 43, 84, 20, 6 (types $E_6$, $E_7$, $E_8$, $F_4$ and $G_2$ when $p \geq 3$ respectively). In particular, in type $A_n$ the corresponding varieties are the hypersurfaces $X_m$ introduced in (\ref{eqn:tiv}).
\end{example}

\begin{corollary}
\label{cor:infinite non isomorphic}
   Over a separably closed field $k$, for any fixed simple adjoint group $G$, the set of isomorphism classes (as $k$-varieties) of $G/P_\varphi$, where $\varphi$ varies through all Wenzel functions, is infinite.
\end{corollary}

\section{Parabolic varieties over an arbitrary field}
\label{sec:pvaf}

The ground field $k$ is arbitrary throughout this
section. Let $G$ be a semisimple $k$-group and let
$X$ be a parabolic $G$-variety, that is,
a complete homogeneous $G$-variety having 
a $k$-rational point $x$. Then $P := \Stab_G(x)$ 
is a parabolic $k$-subgroup of $G$, and the 
orbit map $g \mapsto g \cdot x$ factors
through a $G$-equivariant isomorphism
$G/P \stackrel{\sim}{\longrightarrow} X$.
Together with Remark \ref{rem:norm}, 
it follows that we may view $X$ as the
variety of conjugates of $P$. To obtain 
further information on $X$, we recall 
some notions and results from \cite{BT}
and \cite{Ti65}.

\subsection{Minimal parabolics and roots over the base field}
\label{subsec:mpr}

We first choose a maximal split $k$-torus 
$S \subset G$. By \cite[Th\'eor\`eme 4.21]{BT},
any two such subgroups are conjugate in $G(k)$. 
We then choose a minimal smooth parabolic $k$-subgroup 
$P_{\min} \subset G$ containing $S$. By loc. cit.,
Corollaire 5.9, Th\'eor\`eme 5.15, 
any two such subgroups are conjugate in 
$N_G(S)(k)$. Thus, all such pairs $(S,P_{\min})$
are conjugate in $G(k)$. Also, recall 
that $L := C_G(S)$ (the centralizer of $S$ in $G$) 
is a Levi subgroup of $P_{\min}$, 
and the corresponding opposite parabolic
$P_{\min}^-$ is conjugate to $P_{\min}$
in $G(k)$ (loc.~cit., Th\'eor\`eme 4.13,
Corollaire 4.14, Corollaire 4.16).
The derived subgroup of $C_G(S)$ is the (semisimple)
\emph{anisotropic kernel}, see \cite[\S 2.2]{Ti65}.

Next, we choose a maximal $k$-torus 
$T \subset P_{\min}$ containing $S$.
Then $T^{\sep}$ is a maximal torus 
of the parabolic subgroup 
$P_{\min}^{\sep} \subset G^{\sep}$, 
stable by the action of the absolute Galois group
$\Gamma_k := \Gal(k^{\sep}/k)$. So we may choose 
in addition a Borel $k^{\sep}$-subgroup 
$B \subset G^{\sep}$ such that 
$T^{\sep} \subset B \subset (P_{\min})^{\sep}$ 
(there may be no $\Gamma_k$-stable such $B$).

With the notation of Section \ref{sec:pvscf}, we have
\[ (P_{\min})^{\sep} = P_{\Delta_{\min}} \]
for a unique subset $\Delta_{\min} \subset \Delta$
that may be identified with the set of 
simple roots of the anisotropic kernel. 
(Let us mention that in the original work \cite{Ti65} 
on the subject, such a subset is denoted as $\Delta_0$). 
Thus, the maximal smooth parabolic subgroups of $G^{\sep}$ 
containing $(P_{\min})^{\sep}$ are exactly the $P^\alpha$, 
where $\alpha \in \Delta \setminus \Delta_{\min}$;
they are permuted by the absolute Galois group $\Gamma_k$. 

We now recall the action of $\Gamma_k$ on the 
set of simple roots known as the $*$-\emph{action} 
(see \cite[\S 2.3]{Ti65}, \cite[pp.~527--528]{MPW}).

\begin{definition}\label{def:*action}
Given any $\sigma \in \Gamma_k$, there exists a 
$n_{\sigma} \in N_{G^{\sep}}(T^{\sep})$ such that
$\Int(n_{\sigma}) \circ  \sigma$ stabilizes
the pair $(T^{\sep},B)$ (where $\Int(n_{\sigma})$ 
denotes conjugation by $n_{\sigma}$), and the
coset 
$w_{\sigma} := n_{\sigma} T^{\sep} \in W^{\sep}$
is uniquely determined
(where $W := N_G(T)$ denotes the Weyl group). Thus, 
$\Int(n_{\sigma}) \circ \sigma$ induces
a permutation $\sigma^*$ of the root system $R$
which stabilizes the subset $R^+$ of positive roots, 
and hence the basis $\Delta$. Since we have 
$\sigma^*(\gamma_1 + \gamma_2) = 
\sigma^*(\gamma_1) + \sigma^*(\gamma_2)$ 
for any $\gamma_1,\gamma_2 \in R$ such
that $\gamma_1 + \gamma_2 \in R$, we see
that $\sigma^*$ is uniquely defined by its
restriction to $\Delta$, given by 
$\alpha \mapsto w_{\sigma} \sigma (\alpha)$. 

The permutation $\sigma^*$ is defined, more generally, 
for any maximal $k$-torus $T \subset G$ (not necessarily
containing a maximal split $k$-torus) and any Borel 
$k^{\sep}$-subgroup $B \subset G^{\sep}$ 
containing $T^{\sep}$. By \cite[p.~527--528]{MPW},
the assignment $\sigma \mapsto \sigma^*$
yields a $\Gamma_k$-action on $\Delta$,
which is independent of the choice of the pair
$(T,B)$.

If the $*$-action on $\Delta$ is trivial, 
we say that $G$ is of \emph{inner type}; otherwise 
$G$ is said to be of \emph{outer type}.
\end{definition}

\begin{proposition}\label{prop:*action}
The $*$-action stabilizes $\Delta_{\min}$ and
induces the $\Gamma_k$-action on 
$\Delta \setminus \Delta_{\min}$,
where the latter is viewed as the set
of maximal smooth parabolic subgroups
containing $(P_{\min})^{\sep}$. 
\end{proposition}

\begin{proof}
Let $\sigma \in \Gamma_k$; then
$\sigma$ stabilizes both $T^{\sep}$
and $(P_{\min})^{\sep}$, and hence it
stabilizes the Levi subgroup of
$(P_{\min})^{\sep}$ containing $T^{\sep}$.
As a consequence, $\sigma$ stabilizes
$\Delta_{\min}$. Let $w$ be the unique 
element of $W^{\sep}$ such that 
$w\sigma(\Delta) = \Delta$.
Then $B$ and $w^{-1} B w = \sigma(B)$
are two Borel subgroups of $G^{\sep}$
containing $T^{\sep}$ and contained in 
$(P_{\min})^{\sep}$. It follows that
$w \in (W_{\Delta_{\min}})^{\sep}$
with an obvious notation. So 
$w \sigma$ stabilizes 
$\Delta \cap R_{\Delta_{\min}} 
= \Delta_{\min}$.

Next, let 
$\alpha \in \Delta \setminus \Delta_{\min}$.
Then $\sigma(P^{\alpha}) =P^{\beta}$
for a unique $\beta \in 
\Delta \setminus \Delta_{\min}$. 
Since $w\sigma$ stabilizes
$\Delta \setminus \Delta_{\min}$, 
it follows that
$w(\beta) = w\sigma(\alpha) 
\in \Delta \setminus \Delta_{\min}$.
But $w(\beta) - \beta$ is a linear
combination of roots in $\Delta_{\min}$,
since $w \in W_{\Delta_{\min}}^{\sep}$.
By linear independence of the simple
roots, this forces $w(\beta) = \beta$
and $w\sigma(P^{\alpha}) = P^{\beta}$.
Thus, $\sigma^*(\alpha) = \beta$;
this completes the proof.
\end{proof}

\subsection{Wenzel type}
\label{subsec:wenzel}

Until the end of this section, we consider $k^{\sep}$-subgroups 
$P \subset G^{\sep}$ 
\emph{that contain the chosen Borel subgroup $B$},
and we describe combinatorially (in terms of the simple roots 
and the Wenzel function) when $P$ is defined over $k$. 
One may also ask when the twisted parabolic variety 
$G^{\sep}/P$ descends over $k$; a combinatorial criterion 
for this will be obtained in the next section
(Proposition \ref{prop:definedfunction}).

\begin{proposition}\label{prop:descent}
Let $P \subset G^{\sep}$ be a parabolic 
$k^{\sep}$-subgroup of the form (\ref{eqn:wenzel}),
that is,
\[ P = {_{m_1}}G \cdot P^{\alpha_1} \cap \cdots
\cap {_{m_r}}G \cdot P^{\alpha_r}. \] 
Then $P$ is defined over $k$ if and only if: 
$\{ \alpha_1, \ldots, \alpha_r \}$
is contained in 
$\Delta \setminus \Delta_{\min}$ and
is stable under the $*$-action, and 
the $m_i$ are invariant under that action.
\end{proposition}

\begin{proof}
Assume that $P$ is defined over $k$,
i.e., stable under the action of $\Gamma_k$. 
Then $\Gamma_k$ permutes the 
${_{m_i}} G \cdot P^{\alpha_i}$
by the uniqueness in (\ref{eqn:wenzel}),
and hence permutes their largest smooth 
$k^{\sep}$-subgroups $P^{\alpha_i}$.
In particular, $\Gamma_k$ stabilizes
$Q := P^{\alpha_1} \cap \cdots \cap P^{\alpha_r}$,
the largest smooth $k^{\sep}$-subgroup 
of $P$. Thus, $\Gamma_k$ stabilizes
the centralizer $C_Q(S^{\sep})$.
The latter is a smooth parabolic 
$k^{\sep}$-subgroup of 
$C_{G^{\sep}}(S^{\sep}) = C_G(S)^{\sep}$.
Since $C_G(S)$ contains no proper 
parabolic $k$-subgroup, it follows
that $Q$ contains $C_G(S)^{\sep}$.
As a consequence, 
$\{ \alpha_1,\ldots,\alpha_r \}$
is contained in 
$\Delta \setminus \Delta_{\min}$
and is stable under the action of
$\Gamma_k$; moreover, the $m_i$ are 
invariant under that action.
This proves the ``only if'' implication;
the ``if'' implication follows from Galois 
descent.
\end{proof}

As a direct consequence, we get back the following
well-known result (see \cite[2.5.4]{Ti65}):

\begin{corollary}
\label{cor:smoothdescent}
Consider a \emph{smooth} parabolic $k^\sep$-subgroup
$P \subset G^{\sep}$ containing $B$,
\[
P = \bigcap_{i=1}^r P^{\alpha_i}.
\]
Then $P$ is defined over $k$ if and only if: 
$\{ \alpha_1, \ldots, \alpha_r \}$
is contained in 
$\Delta \setminus \Delta_{\min}$ and
is stable under the $*$-action.
\end{corollary}

Proposition \ref{prop:descent} may be reformulated 
in terms of the $*$-action on Wenzel functions, 
that we define as follows. 

\begin{definition}\label{def:*varphi}
Let $P \subset G^{\sep}$ be a parabolic 
$k^{\sep}$-subgroup containing $B$, and
$\varphi$ the corresponding Wenzel function. 
Given $\sigma \in \Gamma_k$, 
recall that there is a 
$n_{\sigma} \in N_{G^{\sep}}(T^{\sep})$
such that 
$B = n_{\sigma} \sigma(B) n_{\sigma}^{-1}$;
moreover, the coset 
$w_{\sigma} := n_{\sigma} T^{\sep}$
is unique. So
$n_{\sigma} \sigma(P) n_{\sigma}^{-1}$
is a parabolic $k^{\sep}$-subgroup of $G^{\sep}$
containing $B$. We denote the Wenzel function 
of this parabolic by $\sigma^* \varphi$, so that
\[ (\sigma^*\varphi)(\gamma) =
\varphi\big( w_{\sigma} \sigma (\gamma) \big) =
\varphi\big( \sigma^*(\gamma) \big)\]
for all $\gamma \in R^+$.
We say that $\varphi$ is 
\emph{invariant under the $*$-action} 
if $\sigma^*\varphi = \varphi$ for all 
$\sigma \in \Gamma_k$.
\end{definition}

\begin{corollary}\label{cor:wenzel}
The parabolic $k^{\sep}$-subgroup 
$P \subset G^{\sep}$ as in (\ref{eqn:wenzel})
is defined over $k$ if and only if its Wenzel 
function $\varphi$ satisfies:
$\varphi(\alpha) = \infty$ for all 
$\alpha \in \Delta_{\min}$, and 
$\varphi$ is invariant under the $*$-action.
\end{corollary}

\begin{proof}
This follows from Proposition \ref{prop:descent}
together with the fact that $\varphi$ is uniquely
determined by its restriction to $\Delta$
(Subsection \ref{subsec:wenzeltype}).
\end{proof}

\begin{example}
Let us illustrate Proposition \ref{prop:descent} on the simplest example with a non-trivial two-orbit Galois action, that is in type $A_3$. Assume that $p \geq 5$ and consider a 
quadratic extension $k^\prime/k$; denote the non-trivial involution in the corresponding Galois group as $\sigma$. Let $h$ be a non-degenerate $\sigma$-hermitian form on the $k^\prime$-module $V = (k^\prime)^4$, and denote as $\tau$ the associated adjoint involution on $A= M_4(k^\prime)$. Attached to $(A,\tau)$ there is an outer $k$-form 
$G := \SU(A,\tau)$ of 
$G^{\sep} = \SL_{4,k^{\sep}}$. Under the standard ordering of simple roots $(\alpha_1,\alpha_2,\alpha_3)$ in type $A_3$, the Galois $*$-action induced by $\tau$ is the involution exchanging $\alpha_1$ and $\alpha_3$, and fixing $\alpha_2$. Next, consider the subgroup scheme
\[
Q := ({}_1 G^{\sep} \cdot P^{\alpha_1}) \cap P^{\alpha_2} 
\cap ({}_1 G^{\sep} \cdot P^{\alpha_3}) \subset G^{\sep}.
\]
By Proposition \ref{prop:descent}, there is a unique parabolic $k$-subgroup $P \subset G$ such that $P^{\sep} = Q$; in particular, the non-trivial Galois orbit is 
$\Omega = \{\alpha_1,\alpha_3\}$ with  $\varphi (\alpha_i) = 1$ for $\alpha_i \in \Omega$,
while the middle simple root has no Frobenius fattening. Geometrically, $G^{\sep}/Q$ is the twisted full flag variety parameterizing triples of vector subspaces $(V_1,V_2,V_3)$ satisfying $\dim_{k^{\sep}} V_i=i$ and the incidence relations $V_1 \subset F(V_2) \subset V_3$. The $k$-form $G/P$ parametrizes the full flags in the $k^\prime$-vector space $V$ satisfying the above incidence relations, and moreover that are compatible with the $(A,\tau)$-descent. 
\end{example}

\subsection{Small characteristics} Let $p=2$ or $p=3$; we now get analogous structure results, under the assumption that not all parabolic subgroups over $k^\sep$ are of Wenzel type. Two cases have to be treated separately: the one of quasi-Wenzel type, taking into account the presence of the very special isogeny, and finally the one in which the group $G$ has some simple factors of type $G_2$ in characteristic $p=2$.

\begin{remark}
\label{rem:subdiagrams}
Since the $*$-action is defined as being an action by automorphisms of the Dynkin diagram, it preserves the Cartan matrix, hence adjacency relations between simple roots. In particular, it sends connected subdiagrams to isomorphic connected subdiagrams, and it preserves the lengths of the roots (in the case where the group is not simply laced).  
\end{remark}

\begin{lemma}
\label{lem: kernels descend}
The kernels in (\ref{eqn:quasistandard}) are defined over $k$.
\end{lemma}

\begin{proof}
    First, the subgroup $N_{G^{k^\sep}}$ descends over $k$: indeed, it is uniquely characterized as minimal, non-central, normal and killed by Frobenius, and any element $\sigma \in \Gamma_k$ will send it to another subgroup with the same properties. Then the factorisation of isogenies \ref{lem: factorisation isogenies} also descends from $k^\sep$ to $k$, and the claim follows.
\end{proof}



The following is the generalization of Proposition 
\ref{prop:descent} to the case of parabolic subgroups 
of quasi-Wenzel type.

\begin{corollary}
\label{cor:parabolicquasiWenzel}
A parabolic $k^\sep$-subgroup $P\subset G^\sep$ 
of the form (\ref{eqn:quasistandard}) is defined
over $k$ if and only if: 
$\{ \alpha_1, \ldots, \alpha_r \}$
is contained in 
$\Delta \setminus \Delta_{\min}$ and
is stable under the $*$-action, and 
the $(m_i,\varepsilon_i)$ are invariant under 
that action.
\end{corollary}

\begin{proof}
The descent of the reduced part follows from 
Corollary \ref{cor:smoothdescent}. It remains to show that, when $P$ is of quasi-Wenzel type but not of Wenzel type, then the $(m_i,\epsilon_i)$ are invariant under the action of $\Gamma_k$. Let $\Omega \subset \{\alpha_1,\ldots, \alpha_r\}$ be an orbit and consider the corresponding $k^\sep$-subgroup
\[
Q := \bigcap_{\alpha_i \in \Omega} (\ker \xi_i)P^{\alpha_i}, \qquad Q_{\text{red}}= (P^\Omega)^\sep.
\]
In order for $Q$ to be $\Gamma_k$-stable, its intersection with the root $k^\sep$-subgroup $U_{-\alpha_i} \subset G^\sep$ needs to be of constant height (for all the $\alpha_i \in \Omega$) equal to $m_\Omega$. Thus, this means that 
\[
\Omega = \Omega_\rho \cup \Omega_0,
\]
where the union is disjoint and $(m_i,\epsilon_i) = (m_\Omega,0)$ for all $\alpha_i \in \Omega_0$ while
\begin{align*}
    & (m_i,\epsilon_i) = (m_\Omega-1,1)\, \text{for } \alpha_i \in \Omega_\rho, \,  \text{ if the roots in $\Omega$ are short},\\
    & (m_i,\epsilon_i) = (m_\Omega,1) \, \text{for } \alpha_i \in \Omega_\rho, \,  \text{ if the roots in $\Omega$ are long.}
\end{align*}
If the simple factor of $G$ whose root system contains $\alpha_i$ satisfies the edge hypothesis, we denote as $\rho_i$ the corresponding very special isogeny (and extend it to $G$ as the identity on all other simple factors); see Definition \ref{def:edgehypothesis} and Lemma \ref{lem: very special kernel}.

Up to quotienting by $\rho_i^{-1}({}_{m_\Omega-1}G^\sep)$ in the first case, and by ${}_{m_\Omega}G^\sep$ in the second (which does not affect the property of whether $Q$ is defined over $k$, thanks to Lemma \ref{lem: kernels descend}), we can assume that the $\alpha_i$s in the orbit are all long and write
\[
Q = P^{\sep,\Omega_0} \cap N P^{\sep,\Omega_\rho},
\]
where $\Omega_0$ is non-empty. The original claim translates into saying that all the $\xi_i$ are all equal in the orbit: in other words, what is left to prove is that $\Omega_\rho$ must be empty. Assume that $\Omega_\rho$ is non-empty and let $\nu \in \Omega$ be the closest one (in terms of edges of the Dynkin diagrams) to the long roots. Let $\gamma$ be a positive short root whose support is a segment joining some short simple root to $\nu$; in particular, $\mathrm{Supp}(\gamma)\cap \Omega = \{\nu\}$. If $\nu \in \Omega_\rho$ (resp if $\nu \in \Omega_0$) consider $\alpha \in \Omega_0$ (resp $\alpha \in \Omega_\rho$), and consider $\sigma \in \Gamma_k$ which sends $\nu$ to $\alpha$. By \Cref{rem:subdiagrams}, the support of $\gamma$ is sent to another segment $J$ in the Dynkin diagram satisfying $J \cap \Omega = \{\alpha\}$ and containing one multiple edge. Thus, there is some $\delta$ positive long root with $\mathrm{Supp}(\delta) = J$. Then we have the following:
\begin{align*}
    & \text{if } \nu \in \Omega_\rho, \quad \mathrm{ht}(Q \cap U_{-\delta}) = 0 \text{ and } \mathrm{ht}(\sigma Q \cap U_{-\delta}) = 1,\\
     & \text{if } \nu \in \Omega_0, \quad \mathrm{ht}(Q \cap U_{-\delta}) = 1 \text{ and } \mathrm{ht}(\sigma Q \cap U_{-\delta}) = 0.
\end{align*}
This leads to $Q$ not being stable by $\Gamma_k$, hence if we want it to descend over $k$, $\Omega_\rho$ must be empty, as wanted.
\end{proof}

The above discussion on parabolics of quasi-Wenzel type leaves out only one case: let us assume that $p=2$ and that the group $G^\sep$ has at least one simple factor of type $G_2$. Since the $*$-action preserves the length of the simple roots, it cannot act non-trivially on a single such factor, and it can only permute the different factors of type $G_2$; so we can assume that \emph{all} are of type $G_2$, and there is at least two of them. Thus, let us denote as $G_{(i)}$ the simple factors of $G^\sep$, and let $\Delta_i = \{\beta_i,\gamma_i\}$ be its root basis over $k^\sep$, where $\beta_i$ is the short simple root and $\gamma_i$ the long one.

\begin{lemma}
\label{lem:parabolicsG2}
    Under the above assumptions, 
    a parabolic $k^\sep$-subgroup $P \subset G^\sep$ containing $B$ is defined over $k$ if and only if: $\{ \alpha_1, \ldots, \alpha_r \}$
is contained in 
$\Delta \setminus \Delta_{\min}$ and
is stable under the $*$-action, and moreover if
    \[
    P \cap G_{(i)} \simeq P \cap G_{(j)}
    \]
    for any $i,j$ such that $\beta_i$ and $\beta_j$ (or $\gamma_i$ and $\gamma_j$) are in the same orbit for the $*$-action.
\end{lemma}

\begin{proof}
Once again, the first assertion follows from Corollary
\ref{cor:smoothdescent}. Next, let $\Omega \subset \{\alpha_1,\ldots ,\alpha_r\}$ be an orbit of $\Gamma_k$, 
and let us assume that its roots are short. If for some $i$ such that $\beta_i \in \Omega$ we also have that $\gamma_i \in \{\alpha_1,\ldots, \alpha_r\}$, then (since the $*$-action preserves adjacency relations in the Dynkin diagram) the same must hold for all other $\beta_j$ in the same orbit $\Omega.$ The same reasoning holds if we start with an orbit made of long roots instead. So, we are in one among the following: 
\begin{align*}
    & (a)\,\, \Delta_i \cap \{\alpha_1,\ldots,\alpha_r\} = \{\gamma_i\}, & \text{for all $i$ such that } \Delta_i \cap \Omega \neq \emptyset;\\
    &(b)\,\, \Delta_i \cap \{\alpha_1,\ldots,\alpha_r\} = \{\beta_i\}, & \text{for all $i$ such that } \Delta_i \cap \Omega \neq \emptyset;\\
    &(c) \,\,\Delta_i \cap \{\alpha_1,\ldots,\alpha_r\} = \Delta_i, & \text{for all $i$ such that } \Delta_i \cap \Omega \neq \emptyset.
\end{align*}
In order for $P$ to be defined over $k$, we need $Q:= \langle P,P^\Omega\rangle$ (the smallest subgroup scheme containing both) to be $\Gamma_k$-stable. Assume we are in case $(a)$; then as illustrated in Subsection \ref{subsec:G2p2}, we can write
\[
Q = \bigcap_{i \in \Omega} {}_{m_i}G^\sep \cdot P^{\gamma_i}
\]
in a unique way, so for it to descend we need the $m_i$s to be constant and equal to some $m_\Omega$. This yields the desired equality, since all of the corresponding $\beta_i$s are in the Levi of $P_{\mathrm{red}}$, thus \[
P \cap G_{(i)} = Q \cap G_{(i)} =  ({}_{m_\Omega}G^\sep \cdot (P^{\gamma_i}) \cap G_{(i)})
\]
for any $i\in \Omega$. Next, consider case $(b)$: then each $Q^{\beta_i} := \langle P,P^{\beta_i} \rangle$ is not just described by one integer (the height of the Frobenius kernel), but it is uniquely determined by the heights of 
\[
P \cap U_{-\beta_i}, \quad P \cap U_{-2\beta_i-\gamma_i}, \quad P \cap U_{-3\beta_i -2\gamma_i}.
\]
Since $\Gamma_k$ acts transitively on all the $\beta_i$ in the orbit, then we must have that these three integers coincide, which gives exactly $P \cap G_{(i)} \simeq P \cap G_{(j)}$ for any $i,j$ such that $\beta_i,\beta_j \in \Omega$. Case $(c)$ is analogous, considering additionally the height of each $P \cap U_{-\gamma_i}$.
\end{proof}

The following is a uniform statement (independent of the Dynkin diagram and of the characteristic) which summarizes Proposition \ref{prop:descent}, Corollary \ref{cor:parabolicquasiWenzel} and Lemma \ref{lem:parabolicsG2}.

\begin{proposition}
\label{prop:parabolicsuniform}
Let $P \subset G^{\sep}$ be a parabolic $k^{\sep}$-subgroup containing $B$, and $\{\alpha_1,\ldots, \alpha_r\}$ the simple roots associated to its reduced part. Let $Q^{\alpha_i}$ be the smallest parabolic $k^\sep$-subgroup containing both $P$ and $P^{\alpha_i}$. Then $P$ is defined over $k$ if and only if: $\{\alpha_1,\ldots,\alpha_r\}$ is contained in $\Delta\backslash \Delta_{\mathrm{min}}$ and is stable under the $*$-action, and the equality
\[
\sigma \cdot Q^{\alpha_i} = Q^{\sigma^* \alpha_i}
\]
holds for any $i$ and for any $\sigma \in \Gamma.$ 
The last condition is equivalent to asking that the Wenzel function of $P$ is $*$-invariant.
\end{proposition}

\subsection{Rationality of parabolic varieties}
\label{subsec:rational}

We still consider a semisimple $k$-group $G$
and a parabolic $G$-variety $X$. 
Assuming that $k$ is \emph{perfect}, we
first obtain an analogue of Proposition
\ref{prop:rational}:

\begin{proposition}\label{prop:perfect}

Keep the above notation and assumptions.

\begin{enumerate}

\item[{\rm (i)}]
$X(k)$ is the unique orbit of $G(k)$.

\item[{\rm (ii)}]
The stabilizer of every $k$-rational point of $X$ 
contains a minimal smooth parabolic $k$-subgroup.
Moreover, $X$ has a unique $k$-rational 
$P_{\min}$-fixed point $x$. 

\item[{\rm (iii)}]
The orbit $V := R_u(P_{\min}^-) \cdot x$ is open 
in $X$ and isomorphic to an affine space. 

\item[{\rm (iv)}]
$X$ is covered by the translates $g \cdot V$,
where $g \in G(k)$.

\end{enumerate}

\end{proposition}

\begin{proof}
We prove (i) and (ii) simultaneously. 
We may assume that $X(k)$ is non-empty, and hence 
$X = G/P$ for some parabolic $k$-subgroup $P \subset G$.
Then the reduced subscheme $P_{\red} \subset P$ 
is a smooth $k$-subgroup. Moreover,
$P_{\red} \subset G$ is parabolic, since the
natural morphism $f : G/P_{\red} \to G/P$
is finite. By Proposition \ref{prop:hss},
$f$ induces a bijection on $k$-rational points.
Also, $G(k)$ acts transitively on 
$(G/P_{\red})(k)$ as the quotient morphism
$G \to G/P_{\red}$ is a Zariski locally
trivial $P_{\red}$-torsor. This yields (i).
Also, the smooth parabolic $P_{\red}$ contains 
a minimal smooth parabolic, i.e., a 
$G(k)$-conjugate of $P_{\min}$. This implies
readily the existence of a $P_{\min}$-fixed 
point $x \in X(k)$. The uniqueness follows 
from Proposition \ref{prop:rational} (i).

(iii) We adapt the argument of loc.~cit. (ii).
The multiplication map 
$R_u(P_{\min}^-) \times P_{\min} \to G$
is an open immersion, and hence the orbit
$V$ is open in $G$. It is isomorphic to 
an affine space as $R_u(P_{\min})$ is 
a split unipotent $k$-group.

(iv) For simplicity, let
$U := R_u(P_{\min})$ and 
$U^- := R_u(P_{\min}^-)$. 
By \cite[Proposition 6.25]{BT},
we have $G = U(k) \cdot U^- \cdot P_{\min}$, 
and hence 
$X = G(k) \cdot U^- \cdot x = G(k) \cdot V$.
\end{proof}

Note that $P_{\min}$ may be replaced
with its opposite $P_{\min}^-$ in the above
statement, since both are conjugate in $G(k)$.
Also, note that the statement (i) also follows 
from \cite[Lemma 1]{Srimathy} when $G$ is of inner 
type (in the sense of Definition \ref{def:*action}), 

We now consider an \emph{arbitrary} field $k$.
Then $X$ contains at most one $k$-rational 
$P_{\min}$-fixed point; moreover, if such a point 
$x$ exists, its orbit $G(k) \cdot x$ consists of 
the $k$-rational points of $X$ whose stabilizer 
contains a smooth $k$-parabolic (this may be checked 
by adapting the proof of Proposition 
\ref{prop:rational}). But we will see in Section
\ref{sec:example} that such a point $x$ does not exist 
in general, nor does there exist $x' \in X(k)$ 
such that $R_u(P_{\min}) \cdot x'$ is open in $X$; 
in addition, $X(k)$ may contain infinitely many 
orbits of $G(k)$. On the positive side, we now show:

\begin{theorem}\label{thm:krational}
If $p \geq 3$ (and $k$ is possibly imperfect),
then $X$ has an open covering by affine spaces.
\end{theorem}

\begin{proof}
Let $x \in X(k)$. It suffices to show that $x$ 
has an open neighborhood $V$ isomorphic to
an affine $k$-space: then one may check that
$X$ is covered by the translates $g \cdot V$ 
where $g \in G(k)$, by arguing as in the proof 
of Proposition \ref{prop:rational} (iii).

To construct the desired neighborhood $V$, 
we may assume that $G$ acts faithfully on $X$.
In particular, $G$ is adjoint.
Let $P := \Stab_G(x) \subset G$.
Consider the Lie algebras
$\fp := \Lie(P) \subset \Lie(G) =: \fg$
and the normalizer $Q := N_G(\fp)$.
Then $Q \supset P$ and hence $Q$ is 
a parabolic $k$-subgroup of $G$.

We claim that $Q$ is smooth and is not 
the whole $G$. To check this, we may assume $k$ 
algebraically closed, since the formations of 
the Lie algebra and the normalizer commute with 
field extensions. Then there exists a smallest
smooth parabolic $k$-subgroup $P^{\sm} \supset P$; 
moreover, we have $\Lie(P^{\sm}) = \fp$ and hence 
$P^{\sm} \neq G$
(see \cite[Remarks 1.6 and 3.14]{Maccan2}). 
Since $G$ is adjoint semisimple, the
$k$-group $N_G \big( \Lie(P^{\sm})\big)$ 
is smooth with neutral component $P^{\sm}$, see 
\cite[Expos\'e XXII, Proposition 5.3.1]{SGA3}.
It follows that $Q = P^{\sm}$; this proves 
the claim.

By that claim, we have a $G$-equivariant 
morphism $f : X = G/P \to G/Q$. Its fiber
at $y := f(x)$ is isomorphic to $Q/P$,
a projective homogeneous variety under $Q$, 
and hence under its semisimple quotient $Q/R(Q)$.
Moreover, the morphism
$R_u(Q^-) \times Q/P \to G/P$,
$(g,z) \mapsto g \cdot z$ is an open
immersion, where $Q^- \subset G$ denotes 
a smooth parabolic $k$-subgroup which is opposite
to $Q$. Since $R_u(Q^-)$ is an affine space,
we conclude by induction on the dimension
of $X$.
\end{proof}

\begin{remark}\label{rem:krational}
The above theorem holds more generally if all 
parabolic $k^{\sep}$-subgroups of $G^{\sep}$
that contain a Borel subgroup are of quasi-Wenzel type.
The proof extends unchanged to this setting.
\end{remark}


\section{Twisted parabolic varieties over an arbitrary field} 
\label{sec:tpvaf}

In this section, we still consider a semisimple $k$-group $G$,
where $k$ is an arbitrary field. Using the notions and results 
on root systems from Subsection \ref{subsec:mpr}, we first 
characterize those homogeneous spaces $G^{\sep}/P$ that are 
defined over $k$ in the following sense:

\begin{definition}\label{def:form}
A homogeneous space 
$G^{\sep}/H$ is said to be \emph{defined over} $k$ if  
there exists a $G$-variety $X$ over $k$ and a 
$G^{\sep}$-equivariant isomorphism 
$X^{\sep} \stackrel{\sim}{\longrightarrow} G^{\sep}/H$.
We then say that $X$ is a $k$-\emph{form} of
$G^{\sep}/H$.
\end{definition}

Since homogeneous spaces are quasi-projective, 
$G^{\sep}/H$ is defined over $k$ if and only if
it admits an action of $\Gamma_k$ which lifts
its action on $\Spec(k^{\sep})$ and is
compatible with the action of $G^{\sep}$.

\begin{lemma}
\label{lem:defined}
Let $P \subset G^{\sep}$ be a $k^{\sep}$-subgroup 
containing a Borel $k^{\sep}$-subgroup $B$. 

\begin{enumerate}

\item[{\rm (i)}]
The homogeneous space $G^{\sep}/P$ is defined 
over $k$ if and only if:
For any $\sigma \in \Gamma_k$, there exists
$g_{\sigma} \in G(k^{\sep})$ such that
$\sigma(P) = g_{\sigma} P g_{\sigma}^{-1}$.

\item[{\rm (ii)}]
Under these conditions, $G^{\sep}/P$ admits a 
unique $k$-form up to unique equivariant isomorphism.

\end{enumerate}

\end{lemma}

\begin{proof}
(i) Recall that the base point 
$o \in (G^{\sep}/P)(k^\sep)$
is the unique fixed point of $B$
(Proposition \ref{prop:rational}).
Assume that $G^{\sep}/P$ is defined over 
$k$ and let $\sigma \in \Gamma_k$. Then 
$\sigma(o)$ is fixed by $\sigma(B)$,
and the latter is a Borel $k^{\sep}$-subgroup
of $G^{\sep}$. Since any two such subgroups
are conjugate in $G(k^\sep)$, we have
$\sigma(B) = g_{\sigma} B g_{\sigma}^{-1}$
for some $g_{\sigma} \in G(k^{\sep})$. Then 
$g_{\sigma}^{-1} \sigma(o) \in (G^{\sep}/H)(k^{\sep})$
is fixed by $g_{\sigma}^{-1} \sigma(B) = B$, 
and hence equals $o$.
Equivalently, $\sigma(o) = g_{\sigma} \cdot o$.
Taking stabilizers in $G^{\sep}$, we obtain
$\sigma(P) = g_{\sigma} P g_{\sigma}^{-1}$.

For the converse, note that the coset
$g_{\sigma} P$ is unique, since $P$ is
its own normalizer (Remark \ref{rem:norm}). 
As a consequence, we have the cocycle condition
\[ g_{\sigma \tau} = 
\sigma(g_{\tau}) g_{\sigma} \]
for any $\sigma,\tau \in \Gamma_k$.
We define an action of $\Gamma_k$
on $G^{\sep}/P$ via 
\[ \sigma(g \cdot o) = 
\sigma(g) g_{\sigma} \cdot o 
\quad (g \in G^{\sep}). \]
If $h \in P$ then 
\[ \sigma(gh \cdot o) = 
\sigma(g) g_{\sigma}
(g_{\sigma}^{-1} \sigma(h) g_{\sigma}) \cdot o 
= \sigma(g) g_{\sigma} \cdot o, \]
as 
$g_{\sigma}^{-1} \sigma(h) g_{\sigma} \in H$.
Thus, $\sigma(g \cdot o)$ is well defined.
One may check that this yields an action 
of $\Gamma_k$ by using the cocycle condition.
Moreover, we have for any $g, g' \in G^{\sep}$:
\[ \sigma\big( g \cdot (g' \cdot o) \big) =
\sigma(gg' \cdot o) = 
\sigma(gg') g_{\sigma} \cdot o = 
\sigma(g) (g' \cdot o) \]
and hence 
$\sigma(g \cdot x) = \sigma(g) \cdot \sigma(x)$
for any $g \in G^{\sep}$, $x \in G^{\sep}/P$.
So the $\Gamma_k$-action is compatible
with the $G^{\sep}$-action.

(ii) It suffices to show that 
$G^{\sep}/P$ admits a unique action of
$\Gamma_k$ which lifts the action on 
$\Spec(k^{\sep})$ and is compatible with
the $G^{\sep}$-action. For this, we adapt
the above argument: the action of every 
$\sigma \in \Gamma_k$ satisfies
$\sigma(g \cdot o) = \sigma(g) \cdot \sigma(o)$
for any $g \in G_{k^{\sep}}$, and hence is
uniquely determined by $\sigma(o)$. Moreover,
$\sigma(o) \in (G^{\sep}/P)(k^{\sep})$ 
is fixed by 
$\sigma(P) = g_{\sigma} P g_{\sigma}^{-1}$,
and hence $\sigma(o) = g_{\sigma} \cdot o$
as $P$ is its own normalizer. Since the coset
$g_{\sigma} P$ is unique, so is $\sigma(o)$.
\end{proof}

\begin{proposition}
\label{prop:definedfunction}
Given a parabolic $k^{\sep}$-subgroup $P \subset G^{\sep}$ 
containing $B$ and with associated Wenzel function 
$\varphi \colon R^+ \rightarrow \bN \cup \{\infty\}$, 
the homogeneous space $G^{\sep}/P$ is defined over $k$ 
if and only if $\varphi$ is $*$-invariant.
Under these conditions, every $\sigma \in \Gamma_k$ acts 
on $G^{\sep}/P$ via 
$\sigma(g \cdot o) := \sigma(g) n_{\sigma}^{-1} \cdot o$,
where $n_{\sigma} \in N_{G^{\sep}}(T^{\sep})$ satisfies
$B = \Int(n_{\sigma}) \circ \sigma(B)$.
\end{proposition}

\begin{proof}
By Lemma \ref{lem:defined} and 
Proposition \ref{prop:authom}, 
$G^{\sep}/P$ is defined over $k$ if and only if 
for any $\sigma \in \Gamma_k$, there exists
$g_\sigma \in G(k^{\sep})$ 
such that $\sigma(P) = g_\sigma P g_\sigma^{-1}$.
Since $T^{\sep} \subset B \subset P$ and 
$g_{\sigma}^{-1} T^{\sep} g_{\sigma} \subset 
g_{\sigma}^{-1} \sigma(B) g_{\sigma} \subset P$,
there exists $h_{\sigma} \in P(k^{\sep})$
such that
$g_{\sigma}^{-1} T^{\sep} g_{\sigma} =
h_{\sigma} T^{\sep} h_{\sigma}^{-1}$ and
$g_{\sigma}^{-1} \sigma(B) g_{\sigma} =
h_{\sigma} B h_{\sigma}^{-1}$. 
Replacing $g_{\sigma}$ with 
$g_{\sigma} h_{\sigma}$,
we may thus assume that 
$T^{\sep} = g_{\sigma} T^{\sep} g_{\sigma}^{-1}$ 
and $\sigma(B) = g_{\sigma} B g_{\sigma}^{-1}$.
Then $g_{\sigma} \in n_{\sigma}^{-1} T^{\sep}$
and the coset 
$w_{\sigma} := g_{\sigma}^{-1} T^{\sep}$ satisfies
$\sigma^* = w_\sigma \sigma$. By the classification
of parabolic $k$-subgroups containing $B$ in terms
of Wenzel functions (Subsection \ref{subsec:rswf}),
it follows that $P_{\varphi}$ is defined over $k$
if and only if $\varphi$ is invariant under the
$*$-action.
\end{proof}

The above lemma combined with Proposition
\ref{prop:authom} easily yields an analogue of 
Remark \ref{rem:bij} over our ground field $k$ 
(not necessarily separably closed).

\begin{corollary}\label{cor:unique_k}
Let 
$\varphi: R^+ \rightarrow \mathbb{N} \cup \{\infty\}$ 
be a Wenzel function invariant under the $*$-action. 
Then the twisted parabolic $G$-variety supplied by 
Proposition \ref{prop:definedfunction} is unique up 
to a unique $G$-equivariant isomorphism. 
\end{corollary}

\begin{definition}
    Given a Wenzel function $\varphi: R^+ \rightarrow \mathbb{N} \cup \{\infty\}$, that is invariant under the $*$-action, we denote by $X_{\varphi}$ the unique twisted parabolic  $G$-variety given by Corollary \ref{cor:unique_k}.
\end{definition}

\begin{remark}\label{rem:leq}
Consider two Wenzel functions $\varphi,\psi$ such that 
$\varphi \leq \psi$. Then $P_{\varphi} \subset P_{\psi}$
and hence we have a $G^{\sep}$-equivariant morphism 
$G^{\sep}/P_{\varphi} \to G^{\sep}/P_{\psi}$
(such a morphism is unique in view of Proposition 
\ref{prop:rational} (i)).
If $\varphi,\psi$ are $*$-invariant, then this morphism
is also $\Gamma_k$-equivariant as a consequence
of Proposition \ref{prop:definedfunction}.
Thus, it descends to a $G$-equivariant morphism
$X_{\varphi} \to X_{\psi}$ (such a morphism is 
unique as well).
\end{remark}

\begin{lemma}
\label{lem:degree Xhat}
Assume that $k$ is perfect and let $\varphi$ be a Wenzel function 
such that $\varphi(\gamma) = 0$ for some
$\gamma \in R^+$. Then the degree of the universal 
purely inseparable cover 
$f : \widehat{X_\varphi} \rightarrow X_\varphi$ 
is equal to $p^M$, where 
    \[
    M := \sum_{\varphi(\gamma) \neq \infty} \varphi(\gamma).
    \]
\end{lemma}

\begin{proof}
Such a degree can be computed by passing to the algebraic closure. 
Then the universal purely inseparable cover becomes 
$f : \bar{G}/(P_\varphi)_{\red} \rightarrow 
\bar{G}/P_\varphi$.
We now use the notation of Subsection \ref{subsec:rswf},
and let $P := P_{\varphi}$, $Q := P_{\red}$, 
$U^- := R_u(Q^-)$ and $U^-_P := U^- \cap P$. 
Then the pullback of $f$ to the open orbit 
$U^- \cdot P/P \subset \bar{G}/P$ may be identified with
the quotient map $U^- \to U^-/U^-_P$.
So the degree of $f$ is the order of the infinitesimal
group scheme $U^-_P$. This gives the claimed formula
in view of Theorem \ref{thm:wenzel}.
\end{proof}

\begin{question} \label{q:birational}
    Let $k$ be  a perfect field of characteristic $p>0$. Are the twisted parabolic $G$-varieties $X_{\varphi}$ and $X_{\psi}$ birational if the underlying reduced parabolic subgroups over $\bar{k}$ corresponding to $\varphi$ and $\psi$ are  the same?
\end{question}

\begin{remark}\label{rmk:amitsur}  In other words, Question \ref{q:birational} asks  whether the $G$-varieties $X_{\varphi}$ and $X_{\psi}$ are birational as $k$-varieties if $\Supp(\varphi) = \Supp(\psi)$.  First, we note that the answer to this question is negative if $k$ is not assumed to be perfect.  To see this,  let $\Br(k)$ denote the Brauer group of $k$. For a central simple algebra $A$ over $k$ and a non-negative integer $m$,  define $A^{(m)}:= A \otimes_{Frob^m} k$
where $Frob^m$ denotes the $m$th absolute Frobenius.  
Recall from \cite[Proposition 3.2]{florence} that $A^{(m)}$ lies   in the Brauer class of $A^{\otimes p^m}$  i.e.,  $[A^{(m)}] = [A^{\otimes {p^m}}]$. Now suppose  $k$ is  an imperfect field such that  $p$-torsion of the Brauer group  $\Br(k)[p]\neq0$. Let  $0 \neq [A]\in \Br(k)[p]$ be a non-trivial class corresponding to the central simple algebra $A$ over $k$.   Let $G = \mathbf{PGL}_1(A) := \mathbf{GL}_1(A)/\bG_m$.    In our setting, $[A] \in \Br(k)[p]$ means that  $[A^{\otimes p}] = 1 \in \Br(k)$  and  hence $A^{(1)}$ is split.  The functoriality of the Frobenius implies  $\mathbf{PGL}_1(A)^{(1)} \simeq \mathbf{PGL}_1(A^{(1)})$, and therefore $G^{(1)}$ is split.  Let $Y = \mathbf{SB}(A)$ be the Severi-Brauer variety associated to $A$ with the canonical $G$-action so that $Y^{\sep} \simeq G^{\sep}/P$ for a maximal smooth parabolic $P \subset G^{\sep}$. Let $X = \mathbb{P}^r$ where $r = \deg(A) -1$ and the $G$-action is via $G \rightarrow G^{(1)}$. Then $X^{\sep} \simeq G^{\sep}/(_{1}G^{\sep} \cdot P)$ and therefore, the supports of the Wenzel functions corresponding to $X$ and $Y$ are the same. However, by  \cite[Theorem 5.1.3]{gille_szamuely}, $X$ is not birational to $Y$.  Still, it is not known  whether the answer to the above question would be affirmative if we restrict to twisted parabolic varieties with faithful $G$-action when $k$ is  imperfect.

  In fact, it is not clear whether the answer to the above question is affirmative over perfect fields even if  $X_{\varphi}$ and $X_{\psi}$ are isomorphic (as varieties forgetting the $G$-action) to classical twisted flag varieties for some (not necessarily the same) groups. For example, assume  $k$  is perfect.   Let  $G=\mathbf{PGL}_1(A)$ for a  non-split central simple algebra $A$ over $k$ and let $X = \mathbf{SB}(A)$ with the canonical $G$-action. Let  $X^{(m)}  \simeq \mathbf{SB}(A^{(m)})$ be the Frobenius twisted variety with $G$-action  via $G\rightarrow G^{(m)}$.  Note that $X^{\sep} \simeq G^{\sep}/P$  for some parabolic $P \subset G^{\sep}$ and that $(X^{(m)})^{\sep} \simeq G^{\sep}/({}_mG^{\sep}\cdot P)$  so that the supports of the respective Wenzel functions are the same. Moreover, $X$ and $X^{(m)}$ are  isomorphic to classical twisted flag varieties under the action of $G$ and $G^{(m)}$ respectively.       Since $k$ is perfect, the $p$-torsion of the Brauer group of $k$ is trivial. Therefore, $[A]$ and $[A^{(m)}]$ generate the same subgroup of the Brauer group. Hence, the birationality of $\mathbf{SB}(A)$ and $\mathbf{SB}(A^{(m)})$ would provide  supporting evidence for  Amitsur's conjecture (see \cite[Remark 5.4.3]{gille_szamuely}). 
\end{remark}

Next, we obtain a structure result for twisted parabolic
$G$-varieties. Given such a variety $X$, recall from
Proposition \ref{prop:picQ} and Remark \ref{rem:picQ} that 
the Picard group $\Pic(X)$ is free of finite rank, say $s$.
Moreover, $\Pic(X)_{\bQ}$ has a canonical basis 
$\{ D_1,\ldots,D_s \}$
which generates the cone of base-point-free classes, and is 
indexed by the set of Galois orbits in the canonical basis 
of $\Pic(X^{\sep})_{\bQ}$. This yields well-defined contractions
\[
\pi_j : X \longrightarrow X'_j \quad (j = 1, \ldots,s), \]
since the contraction associated with a base-point-free divisor
$D$ is unchanged when $D$ is replaced with $nD$ for some positive
integer $n$. We may now state:

\begin{theorem}\label{thm:structure}
Keep the above notation and assume that all parabolic 
$k^\sep$-subgroups of $G^\sep$ containing a Borel subgroup
are of Wenzel type.

\begin{enumerate}

\item[{\rm (i)}] 
Each $X'_j$ is a twisted parabolic $G$-variety
for a unique action of $G$ such that 
$\pi_j$ is equivariant.

\item[{\rm (ii)}]
The morphism 
\[ (\pi_1,\ldots,\pi_s) : X \longrightarrow 
X'_1 \times \cdots \times X'_s \]
is a $G$-equivariant immersion with image
being the smallest closed $G$-subvariety.

\item[{\rm (iii)}]
The $k$-group $G$ acts on each $X'_j$ via a unique
quotient $G^{(m_j)} = G/{_{m_j}} G$ such that
the homogeneous $G^{(m_j)}$-variety $X'_j$ has
a smooth stabilizer.

\item[{\rm (iv)}] 
We have $\Pic(X'_j) \simeq \bZ$ and the pull-back
of the ample generator of $\Pic(X'_j)$ is a
positive multiple of $D_j$.

\end{enumerate}

\end{theorem}

\begin{proof}
The assertion (i) follows from Blanchard's Lemma
(see \cite[Theorem 7.2.1]{Brion1}).

We now prove the remaining assertions simultaneously
by providing an alternative construction of the
contractions $\pi_j$. We will freely use the notation 
and conventions of Subsection \ref{subsec:mpr}. 

By Proposition \ref{prop:rational}, the Borel 
$k^{\sep}$-subgroup $B \subset G^{\sep}$ has a unique
$k^{\sep}$-rational fixed point in $X$. Therefore, 
there is a unique $G^{\sep}$-equivariant isomorphism
$X^{\sep} \simeq G^{\sep}/P$, where $P$ is a 
parabolic $k^{\sep}$-subgroup containing $B$.
We thus have 
\[ P = {_{m_1}} G^{\sep} \cdot P^{\alpha_1} \cap \cdots 
\cap {_{m_r}} G^{\sep} \cdot P^{\alpha_r} \]
as in (\ref{eqn:wenzel}). Since $G^{\sep}/P$
is defined over $k$, the set 
$\{ \alpha_1,\ldots,\alpha_r \}$
is stable under the $*$-action and the 
$m_i$ are invariant under this action
(Proposition \ref{prop:definedfunction}).
So we have 
\[ P = \bigcap_{\Omega} 
{_{m_{\Omega}}} G^{\sep} \cdot P^{\Omega}, 
\]
where $\Omega$ runs over the orbits of 
$\{ \alpha_1,\ldots,\alpha_r \}$ under the
$*$-action, 
$P^{\Omega} = 
\bigcap_{\alpha \in \Omega} P^{\alpha}$,
and $m_{\Omega} = m_i$ for any $i$ such that
$\alpha_i \in \Omega$. The natural morphism
\[  f_{\Omega} : Y := G^{\sep}/P \longrightarrow 
G^{\sep}/{_{m_{\Omega}}} G^{\sep} \cdot P^{\Omega} 
=: Y'_{\Omega} \]
is a contraction by Lemma \ref{lem:contraction},
and the product morphism
\[ \prod_{\Omega} f_{\Omega} : 
Y \longrightarrow \prod_{\Omega} Y'_{\Omega} \]
is an immersion with image being the unique closed
$G^{\sep}$-orbit, by Borel's fixed point theorem.
Also, we have canonical isomorphisms
\[
Y'_{\Omega} \simeq 
(G^{\sep})^{(m_{\Omega})}/(P^{\Omega})^{(m_{\Omega})}
\simeq 
(G^{\sep}/P^{\Omega})^{(m_{\Omega})},\]
and $Y'_{\Omega}$ is equipped with morphisms to
the homogeneous spaces 
$(G^{\sep}/P^{\alpha})^{(m_{\Omega})}$
($\alpha \in \Omega$), all of which have Picard
group $\bZ$. Moreover, the pull-backs of the 
ample generators form the canonical basis of
$\Pic(Y'_{\Omega})$.

Each homogeneous space $Y'_{\Omega}$
is defined over $k$ by Proposition 
\ref{prop:definedfunction} again. Moreover, 
the contraction $f_{\Omega}$
descends to a $G$-equivariant morphism
\[ \pi_{\Omega} : X \longrightarrow X'_{\Omega}, \]
where $X'_{\Omega}$ is a twisted parabolic
$G$-variety (Remark \ref{rem:leq}). 
Note that $\pi_{\Omega}$ is a contraction, 
since so is its base change to $k^{\sep}$. 
For the same reason, 
$G$ acts on $X'_{\Omega}$ via its quotient
$G^{(m_\Omega)}$, and the homogeneous
$G^{(m_\Omega)}$-variety $X'_{\Omega}$
has a smooth stabilizer. In addition,
the product morphism 
$\prod_{\Omega} \pi_{\Omega}$
is an immersion with image being the smallest
closed $G$-stable subvariety. Finally,
we have $\Pic(X'_{\Omega}) \simeq \bZ$
by Remark \ref{rem:picQ}, since
the canonical basis of $\Pic(Y'_{\Omega})$
is permuted transitively by $\Gamma_k$.

It remains to identify the contractions
$\pi_{\Omega}$ with those defined via
the canonical basis of $\Pic(X)_{\bQ}$.
For this, we recall results from Subsection 
\ref{subsec:wenzeltype} with a slightly
different notation: the canonical morphisms
\[ Y = G^{\sep}/P \longrightarrow 
G^{\sep}/{_{m_i}} G^{\sep} \cdot P^{\alpha_i}
= Z_i \quad (i = 1,\ldots,r) \]
are contractions, and 
the pull-backs of the ample generators 
of the $\Pic(Z_i)$ form the canonical basis
$(D_1,\ldots,D_r)$ of $\Pic(Y)$. 
In view of Remark \ref{rem:picQ}, it follows that 
the canonical basis of $\Pic(X)_{\bQ}$
consists of the 
$D_{\Omega} = \sum_{i, \alpha_i \in \Omega} D_i$ 
and the corresponding contractions are obtained
by Galois descent from the natural morphisms
$Y \to \prod_{i, \alpha_i \in \Omega} Z_i$,
or equivalently from the contractions
$f_{\Omega} : Y \to Y'_{\Omega}$.
\end{proof}

The special case where $X$ has Picard rank $1$
(i.e., $s = 1$) is worth mentioning.
This is equivalent to $\Pic(X)$ being 
isomorphic to $\bZ$.

\begin{corollary}\label{cor:picZ}
We still assume that all parabolic $k^\sep$-subgroups of 
$G^\sep$ containing a Borel subgroup are of Wenzel type.
Let $X$ be a twisted parabolic $G$-variety on which 
$G$ acts faithfully and having Picard rank $1$. 
Then $X$ has a smooth stabilizer. In particular, 
if $X$ has a $k$-rational point then $X \simeq G/P$,
where $P \subset G$ is a maximal smooth parabolic 
$k$-subgroup.
\end{corollary}

A further consequence of Theorem \ref{thm:structure}
is a structure result for all parabolic $k$-subgroups
that generalizes the Wenzel decomposition
(\ref{eqn:wenzel}):

\begin{corollary}\label{cor:fpb}
Keep the assumption of the above corollary,
and consider a parabolic $k$-subgroup $P \subset G$.
Then there exist integers $m_1,\ldots, m_s \geq 0$
and smooth parabolic $k$-subgroups 
$Q_i \subset G^{(m_i)}$ ($i = 1,\ldots,s$)
such that 
\[ P = \bigcap_{i = 1}^s (F^{m_i})^{-1}(Q_i)
\qquad \text{and} \qquad
\Pic(G^{(m_i)}/Q_i) \simeq \bZ
\quad (i = 1,\ldots,s). \]
Moreover, the pairs
$(m_1,Q_1), \ldots,(m_s,Q_s)$
are unique up to reordering.

If $k$ is perfect, then we have 
\begin{equation}\label{eqn:perfect} 
(F^{m_i})^{-1}(Q_i) = {_{m_i}}G \cdot P_i 
\qquad \text{and} \qquad
Q_i = P_i^{(m_i)}
\end{equation}
for unique maximal smooth parabolic $k$-subgroups 
$P_i \subset G$ ($i = 1,\ldots,s$).
\end{corollary}

\begin{proof}
The first assertion follows from Theorem 
\ref{thm:structure} applied to the parabolic
variety $G/P$.

For the second assertion, let 
$P_i = \big( (F^{m_i})^{-1}(Q_i) \big)_{\red}$.
Then $P_i \subset G$ is a smooth $k$-subgroup,
and $F^{m_i} : G \to G^{(m_i)}$
yields a homomorphism $P_i \to Q_i$ 
which induces a bijection on $\bar{k}$-points,
and hence is faithfully flat as $Q_i$ is smooth. 
This yields the equalities (\ref{eqn:perfect}).
Moreover, $P_i \subset G$ is a maximal smooth
parabolic, since so is $Q_i \subset G^{(m_i)}$
by the above corollary.

Conversely, if the first equality in 
(\ref{eqn:perfect}) holds then we have
$\big( (F^{m_i})^{-1}(Q_i) \big)_{\red} = P_i$
in view of Theorem \ref{thm:wenzel} (iv).
\end{proof}

\subsection{Small characteristics} We now turn to characteristic-free versions of the preceding statements, as done in the previous sections. Although they are slightly less transparent than the corresponding statements in characteristic $p\geq 5$, they hold in complete generality.

\begin{remark}
    The key point is the following: if we look at the proof of Theorem \ref{thm:structure}, we see that it relies on one fundamental argument. That is, that over the separable closure we recover any parabolic subgroup as an intersection of smooth parabolics, each fattened by the corresponding Frobenius kernel, by the unique decomposition (\ref{eqn:wenzel}). Moreover, the corresponding quotient maps are contractions, as proven in Lemma \ref{lem:contraction}, and this allows to use Blanchard's lemma. Now, over $k^\sep$, if we remove the assumption that $p \geq 5$, the same situation holds considering instead $Q^{\alpha_i}$, which we recall is defined as the subgroup generated by $P$ and $P^{\alpha_i}$ (the maximal smooth one associated with the simple root). More precisely, $P$ is uniquely determined as intersection of the $Q^{\alpha_i}$s, and the corresponding quotient maps $\pi_i$ defined in (\ref{maximal contraction}) are contractions by Lemma \ref{lem:contraction any char}.
\end{remark}

Let us recall that over $k^\sep$, a homogeneous space is uniquely determined by its Wenzel function. In particular, the only two homogeneous varieties having Picard rank one, with a faithful action and which have a non-smooth stabilizer are those corresponding to $\mathcal{Q}_1$ and $\mathcal {Q}_2$, following the notations of Subsection \ref{subsec:G2p2}. By the geometric description of \cite[Section 2.6]{Maccan1}, we know that the first one corresponds to the projective space $\bP^5$, while the second variety is a smooth hyperplane section, homogeneous under a group of type $G_2$, of the $\Sp_6$-homogeneous variety $\mathrm{LG(3,6)} = \Sp_6/P^\alpha$, which parametrizes isotropic vector subspaces of dimension 3 in a symplectic vector space of dimension 6 (here $\alpha$ is the long simple root).

So with this notation we just recalled, the above remark means that, by following exactly the proof of Theorem \ref{thm:structure}, replacing each ${}_{m_i}G \cdot P^{\alpha_i}$ with $Q^{\alpha_i}$, we get the following.

\begin{theorem}\label{thm:structurebis}
Let $G$ be semisimple over an arbitrary field $k$. Let $X$ be a twisted parabolic $G$-variety and denote as $\pi_j \colon X \rightarrow X_j^\prime$ the contractions associated to the canonical basis of $\Pic(X)_{\bQ}$. Then assertions (i), (ii) and (iv) of Theorem \ref{thm:structure} hold unchanged. Moreover, the $k$-group $G$ acts faithfully on each $X_j$
\begin{enumerate}
    \item [{\rm (iii)}]
    either via a unique quotient $G^{(m_j)}$ such that the homogeneous $G^{(m_j)}$-variety $X_j^\prime$ has a smooth stabilizer;
    \item [{\rm (iii')}] or (and this can only occur when the edge hypothesis is satisfied) via a unique quotient of the form $(G^\flat)^{(m_j)}$ such that the homogeneous $(G^\flat)^{(m_j)}$-variety $X_j^\prime$ has a smooth stabilizer;
    \item [{\rm (iii'')}] or (and this can only occur when $p=2$ and $G^\sep$ has at least one factor of type $G_2$) via a unique quotient  $G^{(m_j)}$ such that the homogeneous $G^{(m_j)}$-variety $X_j^\prime$ is either a form of $\bP^5$ or of $\mathrm{LG}(3,6).$
\end{enumerate}
\end{theorem}

For the sake of completeness, we also state the characteristic-free version of Corollary \ref{cor:fpb}. If $G^\sep$ has no simple factor of type $G_2$ when $p=2$, then any parabolic $k$-subgroup $P\subset G$ writes as 
\[
P = \bigcap_{i=1}^s \xi_i^{-1}  (Q_i),
\]
where $\xi_i \in \{ F^{m_i}_G, F^{m_i}_{G^\flat} \circ \rho_G \}$  are purely inseparable isogenies -- whose kernel is minimal with respect to inclusion, among those satisfying such an equality, see Remark \ref{rem:minimalkernels} -- and each $Q_i$ is a smooth parabolic $k$-subgroup either of $G^{(m_i)}$ or of $(G^{\flat})^{(m_i)}$. This decomposition is unique up to reordering, and the homogeneous spaces whose stabilizers are the $Q_i$s are all of Picard rank one. Finally, if $k$ is perfect then \[
 (\xi_i)^{-1}  (Q_i) = (\ker \xi_i)\cdot P_i
\]for unique maximal smooth parabolic $k$-subgroups $P_i \subset G$.

\section{Local rigidity of parabolic varieties}
\label{sec:rigidity}

By a result of Demazure (see \cite[Proposition 5]{Demazure}), 
flag varieties are locally rigid in the following sense: 
given a smooth proper morphism of schemes $\cX \to S$ 
and a point $s \in S$ such that the geometric fiber 
$\cX_{\bar{s}}$ is a flag variety, there exists an open 
neighborhood $U$ of $s$ over which every geometric fiber 
is a flag variety. A key ingredient of the proof is
the vanishing of $H^1(\cX_{\bar{s}},\cT_{\cX_{\bar{s}}})$; 
equivalently, every flag variety is infinitesimally rigid. 

This does not extend to parabolic varieties. Indeed, 
the vanishing of $H^1(X,\cT_X)$ for a smooth complete variety
$X$ implies that the $k$-group $\bfAut_X$ is smooth (see 
\cite[Expos\'e III, Lemme 4.2, Corollaire 5.6]{SGA1}).
But there are many examples of parabolic varieties $X$
such that $\bfAut_X$ is non-smooth, as follows from the main 
results of \cite{Maccan3} that determine the automorphism
group scheme of all parabolic varieties over an algebraically 
closed field. 

Further examples of parabolic varieties that are not
infinitesimally rigid include the hypersurfaces 
$X_m$ of bi-degree $(p^m,1)$ in $\bP^n \times \bP^n$ 
defined in (\ref{eqn:tiv}), where $m \geq 1$:
one may check that the family of all hypersurfaces of 
this bi-degree gives rise to nontrivial infinitesimal 
deformations.

Still, parabolic varieties $G/P$ satisfy local rigidity 
properties \emph{as $G$-varieties}, as shown by Sancho 
de Salas over an algebraically closed field again
(see \cite[Theorem 2.2, Theorem 3.3]{SS}). 
We will obtain analogous results for parabolic varieties 
over an arbitrary field $k$, starting with the case where 
the acting group is split.

\begin{theorem}\label{thm:rigid}
Let $G$ be a split semisimple $k$-group,
$B \subset G$ a Borel $k$-subgroup, 
$\cX$ a scheme equipped with an action of $G$, and 
$f : \cX \to S$ a proper morphism of schemes
such that $f$ is $G$-invariant and its fiber $\cX_s$
at every $s \in S$ is homogeneous under $G_{\kappa(s)}$. 
Assume that $S$ is reduced, connected and locally noetherian,
and $f$ has a section. Then there exists a unique parabolic 
$k$-subgroup $P \subset G$ containing $B$ and a unique 
$G$-equivariant isomorphism of $S$-schemes 
\[ F : G/P \times S \stackrel{\sim}{\longrightarrow} \cX \]
(where $G/P \times S$ is viewed as an $S$-scheme via 
the second projection).
\end{theorem}

\begin{proof}
We first show the uniqueness assertions. That on $P$ is 
straightforward, and that on $F$ follows from the triviality
of the group of $G$-equivariant automorphisms of 
$G/P \times S$ over $S$ (which in turn follows from 
Lemma \ref{lem:autnorm} and Proposition \ref{prop:authom}).

For the existence assertions, we will freely use 
the fact that $\cX$ is locally noetherian, since 
$f$ is of finite type and $S$ is locally noetherian.

As $f$ has a section, every fiber $\cX_s$ has a 
$\kappa(s)$-rational point. Since $\cX_s$ is a complete
$G_{\kappa(s)}$-variety, it has a $\kappa(s)$-rational
point that is fixed by the split solvable subgroup
$B_{\kappa(s)}$, or equivalently by $B$. Consider
the fixed point scheme $\cX^B$ (which is a closed subscheme 
of $\cX$ by \cite[II.1.3.6]{DG}), its reduced
subscheme $S' := (\cX^B)_{\red}$ and the restriction of 
$f$ to that subscheme, $u : S' \to S$. Then
the fiber of $u$ at every $s \in S$ satisfies 
$S'_s \stackrel{\sim}{\longrightarrow} \Spec \kappa(s)$.
Indeed, by the above, $S'_s$ has a unique point, and this
point is $\kappa(s)$-rational; moreover, given a maximal 
torus $T \subset B$, we have 
$S'_s \subset (\cX^B)_s \subset (\cX^T)_s = (\cX_s)^T$,
and the latter is \'etale over $\kappa(s)$
as a consequence of Lemma \ref{lem:fixed}.

We now claim that \emph{$u$ is an isomorphism}.
To show this, note that $u$ is bijective,
and proper as $S'$ is closed in $X$. By
\cite[Theorem 8.11.1]{EGAIV3}, it follows that
$u$ is finite. For any $s \in S$ with
preimage $s' \in S'$, the homomorphism
of local rings 
$u^{\#}_s : \cO_{S,s} \to \cO_{S',s'}$
satisfies 
$u^{\#}_s(\fm_s) \cO_{S',s'} = \fm_{s'}$
as $S'_s = \Spec \kappa(s)$.
By Nakayama's lemma, it follows that
$u^{\#}_s$ is surjective, and hence
$u$ is a closed immersion. As $u$
is surjective and $S$ is reduced, this implies
our claim. 

We will identify $S'$ with $S$ via $u$,
and hence equip every fiber $\cX_s$ with 
its $B$-fixed base point $s$.
Next, let $\eta \in S$ be a maximal point,
and $K := \kappa(s)$ its residue field.
Then we have a $G_K$-equivariant isomorphism
$G_K/H \stackrel{\sim}{\longrightarrow} \cX_{\eta}$,
where $H \subset G_K$ is a $K$-subgroup containing
$B_K$. 

We now use the structure of $H$ from Theorem 
\ref{thm:wenzel} and the subsequent Remark
\ref{rem:wenzel}.
Since $G$ is split, $H_{\red} = Q_K$
for a unique smooth $k$-subgroup $Q \subset G$
containing $B$. Thus, $Q$ fixes $\eta$ and hence
its closure $C$, an irreducible component of $S$.
Also, for any root $\gamma$ of $R_u(Q)$, 
the isomorphism 
$\bG_a \stackrel{\sim}{\longrightarrow} U_{-\gamma}$
identifies $(U_{-\gamma})_K \cap H$ with
$(\bfalpha_{p^{\varphi(\gamma)}})_K$,
where $\varphi$ denotes the Wenzel function
of $H$. Thus, $\bfalpha_{p^{\varphi(\gamma)}}$
fixes $\eta$ as well. Let
\[ P_{\varphi} := 
\big( \prod_{\gamma} \bfalpha_{p^{\varphi(\gamma)}} \big)
\times Q \subset R_u(Q^-) \times Q \subset G. \]
Then $(P_{\varphi})_K = H$. Thus, $P_{\varphi}$
is a $k$-form of $H$, which fixes $\eta$
and hence $C$. So the $G$-action on $\cX$ yields 
a $G$-equivariant morphism of $C$-schemes
\[ F_C : G/P_{\varphi} \times C
\longrightarrow \cX_C, \]
where $\cX_C := \cX \times_S C$.

We now claim that \emph{$F_C$ is an isomorphism}.
Indeed, the pull-back $f_C : \cX_C \to C$ satisfies
$f_C \circ F_C = \pr_C$,
the projection 
$G/P_{\varphi} \times C \to C$. Moreover,
$\pr_C$ is faithfully flat and the pull-back
$(F_C)_s : (G/P_{\varphi})_{\kappa(s)} \to \cX_s$
is faithfully flat as well for any $s \in S$ 
(since it is an equivariant morphism of homogeneous 
varieties). By the fiberwise flatness criterion (see 
\cite[Corollary 11.3.11]{EGAIV3}), it follows
that $F_C$ and $f_C$ are faithfully flat.
In particular, they are equidimensional
(see \cite[Corollaire 6.1.4]{EGAIV2}).

Since $C$ is integral and $G/P_{\varphi}$
is geometrically integral, 
$G/P_{\varphi} \times C$ is integral.
As a consequence, $\cX_C$ is irreducible, and 
reduced by \cite[Proposition 2.1.13]{EGAIV2}.
Equivalently, $\cX_C$ is integral. Moreover,
$F_C$ is birational by construction.
It is also proper in view of the 
properness of $\pr_C$ together with
\cite[Corollaire 5.4.3]{EGAII}. 
Since $F_C$ is equidimensional, it must be
quasi-finite, and hence finite. But every
finite flat birational morphism of reduced
schemes is an isomorphism (this follows from
Nakayama's lemma as at the beginning of
the proof). This yields our claim, and
proves the theorem in the case where $S$ is
integral.

In the general case (where $S$ is reduced
and connected), consider the family of
irreducible components $(C_i)_{i \in I}$
of $S$. By the above step, each pull-back
$F_{C_i}$ is isomorphic to the projection
$G/P_i \times C_i \to C_i$ for a unique
$k$-subgroup $P_i$ containing $B$. Choose
such a $P := P_{i_0}$ and let 
\[
S_P := \bigcup_{i \in I, P_i = P} C_i,
\quad 
S^P := \bigcup_{i \in I, P_i \neq P} C_i.
\]
Clearly, we have $S = S_P \cup S^P$.
Moreover, $S_P \cap S^P$ is empty, since
$P$ is constant on each irreducible 
component by the above step again. Also,
$S_P$ and $S^P$ are closed in $S$: more 
generally, since $S$ is locally noetherian, 
every union of irreducible components 
$\bigcup_{j \in J} C_j$ is closed
(if $x \in S \setminus \bigcup_{j \in J} C_j$
and $V$ is an open affine noetherian 
neighborhood of $x$ in $S$, then $V$ meets
only finitely many $C_i$'s, and hence
$V \setminus \bigcup_{j \in J} C_j$
is an open neighborhood of $x$ in $S$,
contained in 
$S \setminus \bigcup_{j \in J} C_j$).
As $S$ is connected, it follows that
$S = S_P$. This yields a $G$-equivariant
morphism of $S$-schemes 
$F : G/P \times S \to \cX$,
which induces an isomorphism over
each irreducible component of $S$.
By the fiberwise isomorphism criterion
of \cite[Corollaire 17.9.5]{EGAIV4},
we conclude that $F$ is an isomorphism.
\end{proof}

\begin{corollary}\label{cor:rigid}
Let $S$ be a reduced, locally noetherian
$k$-scheme. Let $\cG$ be a semisimple $S$-group 
scheme, and $\cX$ be a smooth proper $S$-scheme
equipped with an action of $\cG$ such that
each fiber $\cX_s$ is homogeneous under 
$\cG_{\kappa(s)}$. Then every point of $S$ 
admits an open neighborhood $U$ together
with an \'etale surjective morphism $V \to U$
and compatible isomorphisms of $V$-schemes
\[ \cG_V \simeq G \times V, 
\quad \cX_V \simeq G/P \times V, \]
where $G$ is a split semisimple $k$-group
and $P \subset G$ is a parabolic subgroup
containing a Borel subgroup of $G$.
\end{corollary}

\begin{proof}
The assertion on $\cG$ is a consequence of 
\cite[Expos\'e XXII, Corollaire 2.3;
Expos\'e XXV, Corollaire 1.2]{SGA3}.
We may thus assume that $\cG = G \times S$
for a split semisimple $k$-group $G$.
Since $f$ is smooth, we may further assume
that it has a section (see 
\cite[Corollaire 17.15.3]{EGAIV4}). 
Then the assertion on $\cX$ follows from 
Theorem \ref{thm:rigid}.
\end{proof}

The following result is a version of 
\cite[Theorem 3.3]{SS}, with weaker assumptions
(in particular, we do not assume that
orbits have constant ``nilpotency degree'').

\begin{theorem}\label{thm:closed}
Let $G$ be a semisimple group over an
algebraically closed field $k$, and
$B \subset G$ a Borel $k$-subgroup.
Let $X$ be a normal $G$-variety such that 
the $G$-orbits of all $k$-rational points
are complete and have the same dimension. 
Then there is a unique $G$-equivariant 
isomorphism
\[ F : G/P \times S 
\stackrel{\sim}{\longrightarrow} X, \]
where $S \subset X$ is a closed subvariety,
$P \subset G$ is a parabolic $k$-subgroup containing
$B$, and $G$ acts on $G/P \times S$
via its left action on $G/P$.
\end{theorem}

\begin{proof}
The uniqueness assertion is checked
as in the proof of Theorem \ref{thm:rigid}.
We will also adapt this proof to show
the existence assertion, but we have a priori 
no $G$-invariant morphism $f : X \to S$ with 
homogeneous fibers. Still, such a morphism
exists over a dense open $G$-stable subset
of $X$, by a classical result of Rosenlicht
(see \cite[Theorem 2, p.~407]{Ro56}, and 
\cite[Theorem 1.1]{BGR} for a modern version).
More specifically, there exists such a subset $X_0$ 
having a geometric quotient $\pi : X_0 \to Y$, 
i.e., $\pi$ is a surjective morphism of varieties;
for any field extension $k \subset K$
and any $x \in X_0(K)$, the set-theoretic 
fiber of $\pi$ at $x$ is the $G$-orbit
$G \cdot x$; and $\pi$ induces an isomorphism
$k(Y) \stackrel{\sim}{\longrightarrow} k(X_0)^G$. 
By replacing $Y$ with a dense open subset and 
$X_0$ with its preimage, we may assume that $\pi$ 
is faithfully flat 
(see \cite[Th\'eor\`eme 6.9.1]{EGAIV2}).
Also, the field extension
$k(X_0)^G = k(X)^G \subset k(X) = k(X_0)$
is separable by \cite[Lemma, p.~4]{Ro55},
and hence the generic fiber of $\pi$ is
geometrically reduced. Consider the subset
of $X_0$ consisting of those $x$ such that the 
(schematic) fiber $(X_0)_x$ is geometrically 
reduced at $x$. This subset is $G$-stable, 
open in $X$ by 
\cite[Corollaire 12.1.7]{EGAIV3},
and nonempty as it contains the generic fiber.
By shrinking further $X_0$ and $Y$, we may thus 
assume that all fibers of $\pi$ are geometrically 
reduced. Then $\pi$ is smooth, since its 
set-theoretic fibers are homogeneous. Shrinking $Y$ 
once more, we may assume in addition that it is 
smooth; then $X_0$ is smooth as well.

We now follow closely the argument of Theorem 
\ref{thm:rigid}: consider the reduced subscheme 
$Z := (X_0^B)_{\red}$. Then $\pi$ restricts
to a morphism $u : Z \to Y$
which induces isomorphisms 
$Z_y \stackrel{\sim}{\longrightarrow} y$
for all $y \in Y(k)$. Since $k$ is algebraically 
closed, it follows that $u$ is quasi-finite. 
By Zariski's main theorem (see 
\cite[Theorem 8.12.6]{EGAIV3},
$u$ factors as
\[ X_0 \stackrel{i}{\longrightarrow}
X' \stackrel{v}{\longrightarrow} Y, \]
where $i$ is an open immersion
and $v$ is finite. Replacing
$Y$ with its dense open subset
$Y \setminus v \big( X' \setminus i(X_0) \big)$,
we may thus assume that $u$ is finite. 
Now it follows as in the proof of Theorem 
\ref{thm:rigid} that $u$ is a section of $\pi$. 
By that theorem, this yields an isomorphism 
\[ G/P \times Z 
\stackrel{\sim}{\longrightarrow} Y \]
for a unique parabolic $k$-subgroup 
$P \subset G$, which contains $B$ and 
fixes $Z$ pointwise. Thus, $P$
fixes pointwise the closure $S$ of
$Z$ in $X$, and the resulting morphism
$F : G/P \times S \to X$ is birational.
Since all $G$-orbits of $k$-rational
points of $X$ have the same dimension, 
$F$ is quasi-finite. 
It is also proper, since it factors as 
\[ G/P \times S \longrightarrow 
G/P \times X \longrightarrow 
G/P \times X \longrightarrow X, 
\quad (g,s) \longmapsto (g,s) \longmapsto 
(g, g \cdot s) \longmapsto s, \]
where the left map is a closed immersion, the
middle one is an isomorphism, and the right
one is proper. So $F$ is an isomorphism 
by a corollary of Zariski's main theorem
(see \cite[Corollaire 8.12.10]{EGAIV3}).
\end{proof}  

\begin{remark}\label{rem:pinching}
The normality assumption in the above theorem 
cannot be suppressed, as shown by examples
obtained from pinching along a $G$-orbit. 
Specifically, let $G$ be a semisimple $k$-group, 
$P\subsetneq G$ a parabolic $k$-subgroup,
and $Y$ a quasi-projective variety having a
$k$-rational point $y$.
There exists a parabolic $k$-subgroup 
$Q \subsetneq G$ such that $P \subsetneq Q$
and $Q/P$ is finite (indeed, we may take 
$Q = {}_m G \cdot P$ for $m \gg 0$).
The resulting finite morphism $u: G/P \to G/Q$
fits in a cocartesian diagram of varieties 
\[ \xymatrix{ 
G/P \ar[r]^-u \ar[d]_{(\id,y)} & G/Q \ar[d]^v \\
G/P \times Y \ar[r]^-f & X, \\
} \]
where $f$ is finite and $v$ is a closed immersion
(see \cite[Th\'eor\`eme 5.4]{Ferrand}).
Moreover, the natural action of $G$ on 
$G/P \times Y$ induces an action on $X$
such that $f$ and $v$ are equivariant
(as follows from \cite[Lemma 3.5]{Laurent}).
Thus, the $G$-orbits of $k$-rational
points of $X$ are all complete,
and isomorphic to $G/P$ except for one 
isomorphic to $G/Q$. In particular,
there is no equivariant isomorphism
$G/P \times S 
\stackrel{\sim}{\longrightarrow} X$
as in Theorem \ref{thm:closed}. 
Note finally that $X$ is non-normal,
since the morphism $f$ is finite,
birational and not an isomorphism.
\end{remark}

\section{Examples over imperfect fields}
\label{sec:example}

\subsection{}
\label{subsec:ex1}
Throughout this subsection, we assume that the ground field 
$k$ is \emph{imperfect}, and consider a \emph{split} 
adjoint simple $k$-group $G$ of rank at least two. 
We will construct examples of parabolic $G$-varieties 
$X$ on which $G$ acts faithfully, and such that $G(k)$ 
has infinitely many orbits in $X(k)$.

Since $G$ is split, we may choose a pair $(T,B)$, 
where $B \subset G$ is a Borel $k$-subgroup
and $T \subset B$ is a split maximal $k$-torus. 
With the notation of \S \ref{subsec:rswf}, 
let $\alpha_1, \alpha_2$ be distinct simple roots,
$P_1 := P_{\alpha_1}, P_2 := P_{\alpha_2}$
the corresponding smooth parabolic $k$-subgroups 
containing $B$, and $m \geq 2$ an integer. Then 
\[ P = ({}_m G \cdot P_1) \cap P_2 \]
is a non-smooth parabolic $k$-subgroup of $G$,
with largest smooth $k$-subgroup $P_1 \cap P_2 = B$. 
Thus, $X := G/P$ is a parabolic variety 
with non-smooth stabilizer. 

\begin{proposition}
\label{prop:imperfect no torus}

\begin{enumerate}

\item[{\rm (i)}] The $G$-action on $X$
is faithful.

\item[{\rm (ii)}]
There is $x \in X(k)$ such that the parabolic 
$k$-subgroup $Q := \Stab_G(x)$ contains 
no maximal $k$-torus of $G$.

\item[{\rm (iii)}]
The set $X(k)$ contains infinitely many orbits of $G(k)$.

\item[{\rm (iv)}]
The homogeneous variety $X$ admits no universal 
purely inseparable cover.

\end{enumerate}

\end{proposition}

\begin{proof}
(i) The natural morphism
\[ f : X = G/P \longrightarrow G/P_2 \]
is $G$-equivariant, and the $G$-action on $G/P_2$ 
is faithful as $G$ is adjoint simple and
$P_2$ is smooth (see e.g.
\cite[Proposition 1]{Demazure}). 
Thus, $G$ acts faithfully on $X$ as well.

(ii) The above morphism $f$ is a Zariski locally 
trivial fibration with fiber $P_2/P$ at the base 
point $o \in (G/P_2)(k)$. To determine this
fiber, we use the Levi decomposition 
$P_2 = R_u(P_2) L_2$, where $L_2 \supset T$. 
Since $R_u(P_2) \subset P$, we have
\[ P_2/P = L_2/({}_m G \cdot P_1) \cap L_2
= L_2/({}_m L_2) \cdot (P_1 \cap L_2), \]
where the second equality follows from 
\cite[Theorem 10]{Wenzel} and descent. 
Furthermore, $B_{L_2} = P_1 \cap L_2$ 
is a Borel $k$-subgroup of $L_2$ and 
\[ L_2/({}_m L_2 \cdot (P_1 \cap L_2)) =
(L_2/({}_m L_2)/(B_{L_2}/{}_m B_{L_2})
= L_2^{(m)}/B_{L_2}^{(m)}, 
\]
where the second equality follows from the
smoothness of $L_2$ and of $B_{L_2}$. Together, 
this yields identifications
\[ 
P_2/P  = (L_2/B_{L_2})^{(m)} = (\bP^1)^{(m)}.
\]

We now consider $k$-rational points.
The map $f(k) : X(k) \to (G/P_2)(k)$ is 
$G(k)$-equivariant with fiber at $o$ being 
\[ (P_2/P)(k)  = (\bP^1)^{(m)}(k)
=  \bP^1(k^{1/p^m}), \]
where the second equality follows from
the definition of the Frobenius twist as 
a base change. Also, we have 
$(G/P_2)(k) = G(k)/P_2(k)$
by \cite[Th\'eor\`eme 4.13]{BT}
and the smoothness of $P_2$. All of this 
yields an identification of sets of orbits
\begin{equation}\label{eqn:orbits} 
X(k)/G(k) = \bP^1(k^{1/p^m})/P_2(k), 
\end{equation}
where $P_2(k)$ acts on $\bP^1(k^{1/p^m})$
via th composition of homomorphisms
$P_2(k) \to L_2(k) \to \PGL_2(k)$
and the natural action of $\PGL_2(k)$
on $\bP^1$.

Since $k$ is imperfect, we may choose 
$t \in k \setminus k^p$. Then the polynomial
$T^{p^m} - t$ is irreducible in $k[T]$,
i.e., the powers 
\[ 1,t^{1/p^m}, t^{2/p^m}, \ldots, 
t^{(p^m - 1)/p^m} \in k^{1/p^m} \]
are linearly independent over $k$.
Let $x \in \bP^1(k^{1/p^m})$ be the point 
with homogeneous coordinates $(t^{1/p^m},1)$,
and $Q := \Stab_G(x)$. Then $Q$ contains 
the kernel $N$ of the homomorphism 
$P_2 \to \bfAut_{\bP^1} = \PGL_2$,
and $Q \subset P_2$ as $f(x) = o$.
Thus, the exact sequence of $k$-groups
$1 \to N \to P_2 \to \PGL_2$ yields 
an exact sequence
\[ 1 \longrightarrow N \longrightarrow Q 
\longrightarrow \Stab_{\PGL_2}(x). \]
To prove our assertion, it suffices to show 
that $Q(k^{\sep}) = N(k^{\sep})$, and in turn 
that the stabilizer of $x$ in $\PGL_2(k^{\sep})$ 
is trivial. By an elementary computation, 
this stabilizer consists of the images
of the matrices 
$\left(\begin{smallmatrix}
a & b \\ c & d \\    
\end{smallmatrix}\right) \in \GL_2(k^{\sep})$
such that 
$b + (a - d) t^{1/p^m} - c t^{2/p^m} = 0$.
Since $m \geq 2$, then $1,t^{1/p^m},t^{2/p^m}$
are linearly independent over $k$,
and hence over $k^{\sep}$ by MacLane's
criterion. As a consequence, $x$ has 
a trivial stabilizer in $\PGL_2(k^{\sep})$,
as desired.

(iii) In view of (\ref{eqn:orbits}), 
it suffices to show that $\bP^1(k^{1/p^m})$
contains infinitely many orbits of $\PGL_2(k)$.
Let $x,y,z,w \in k$ with $y \neq 0$ and
$w \neq 0$. We claim that the points of 
$\bP^1(k^{1/p^m})$ with homogeneous coordinates 
$(x + y t^{1/p^m},1)$ and $(z t^{1/p^m} + w t^{2/p^m},1)$
are not in the same orbit of $\PGL_2(k)$. Otherwise,
\[ 1, x + y t^{1/p^m}, z t^{1/p^m} + w t^{2/p^m},
(x + y t^{1/p^m}) (z t^{1/p^m} + w t^{2/p^m}) 
\in k^{1/p^m} \]
are linearly dependent over $k$. This is
equivalent to 
$1, t^{1/p^m}, t^{2/p^m}, t^{3/p^m}$
being linearly dependent over $k$, a contradiction
since $m \geq 2$. This proves the claim, 
and in turn the assertion as $k$ is infinite.


(iv) Assume that there is a universal purely 
inseparable cover $f : \widehat{X} \to X$, i.e., 
a $G$-equivariant morphism from a hss $G$-variety
which is universal for this property.
Then the canonical morphism
$u : G/B \to G/P = X$ factors uniquely
through an equivariant morphism
$\widehat{u} : G/B \to \widehat{X}$.
As a consequence, we have
$\widehat{X} \simeq G/\widehat{P}$,
where $B \subset \widehat{P} \subset P$.
Since $\widehat{X}$ is hss, it follows
that $\widehat{P} = B$. Next, consider 
the largest smooth $k$-subgroup 
$Q^+ \subset Q$, which exists 
by \cite[Lemma C.4.1]{CGP}, and
the canonical morphism 
$v : G/Q^+ \to G/Q = X$. Arguing similarly, 
we obtain an isomorphism 
$\widehat{X} \simeq G/Q^+$. Thus,
the smooth $k$-groups $B$ and $Q^+$ 
become isomorphic after base change 
to $\bar{k}$; in particular,
they have the same rank (see 
\cite[Expos\'e XIV, Th\'eor\`eme 1.1]{SGA3}). 
But this contradicts the fact that $Q^+$ 
contains no maximal $k$-torus of $G$.
\end{proof}

\begin{remarks}
(i) The above proposition also holds if $m = 1$
and $p \geq 5$.

(ii) The $G$-variety $X$ is equivariantly isomorphic to 
both $G/P$ and $G/Q$, where $Q := \Stab_G(x)$. But the 
parabolic $k$-subgroups $P,Q$ are not conjugate in 
$G(k)$; their base changes $P_{\bar{k}}, Q_{\bar{k}}$
are conjugate in $G(\bar{k})$. 

(iii) If the center $Z(L_2)$ is smooth and connected, 
then the (non-smooth) parabolic $k$-group $Q$
is integral; in particular, reduced.
Indeed, the exact sequences
\[ 1 \to R_u(P_2) \to P_2 \to 
L_2 \to 1, \quad 
1 \to Z(L_2) \to L_2 \to \PGL_2 \to 1 \]
yield exact sequences
\[ 1 \to R_u(P_2) \to Q \to 
\Stab_{L_2}(x) \to 1, \quad 
1 \to Z(L_2) \to \Stab_{L_2}(x) \to 
\Stab_{\PGL_2}(x) \to 1, \]
and hence it suffices to check that
$\Stab_{\PGL_2}(x)$ is integral.
But the pull-back of this $k$-group
to $\GL_2$ is the zero subscheme of
the polynomial
\[ \big( b + (a - d) t^{1/p^m} - c t^{2/p^m} \big)^{p^m} 
= b^{p^m} + (a^{p^m} - d^{p^m}) t - c^{p^m} t^2 \]
in the matrix coefficients $a,b,c,d$,
and this polynomial is irreducible over $k$. 
The smoothness and connectedness assumptions 
on $Z(L_2)$ hold for example if $G = \PGL_n$;
this yields many examples of reduced, non-smooth
parabolic subgroups.

(iv) By \cite[Corollaire 6.4]{BS},
if $k$ is a \emph{perfect} field having only 
finitely many extensions of any prescribed 
degree, then for any smooth $k$-group $G$ and
any homogeneous $G$-variety $X$, the set
$X(k)$ consists of finitely many orbits
of $G(k)$. This does not extend to
imperfect fields in view of the above
proposition: just take for $k$ a separably
closed field such that $[k:k^p] = 1$.
\end{remarks}

\subsection{}
\label{subsec:ex2}
In this subsection and the next one, we consider 
an imperfect field $k$ and an \emph{isotropic} 
simple adjoint $k$-group $G$. We will construct 
examples of pairs $(G,X)$ where $X$
is a parabolic $G$-variety on which $G$ acts
faithfully, and the stabilizer of every $x \in X(k)$
contains no smooth parabolic $k$-subgroup of $G$.

Our construction goes as follows: choose
a maximal split $k$-torus $S$ with centralizer 
$L := C_G(S)$, and a minimal parabolic $k$-subgroup 
$P_{\min} \subsetneq G$ containing $S$. Recall that 
$L$ is a Levi subgroup of $P_{\min}$. 
Moreover, the derived subgroup of $L$ is anisotropic, 
and hence $L$ contains 
no nontrivial smooth parabolic $k$-subgroup. 
We will construct pairs of $k$-groups $G$
as above and faithfully flat homomorphisms
$\varphi : L \to G'$, where $G'$ is an \emph{isotropic} 
simple adjoint $k$-group. Thus, $G'$ contains 
a nontrivial smooth parabolic $k$-subgroup $P'$. 
This will provide the desired examples in view 
of the following observation.

\begin{proposition}\label{prop:isotropic}
With the above notation, let $P \subset G$
be the inverse image of $P' \subset G'$
under the composition
\[ P_{\min} \stackrel{\psi}{\longrightarrow}
P_{\min}/R_u(P_{\min}) \simeq L 
\stackrel{\varphi}{\longrightarrow} G',
\]
where $\psi$ is the quotient homomorphism. Then, 

\begin{enumerate}

\item[{\rm (i)}]
$P \subset G$ is a parabolic $k$-subgroup, 
strictly contained in $P_{\min}$, and
$G$ acts faithfully on $X := G/P$.

\item[{\rm (ii)}]
For any $x \in X(k)$, the orbit 
$R_u(P_{\min}) \cdot x$ is not open in $X$.
Moreover, the stabilizer of $x$ in $G^{\sep}$
contains no Borel $k^{\sep}$-subgroup of
$G^{\sep}$.

\end{enumerate}

\end{proposition}

\begin{proof}
(i) By construction, $P$ is strictly contained in
$P_{\min}$. This yields a $G$-equivariant
morphism
\[ f : X = G/P \longrightarrow G/P_{\min}.
\]
Since $P_{\min}$ is a smooth parabolic,
$f$ is a Zariski locally trivial fibration.
Its fiber at the base point is 
$P_{\min}/P \simeq L/\varphi^{-1}(P') 
\simeq G'/P'$. 
As a consequence, $f$ is proper, so that 
$P \subset G$ is a parabolic $k$-subgroup. 
Moreover, $G$ acts faithfully on $G/P_{\min}$ 
(see \cite[Proposition 1]{Demazure}) and hence 
on $X$. 

(ii) The first assertion follows from the inequality
\[ \dim(X) = \dim(G/P_{\min}) + \dim(P_{\min}/P)
> \dim(G/P_{\min}) = \dim \big( R_u(P_{\min}) \big). \]
For the second assertion, let $Q := \Stab_G(x)$
and assume that $Q^{\sep}$ contains a Borel 
$k^{\sep}$-subgroup of $G^{\sep}$. Then the
largest smooth $k^{\sep}$-subgroup 
$(Q^{\sep})_{\sm}$ is parabolic in $G^{\sep}$.
Since the formation of the largest smooth
$k$-subgroup commutes with separable field
extensions (see \cite[Lemma C.4.1]{CGP}),
it follows that $Q_{\sm}$ is parabolic
in $G$. Equivalently, $Q$ contains a minimal
smooth parabolic $k$-subgroup. But this 
contradicts the above inequality.
\end{proof}

We now obtain a first class of examples
of pairs $(G, \varphi : L \to G')$ as above,
assuming that there exists a central division 
$k$-algebra $D$ of degree $p$. 
We then have an isomorphism of $k$-algebras
$D^{\otimes p} \simeq M_{p^p}(k)$, as follows e.g.~from 
\cite[Theorem 2.7.1]{Jacobson}.
By loc.~cit., Theorem 4.1.8, it follows
that there is an isomorphism of $k$-algebras
\[ D \otimes_k k^{1/p} \simeq M_p(k^{1/p}). \]
As a consequence, the $k$-group $\mathbf{GL}_1(D)$
(representing the contravariant functor 
from schemes to groups
$S \mapsto \big( D \otimes_k \cO(S) \big)^{\times}$)
satisfies $\mathbf{GL}_1(D)^{(1)} \simeq \GL_p$.
In turn, the $k$-group 
\[ \mathbf{PGL}_1(D) = \mathbf{GL}_1(D)/\bG_m \]
(where $\bG_m = \mathbf{GL}_1(k)$)
satisfies $\mathbf{PGL}_1(D)^{(1)} \simeq \PGL_p$.
This identifies the relative Frobenius homomorphism
$\mathbf{PGL}_1(D) \to \mathbf{PGL}_1(D)^{(1)}$ 
with a homomorphism
\[ \eta : \mathbf{PGL}_1(D) \longrightarrow \PGL_p. \]
Moreover, the $k$-group $\mathbf{PGL}_1(D)$ is smooth
and hence $\eta$ is a purely 
inseparable isogeny. The natural action of $G' := \PGL_p$ 
on the projective space $\bP^{p-1}$ equips it with 
the structure of a parabolic $G'$-variety 
$X' = G'/P'$, where $P'$ denotes the stabilizer
of some point of $\bP^{p-1}(k)$.

Next, consider the central simple $k$-algebra
$A := M_2(D)$, the $k$-group $\mathbf{GL}_1(A)$ 
representing the functor
$S \mapsto \big( A \otimes_k \cO(S) \big)^{\times}$,
and the quotient 
\[ G := \mathbf{PGL}_1(A) 
= \mathbf{GL}_1(A)/\bG_m. \]
Then $G$ is a simple adjoint $k$-group containing 
$\PGL_2 = \mathbf{GL} \big( M_2(k) \big)/\bG_m$.
Moreover, the standard $k$-torus
$S \subset \PGL_2$ (the image of the diagonal
$k$-torus $T \subset \GL_2$) is a maximal 
$k$-split torus of $G$. Its centralizer 
$L := C_G(S)$ is the quotient 
$C_{\mathbf{GL}_1(A)}(T)/\bG_m$, where 
$C_{\mathbf{GL}_1(A)}(T)$ is the $k$-subgroup of
$\mathbf{GL}_1(A)$ consisting of block-diagonal
$2 \times 2$ matrices with entries in 
$\mathbf{GL}_1(D)$. Thus, we have 
\[ L \simeq 
\big( \mathbf{GL}_1(D) \times \mathbf{GL}_1(D) \big)
/\bG_m , \]
where $\bG_m$ is diagonally embedded.
The first projection
$\mathbf{GL}_1(D) \times \mathbf{GL}_1(D) 
\to \mathbf{GL}_1(D)$
yields a faithfully flat homomorphism 
of $k$-groups $L \to \mathbf{PGL}_1(D)$, and hence
the desired faithfully flat homomorphism
$\varphi : L \to G'$ by composing with 
the isogeny $\eta$.

\subsection{}
\label{subsec:ex3}

In this subsection, we construct a further 
class of examples of pairs $(G,\varphi : L \to G')$ 
satisfying the same conditions as in \S \ref{subsec:ex2}: 
$G$ is an isotropic simple adjoint $k$-group, 
$L$ is the centralizer of a maximal split $k$-torus of $G$, 
and $\varphi$ is a faithfully flat homomorphism
to an isotropic simple adjoint $k$-group.

We assume that $k$ has characteristic $2$ 
and admits a quadratic space $(V',q')$, where
$\dim(V') = 2n +1$ and $q'$ is an anisotropic 
non-degenerate quadratic form. 
We first revisit an example of Borel and Tits 
(see \cite[1.2]{BTc}).
The special orthogonal group $\Or^+(V',q')$
is an adjoint semisimple $k$-group, and is anisotropic 
by \cite[4.24]{BT}. Also, the bilinear form 
associated with $q'$ has radical
$\ell$ for some line $\ell \subset V'$. Thus, the quotient 
$V'/\ell$ is a equipped with a nondegenerate alternating
form $b'$. 
This yields a purely inseparable isogeny
\[ \eta : \Or^+(V',q') \longrightarrow \Sp(V'/\ell,b'), \]
a version of the ``very special isogeny'' of 
\cite[Remark 7.1.6]{CGP}.

Next, consider the nondegenerate quadratic form 
on $V := V' \oplus k^2$
\[ 
q(v',x,y) := q'(v') + x y
\]
and its special orthogonal group $G := \Or^+(V,q)$,
an adjoint simple $k$-group. A $k$-split torus 
$S$ of $G$ is the multiplicative group $\bG_m$ 
acting on $V$ via
\[ t \cdot (v',x,y) := (v',t x,t^{-1} y). \]
The centralizer $L := C_G(S)$ is isomorphic
to $\Or^+(V',q') \times \bG_m$,
where $\Or^+(V',q')$ acts on $V'$ and
$\bG_m$ on $x,y$. As a consequence,
$S$ is a maximal $k$-split torus of $G$.
Composing the projection 
$L \to \Or^+(V',q')$ with the isogeny $\eta$
yields the desired faithfully flat homomorphism 
$\varphi : L \to \Sp(V'/\ell,b') =: G'$.

In the case where $n = 1$, this construction 
is related to that in \S \ref{subsec:ex2} 
via the isomorphism of $k$-groups
\[ \Or^+(V',q') 
\stackrel{\sim}{\longrightarrow}
\PGL_1(D), \]
where $D$ denotes the even Clifford
algebra of $q'$. Moreover, $D$ is a
quaternion algebra; in particular,
a division algebra of degree $2$ over $k$,
and every such division algebra is obtained
in this way from an anisotropic
nondegenerate quadratic form of rank $3$
(see \cite[\S 15.A, \S 26.B]{KMRT} for these
facts). But the Frobenius isogeny
$\mathbf{PGL}_1(D) \to \PGL_2$ used in 
\S \ref{subsec:ex2} differs from the
composition of the ``very special'' one,
$\Or^+(V',q') \to \Sp(V'/\ell,b') \simeq \SL_2$,
with the quotient homomorphism 
$\SL_2 \to \PGL_2$.

\bibliographystyle{amsalpha}

\end{document}